\documentclass[reqno,a4paper]{amsart} 
\usepackage[english]{babel}

\usepackage{graphicx,amsmath,amsfonts,amssymb,amsthm, mathrsfs, tikz,tikz-cd, amscd,amsbsy,amstext,mathtools, enumerate, enumitem,stmaryrd, tensor,a4wide, fancyhdr,ifthen, empheq, verbatim, latexsym, epigraph,tabularray}

\mathtoolsset{showonlyrefs,showmanualtags}

\usepackage{mathpazo}

\usepackage{nameref,zref-xr} 
\zxrsetup{toltxlabel} 
\zexternaldocument*[COHA-Yangian-]{COHA_1-dimensional_Yangian} 
\zexternaldocument*[torsion-pairs-]{Torsion_pairs_final} 

\usetikzlibrary{positioning,shapes,shadows,arrows}

\tikzstyle{roundbox} = [rectangle, draw, text centered, rounded corners, minimum height=2em]
\tikzstyle{connector} = [draw, -latex']

\usepackage{color}
\definecolor{MyDarkBlue}{rgb}{0.15,0.25,0.45}

\usepackage[utf8x]{inputenc}

\usepackage[toc,page]{appendix}
\usepackage[mathscr]{euscript} 
\allowdisplaybreaks

\newif\ifpersonal
\personaltrue 
\ifpersonal

\newcommand{\todo}[1]{\textcolor{red}{(Todo: #1)}}
\else
\newcommand\ast{\personal}[1]{\ignorespaces}
\newcommand\ast{\todo}[1]{\ignorespaces}
\fi

\newcommand{\calA}{\mathcal{A}}

\newcommand{\calC}{\mathcal{C}}

\newcommand{\calE}{\mathcal{E}}
\newcommand{\calF}{\mathcal{F}}
\newcommand{\calG}{\mathcal{G}}
\newcommand{\calH}{\mathcal{H}}
\newcommand{\calI}{\mathcal{I}}
\newcommand{\calJ}{\mathcal{J}}

\newcommand{\calQ}{\mathcal{Q}}

\newcommand{\calS}{\mathcal{S}}
\newcommand{\calT}{\mathcal{T}}
\newcommand{\calU}{\mathcal{U}}
\newcommand{\calV}{\mathcal{V}}

\newcommand{\calX}{\mathcal{X}}
\newcommand{\calY}{\mathcal{Y}}

\newcommand{\A}{\mathbb{A}}
\newcommand{\C}{\mathbb{C}}

\newcommand{\E}{\mathbb{E}}

\newcommand{\G}{\mathbb{G}}

\newcommand{\LL}{\mathbb{L}}
\newcommand{\N}{\mathbb{N}}
\newcommand{\PP}{\mathbb{P}}
\newcommand{\Q}{\mathbb{Q}}

\newcommand{\bS}{\mathbb{S}}

\newcommand{\Z}{\mathbb{Z}}

\newcommand{\scrA}{\mathscr{A}}

\newcommand{\scrC}{\mathscr{C}}

\newcommand{\scrE}{\mathscr{E}}
\newcommand{\scrF}{\mathscr{F}}

\newcommand{\scrO}{\mathscr{O}}

\newcommand{\scrT}{\mathscr{T}}
\newcommand{\scrU}{\mathscr{U}}

\newcommand{\scrW}{\mathscr{W}}
\newcommand{\scrX}{\mathscr{X}}
\newcommand{\scrY}{\mathscr{Y}}
\newcommand{\scrZ}{\mathscr{Z}}

\newcommand{\sfB}{\mathsf{B}}

\newcommand{\sfH}{\mathsf{H}}

\newcommand{\sfT}{\mathsf{T}}

\newcommand{\sfZ}{\mathsf{Z}}

\newcommand{\sft}{\mathsf{t}}

\newcommand{\frakU}{\mathfrak{U}}

\newcommand{\frakm}{\mathfrak{m}}

\newcommand{\frakq}{\mathfrak{q}}

\newcommand{\bfD}{\mathbf{D}}

\newcommand{\bfG}{\mathbf{G}}

\newcommand{\bfV}{\mathbf{V}}

\newcommand{\bff}{\mathbf{f}}

\newcommand{\bfv}{\mathbf{v}}

\newcommand{\op}{^{\mathsf{op}}}

\newcommand{\Cat}{\mathsf{Cat}}

\newcommand{\CAlg}{\mathsf{CAlg}}
\newcommand{\dAff}{\mathsf{dAff}}

\newcommand{\PrL}{\mathsf{Pr}^{\mathsf{L}}}

\newcommand{\PrLomega}{\mathsf{Pr}^{\mathsf{L},\omega}}
\newcommand{\PSh}{\mathsf{PSh}}

\newcommand{\BettiD}{\tensor*[^{\mathsf{top}}]{\bfD}{}}

\newcommand{\PreSt}{\mathsf{PreSt}}
\newcommand{\dSt}{\mathsf{dSt}}

\newcommand{\dSch}{\mathsf{dSch}}
\newcommand{\dGeom}{\mathsf{dGeom}}
\newcommand{\formaldSch}{\mathsf{fdSch}}

\newcommand{\Corr}{\mathsf{Corr}}

\newcommand{\indGeom}{\mathsf{indGeom}}
\newcommand{\dGeomqc}{\dGeom^{\mathsf{qc}}}
\newcommand{\dGeomqcqs}{\dGeom^{\qcqs}}
\newcommand{\dGeomqs}{\dGeom^{\mathsf{qs}}}
\newcommand{\indGeomindqcqs}{\indGeom^{\ind\textrm{-}\qcqs}}
\newcommand{\indGeomqcqs}{\indGeom^{\qcqs}}
\newcommand{\indGeomadm}{\indGeom^{\mathsf{adm}}}
\newcommand{\LambdadSt}{\Lambda\textrm{-}\dSt}
\newcommand{\LambdaindGeomadm}{\Lambda\textrm{-}\indGeom^{\mathsf{adm}}}

\newcommand{\proD}{\tensor*[^{\mathsf{pro}}]{\bfD}{}}

\newcommand{\catQCohnil}{\catQCoh^{\mathsf{nil}}}

\newcommand{\ind}{\mathsf{ind}}

\newcommand{\Ind}{\mathsf{Ind}}
\newcommand{\Pro}{\mathsf{Pro}}

\newcommand{\Fun}{\mathsf{Fun}}

\newcommand{\Sym}{\mathsf{Sym}}

\newcommand{\Map}{\mathsf{Map}}

\newcommand{\bfDM}{\mathbf{DM}}

\newcommand{\leqdAff}[2]{\tensor[^{\leqslant #1}]{\dAff}{_{#2}}}
\newcommand{\evccndAff}{\tensor*[^{<\infty}]{\dAff}{}}
\newcommand{\evccnPreSt}{\tensor*[^{<\infty}]{\PreSt}{}}
\newcommand{\evccnj}{\tensor*[^{<\infty}]{j}{}}
\newcommand{\conv}[1]{\tensor*[^{\mathsf{conv}}]{#1}{}}

\newcommand{\Red}[1]{\tensor*[^{\mathsf{red}}]{#1}{}}
\newcommand{\reddAff}{\Red{\dAff}}

\newcommand{\catMod}{\mathsf{Mod}}

\newcommand{\Modd}{\mathsf{Mod}}

\newcommand{\catCoh}{\mathsf{Coh}}
\newcommand{\catCohnil}{\mathsf{Coh}^{\mathsf{nil}}}

\newcommand{\catCohb}{\mathsf{Coh}^{\mathsf{b}}}
\newcommand{\catPerf}{\mathsf{Perf}}
\newcommand{\catAPerf}{\mathsf{APerf}}

\newcommand{\catQCoh}{\mathsf{QCoh}}

\newcommand{\catIndCoh}{\mathsf{IndCoh}}
\newcommand{\catCohbnil}{\mathsf{Coh}^{\mathsf{b}, \mathsf{nil}}}

\newcommand{\ps}{\mathsf{ps}}

\newcommand{\qcqs}{\mathsf{qcqs}}

\newcommand{\stackCoh}{\mathfrak{Coh}}

\newcommand{\stackCohps}{\mathfrak{Coh}_{\mathsf{ps}}}
\newcommand{\stackCohpsnil}{\mathfrak{Coh}_{\mathsf{ps}}^{\mathsf{nil}}}

\newcommand{\stackCohext}{\mathfrak{Coh}^{\mathsf{ext}}}
\newcommand{\stackCohpsext}{\mathfrak{Coh}_{\mathsf{ps}}^{\mathsf{ext}}}

\newcommand{\trunc}[1]{\mathsf{t}_0(#1)}

\newcommand{\dstackCoh}{\mathbf{Coh}}
\newcommand{\dstackCohnil}{\mathbf{Coh}^{\mathsf{nil}}}
\newcommand{\dstackCohps}{\mathbf{Coh}_{\mathsf{ps}}}
\newcommand{\dstackCohpsnil}{\mathbf{Coh}_{\mathsf{ps}}^{\mathsf{nil}}}

\newcommand{\dstackCohext}{\mathbf{Coh}^{\mathsf{ext}}}

\newcommand{\Hilb}{\mathsf{Hilb}}

\newcommand{\gr}{\mathsf{gr}}

\newcommand{\HBM}{\mathsf{H}^{\mathsf{BM}}}

\newcommand{\Hbullet}{\mathsf{H}^\bullet}
\newcommand{\HBMbullet}{\mathsf{H}_\bullet^{\mathsf{BM}}}
\newcommand{\HBMbulletT}{\mathsf{H}_\bullet^T}

\newcommand{\Hmotbullet}{\mathsf{H}_\bullet^{\mathsf{mot}}}
\newcommand{\HmotbulletT}{\mathsf{H}_\bullet^{\mathsf{mot}, T}}

\newcommand{\HBMDGamma}{\mathsf{H}^{\bfD,\Gamma}}

\newcommand{\HBMGamma}{\mathsf{H}^{\mathsf{BM},\Gamma}}

\newcommand{\CBM}{\mathsf{C}^{\mathsf{BM}}_\bullet}

\newcommand{\coha}{\mathbf{HA}}
\newcommand{\evencoha}{\coha^{\mathsf{even}}}

\newcommand{\DGammacoha}{\coha^{\bfD,\Gamma}}
\newcommand{\DGammaTcoha}{\coha^{\bfD,\Gamma, T}}
\newcommand{\motcoha}{\coha^{\mathsf{mot}}}
\newcommand{\motTcoha}{\coha^{\mathsf{mot}, T}}

\newcommand{\Exp}{\mathsf{Exp}}

\newcommand{\id}{{\mathsf{id}}}
\newcommand{\Spec}{\mathsf{Spec}}
\newcommand{\fib}{\mathsf{fib}}
\newcommand{\ch}{\mathsf{ch}}

\newcommand{\bfDelta}{\mathbf{\Delta}}

\newcommand{\Spc}{\mathsf{Spc}}

\newcommand{\ev}{\mathsf{ev}}

\newcommand{\calHom}{\mathcal{H}\mathsf{om}}

\newcommand{\End}{\mathsf{End}}

\newcommand{\Pic}{\mathsf{Pic}}

\DeclareMathOperator*{\fcolim}{``colim''}
\DeclareMathOperator*{\flim}{``lim''}

\newcommand{\bfjmathhat}{\hat{\boldsymbol{\jmath}}}

\newcommand{\overunder}{/\!\!/}
\newcommand{\red}{\mathsf{red}}

\makeatletter
\newcommand{\colim@}{%
	\vtop{\m@th\ialign{##\cr
			\hfil$\operator@font colim$\hfil\cr
			\noalign{\nointerlineskip\kern1.5\ex@}\cr
			\noalign{\nointerlineskip\kern-\ex@}\cr}}%
}
\newcommand{\colim}{%
	\mathop{\mathpalette\colim@{\textstyle}}\nmlimits@
}
\makeatother

\usepackage{thmtools}

\declaretheoremstyle[
spaceabove=1em, spacebelow=1em,
headfont=\bfseries, notefont=\normalfont, bodyfont=\itshape,
headpunct={.}, notebraces={(}{)}, postheadspace={ }
]{basic-theorem}

\declaretheoremstyle[
spaceabove=1em, spacebelow=1em,
headfont=\bfseries, notefont=\normalfont, bodyfont=\normalfont,
headpunct={.}, notebraces={(}{)}, postheadspace={ }, qed=$\oslash$
]{basic-definition}

\declaretheoremstyle[
spaceabove=1em, spacebelow=1em,
headfont=\itshape, notefont=\normalfont, bodyfont=\normalfont,
headpunct={.}, notebraces={(}{)}, postheadspace={ }, qed=$\triangle$
]{basic-remark}

\theoremstyle{basic-theorem}

\newtheorem{theorem}{Theorem}[section]
\newtheorem{corollary}[theorem]{Corollary}
\newtheorem{lemma}[theorem]{Lemma}
\newtheorem{proposition}[theorem]{Proposition}

\declaretheorem[style=basic-definition, numberlike=theorem]{definition}
\declaretheorem[style=basic-definition, numberlike=theorem]{example}
\declaretheorem[style=basic-definition, numberlike=theorem]{notation}
\declaretheorem[style=basic-definition, numberlike=theorem]{warning}
\declaretheorem[style=basic-definition, numberlike=theorem]{variant}
\declaretheorem[style=basic-definition, numberlike=theorem]{construction}
\declaretheorem[style=basic-remark, numberlike=theorem]{remark}
\declaretheorem[style=basic-remark, numberlike=theorem]{recollection}

\newtheorem{assumption}{Assumption}



\numberwithin{equation}{section}



\declaretheorem[style=basic-theorem,name=Theorem]{theoremintroduction}

\usepackage{hyperref}

\hypersetup{anchorcolor=black, linktocpage=true,
	colorlinks=true,
	citecolor=MyDarkBlue,
	linkcolor=MyDarkBlue,
	urlcolor=black,
	pdftitle={Nilpotent Cohomological Hall algebras of surfaces}, 
	breaklinks=true,
	plainpages=true
}
\usepackage[hyperpageref]{backref}

\title[Nilpotent Cohomological Hall algebras of surfaces]{Nilpotent Cohomological Hall algebras of surfaces}

\author[D.-E.~Diaconescu]{Duiliu-Emanuel Diaconescu}
\address[Duiliu-Emanuel Diaconescu]{New High Energy Theory Center - Serrin Building, Rutgers, The State University Of New Jersey, 126 Frelinghuysen Rd., Piscataway, NJ 08854-8019, USA}
\curraddr{}
\email{\href{mailto:duiliu@physics.rutgers.edu}{duiliu@physics.rutgers.edu}}

\author[M.~Porta]{Mauro Porta}
\address[Mauro Porta]{Institut de recherche mathématique avancée (IRMA), Université de Strasbourg, France and Institut Universitaire de France (IUF)}
\curraddr{}
\email{\href{mailto:porta@math.unistra.fr}{porta@math.unistra.fr}}

\author[F.~Sala]{Francesco Sala}
\address[Francesco Sala]{Università di Pisa, Dipartimento di Matematica, Largo Bruno Pontecorvo 5, 56127 Pisa (PI), Italy}
\address{Kavli IPMU (WPI), UTIAS, The University of Tokyo, Kashiwa, Chiba 277-8583, Japan}
\curraddr{}
\email{\href{mailto:francesco.sala@unipi.it}{francesco.sala@unipi.it}}

\author[O.~Schiffmann]{Olivier Schiffmann}
\address[Olivier Schiffmann]{Laboratoire de Math\'ematiques d'Orsay, Universit\'e de Paris-Sud Paris-Saclay, B\^at. 425, 91405 Orsay Cedex, France, UMR8628 (CNRS), and Simion Stoilow Institute of Mathematics, Bucharest, Romania}
\email{\href{mailto:olivier.schiffmann@universite-paris-saclay.fr}{olivier.schiffmann@universite-paris-saclay.fr}}

\author[E.~Vasserot]{Eric Vasserot}
\address[Eric Vasserot]{Université de Paris, 75013 Paris, France, UMR7586 (CNRS) and Institut Universitaire de France (IUF)}
\email{\href{mailto:eric.vasserot@imj-prg.fr}{eric.vasserot@imj-prg.fr}}

\thanks{The work of the first-named author is partially supported by NSF grant DMS-1802410. The work of the third-named author is partially supported by JSPS KAKENHI Grant Numbers JP21K03197 and JP26K06721. The work of the fourth-named author is partially supported by the PNRR Grant `\textit{Cohomological Hall algebras and smooth surfaces and applications}', CF~44/14.11.2022. The work of the fifth-named author is partially supported by ANR-18-CE40-0024 ``Categorification in topology and representation theory - Catore''. This work was also partially supported by the ``National Group for Algebraic and Geometric Structures, and their Applications'' (GNSAGA – INDAM)} 



\begin{document}

\begin{abstract}
		This paper develops a framework for systematically studying cohomological ``Hecke operators'' associated with modifications of coherent sheaves on a smooth surface $X$ along a fixed proper curve $Z \subset X$ (possibly singular and reducible), using the theory of cohomological Hall algebras.
		
		More precisely, we construct a moduli stack of coherent sheaves $\dstackCoh(\widehat{X}_Z)$ on $X$ with \textit{set-theoretic} support $Z$ and we prove that its reduced is an Artin stack locally of finite type. This provides a vast generalization of the global nilpotent cone. Subsequently, we develop the needed background to define the (motivic, $T$-equivariant) cohomological Hall algebra $\coha^{T}_{X,Z}$ of  the moduli stack of coherent sheaves on $X$ with \textit{set-theoretic} support on $Z$, in the setting of a general motivic formalism $\bfD$ in the sense of Khan \cite{Khan_Voevodsky_criterion}. The algebra $\coha^{\bfD, A}_{X,Z}$ is functorial with respect to closed immersions $Z' \subset Z$ and transformations of the motivic formalism $\bfD$, and only depends on the formal neighborhood $\widehat{X}_Z$ of $Z$ in $X$.
		
		In the companion paper \cite{DPSSV-3}, we use the nilpotent COHA $\coha^{T}_{X,Z}$ to answer a question previously raised in \cite{DPS_McKay} about the precise relationship between the COHA of a minimal resolution of a Kleinian singularity and the corresponding preprojective COHA.
\end{abstract}

\maketitle
\thispagestyle{empty}

\tableofcontents

\bigskip\section{Introduction}

In the present paper, we introduce the theory of \textit{nilpotent} cohomological Hall algebras (COHAs) associated to pairs $(X, Z)$, where $X$ is a smooth surface and $Z\subset X$ is a reduced subscheme, possibly singular and reducible. Before turning to the construction of these new COHAs -- together with the technical aspects of the construction and their applications -- we describe the framework in which they naturally fit and our motivations.

COHAs of smooth surfaces provide an example of \textit{2-dimensional} COHAs; the other primary example is given by the COHA of finite-dimensional representations of a preprojective algebra of a quiver \cite{SV_Cherednik, SV_generators}.
The latter play a central role in the theory of quantum groups, due to their relationship with \textit{Maulik–Okounkov Yangians} \cite{MO_Yangian}, as established in \cite{BD_Okounkov, SV_YangiansCOHA}. More generally, it is possible to attach a COHA to any pair $(\scrC, \tau)$, where $\scrC$ is a $\C$-linear $\infty$-category of finite type and $\tau$ is a (sufficiently nice) $t$-structure on $\scrC$, see \cite[Theorem~\ref*{torsion-pairs-thm:COHA-Coh}]{DPS_Torsion-pairs}. In the situation of the McKay correspondence, different choices of the $t$-structure produce respectively a (stacky) surface COHA or a preprojective COHA, see \cite{SSS}. Understanding their precise mutual relation was a question raised first in \cite{DPS_McKay} and that found a final answer in \cite{DPSSV-3}. Nilpotent COHAs are one of the fundamental tools needed in \cite{DPSSV-3}.

From a slightly different perspective, a central challenge in the theory of COHAs is to provide explicit presentations by generators and relations.
From this point of view, the subalgebra of $0$-dimensional operators of a surface COHA has been deeply understood thanks to the work of the fourth- and fifth-named authors together with Mellit and Minets in \cite{MMSV}, but it remains challenging to address the full COHA.
At the same time, in the setting of preprojective COHAs, the fourth- and fifth-named authors introduced \textit{nilpotent} COHAs of quivers in \cite{SV_generators} and constructed a set of generators for them, showing that the problem of constructing explicit presentations for the full algebra becomes more manageable in the nilpotent setting. The goal of this paper is to generalize the construction of nilpotent COHAs from quivers to surfaces; the results of \cite{DPSSV-3} allow to obtain an explicit presentation in the case the surface is a minimal resolution of a Kleinian singularity and $Z$ is the exceptional curve.

Another motivation for our construction concerns the definition of an algebra of \textit{cohomological Hecke operators} governing modifications of coherent sheaves on a smooth surface $X$ along a curve $Z$. Interest in such algebras dates back to the 1990s. For instance, in their work on a `surface analogue' of the geometric Langlands program, the fifth-named author, together with Ginzburg and Kapranov \cite{GKV_Langlands}, introduced an algebra of Hecke operators associated with modifications of coherent sheaves on a smooth surface $X$ along an ADE configuration of $\PP^1$’s, acting on spaces of functions on certain moduli stacks of coherent sheaves on $X$. In the case of an elliptic surface with a special fiber of affine \textit{ADE} type, this algebra was further extended to a two-parameter version of the quantum toroidal algebra of the corresponding ADE quiver. Around the same time, motivated by gauge-theoretic considerations, Nakajima \cite{Nakajima1996} speculated on the existence of an algebra of cohomological Hecke operators associated with curve modifications. More recently, thanks to developments in the theory of \textit{Nakajima operators} from the perspective of derived algebraic geometry, interest in defining such algebras has resurfaced; see, for example, \cite{LJZ_Blown-up, DPSZ}. From this viewpoint, our theory of nilpotent COHAs associated with $(X, Z)$ provides an intrinsic construction of the ``largest'' algebra of \textit{cohomological Hecke operators} governing modifications of coherent sheaves on $X$ along $Z$.

\subsection{Main results}

The main aim of the present paper is to lay the foundations for \textit{nilpotent cohomological Hall algebras}. In order to describe the main results and the challenges that have to be solved, start by fixing a smooth quasi-projective complex surface $X$ and a closed subscheme $Z \subset X$. For simplicity, in this introduction we take the additional assumption that $Z$ is proper. We are then interested in:
\begin{enumerate}\itemsep=0.2cm
	\item define the moduli stack $\stackCoh(\widehat{X}_Z)$ of coherent sheaves on $X$ \textit{set-theoretically} supported on $Z$;
	
	\item define the Borel-Moore homology of such moduli stack;
	
	\item prove that the Borel-Moore homology of such moduli stack carries a canonical cohomological Hall algebra structure.
\end{enumerate}
The first challenge comes from the set-theoretic condition: this essentially says that the moduli stack in question cannot be an algebraic stack (locally) of finite type. From the point of view of the moduli of objects, we are working with the category of coherent sheaves on the formal completion $\widehat{X}_Z$, and the difficulty is that this category is not of finite type. The second challenge is to make sense of the Borel-Moore homology of $\stackCoh(\widehat{X}_Z)$. Both these difficulties are addressed by the following:
\begin{theoremintroduction}[Theorem~\ref{thm:Cohnil_admissible} and Corollary~\ref{cor:Cohnil_as_formal_completion-2}]
	The functor of points
	\begin{align}
		\stackCoh(\widehat{X}_Z) \colon \mathsf{Aff}_k\op \longrightarrow \Spc 
	\end{align}
	that sends $S \in \mathsf{Aff}_k$ to the groupoid of $S$-flat coherent sheaves on $X \times S$ set-theoretically supported on $Z \times S$ is an indgeometric stack, that is it is a filtered colimit of Artin stacks where the transition maps are closed immersions.
	Furthermore:
	\begin{enumerate}\itemsep=0.2cm
		\item its reduced $\Red{\stackCoh(\widehat{X}_Z)}$ is a quasi-separated (but not quasi-compact) Artin stack locally of finite type (in the underived sense);
		
		\item there is a canonical morphism
		\begin{align}
			\Red{\stackCoh(\widehat{X}_Z)} \longrightarrow \stackCohps(X) 
		\end{align}
		which is a closed immersion, and $\stackCoh(\widehat{X}_Z)$ coincides with the formal completion of $\stackCohps(X)$ at this closed substack.
	\end{enumerate}
\end{theoremintroduction}

\begin{example}[Global nilpotent cone]
	Let $C$ be a smooth projective curve and set $X \coloneqq T^\ast C$ be its cotangent bundle. Take $Z \coloneqq C$ to be the curve embedded as its zero section. In this case, $\stackCohps(T^\ast C)$ is the stack of Higgs sheaves on $C$, and $\Red{\stackCohps(\widehat{T^\ast C}_C)}$ coincides with the stack of nilpotent Higgs sheaves (with reduced structure). Therefore $\stackCoh(\widehat{T^\ast C}_C)$ is the formal completion of the stack of Higgs sheaves at the closed substack given by nilpotent Higgs sheaves.
\end{example}

\begin{remark}
	In the case where $C$ is a hyperplane section of a projective $K3$ surface $S$, we expect that our construction recovers the generic fiber of the deformation of the Hitchin system of $C$ constructed in \cite{Donagi_Ein_Lazarsfeld}. However, in the general case we do not know how to define the closed substack $\Red{\stackCoh(\widehat{X}_Z)}$ of $\stackCohps(X)$ \textit{a priori}, that is, without first defining the indgeometric stack $\stackCoh(\widehat{X}_Z)$. A natural possibility, that will make the object of a future work, is to replace the base of the Hitchin system with the reduced Chow variety of Rydh \cite{Rydh}.
\end{remark}

The above theorem allows to define the Borel-Moore homology of $\stackCoh(\widehat{X}_Z)$: indeed, Borel-Moore homology only depends on the underlying reduced stack, and therefore it is sensible to set as a definition
\begin{align}
	\HBM_\ast( \stackCoh(\widehat{X}_Z) ; \Q ) \coloneqq \HBM_\ast( \Red{\stackCoh(\widehat{X}_Z)} ; \Q ) \ , 
\end{align}
where the latter makes sense thanks to the work of A.\ Khan \cite{Khan_VFC}. It is worth remarking that some extra care needs to be applied since $\Red{\stackCoh(\widehat{X}_Z)}$ is not quasi-compact; as already explained in \cite[Appendix~A]{Porta_Sala_Hall}, this leads to consider these Borel-Moore homology groups as \textit{topological} vector spaces, with a topology induced by quasi-compact open exhaustions of $\Red{\stackCoh(\widehat{X}_Z)}$.

\medskip

Nevertheless, this not sufficient to endow $\HBM_\ast( \stackCoh(\widehat{X}_Z) ; \Q )$ with a Hall multiplication. The reason is, as usual, that the canonical map
\begin{align}
	\partial_0 \times \partial_2 \colon \stackCohext( \widehat{X}_Z ) \longrightarrow \stackCoh(\widehat{X}_Z) \times \stackCoh(\widehat{X}_Z) 
\end{align}
is too singular to allow to define on the nose a refined Gysin pullback. However, it turns out that the square
\begin{align}
	\begin{tikzcd}[ampersand replacement=\&]
		\stackCohext(\widehat{X}_Z) \arrow{r} \arrow{d}{\partial_0 \times \partial_2} \& \stackCohpsext(X) \arrow{d}{\partial_0 \times \partial_2} \\
		\stackCoh(\widehat{X}_Z) \times \stackCoh(\widehat{X}_Z) \arrow{r} \& \stackCohps(X) \times \stackCohps(X) 
	\end{tikzcd}
\end{align}
is a pullback.
As explained in \cite{Porta_Sala_Hall}, the right vertical map admits a canonical derived enhancement that is derived lci. The theory developed by A.\ Khan \cite{Khan_VFC} allows therefore to define a refined Gysin pullback. One \textit{could} try to exploit this to define the Hall multiplication for $\HBM_\ast( \stackCoh(\widehat{X}_Z) )$ by hand, but the problem is proceeding in this way, the Hall multiplication would a priori depend on the ambient surface $X$.

\medskip

Since in explicit computations it is useful to reduce to `local' computations depending only on the formal completion $\widehat{X}_Z$ of $X$ along $Z$ (as done in \cite[\S\ref*{COHA-Yangian-sec:explicit-computations}]{DPSSV-3}), we take a more fundamental approach. Namely, we define a derived enhancement $\dstackCoh( \widehat{X}_Z )$ of $\stackCoh(\widehat{X}_Z)$, and more precisely we construct a full $2$-Segal derived stack $\calS_\bullet \dstackCoh(\widehat{X}_Z)$ encoding the Hall multiplication as an algebra structure in correspondences. Our construction makes more generally sense for coherent sheaves on formal schemes, and does not require a priori the embedding in an ambient variety. In particular, the Hall algebra structure encoded by $\calS_\bullet \dstackCoh(\widehat{X}_Z)$ only depends, in a tautological way, on the formal completion $\widehat{X}_Z$. We then prove:
\begin{theoremintroduction}[Theorem~\ref{thm:embedded_2_Segal}]\label{thmintro:B}
	\hfill
	\begin{enumerate}\itemsep=0.2cm
		\item As a \textit{derived} stack, $\dstackCoh(\widehat{X}_Z)$ is the formal completion of $\dstackCoh(X)$ at $\Red{\stackCoh(\widehat{X}_Z)}$.
		
		\item Both squares
		\begin{align}
			\begin{tikzcd}[ampersand replacement=\&]
				\dstackCohext(\widehat{X}_Z) \arrow{r} \arrow{d}{\partial_0 \times \partial_2} \& \dstackCohext(X) \arrow{d}{\partial_0 \times \partial_2} \\
				\dstackCoh(\widehat{X}_Z) \times \dstackCoh(\widehat{X}_Z) \arrow{r} \& \dstackCoh(X) \times \dstackCoh(X)
			\end{tikzcd} 
		\end{align}
		and
		\begin{align}
			\begin{tikzcd}[ampersand replacement=\&]
				\dstackCohext(\widehat{X}_Z) \arrow{d}{\partial_1} \arrow{r} \& \dstackCohext(X) \arrow{d}{\partial_1} \\
				\dstackCoh(\widehat{X}_Z) \arrow{r} \& \dstackCoh(X)
			\end{tikzcd} 
		\end{align}
		are pullback. Here, we refer to Formula~\eqref{eq:partial} for the definition of the maps $\partial_i$ for $i=1,2$.
	\end{enumerate}
\end{theoremintroduction}

This theorem implies that the map
\begin{align}
	\partial_0 \times \partial_2 \colon \dstackCohext(\widehat{X}_Z) \longrightarrow \dstackCoh(\widehat{X}_Z) \times \dstackCoh(\widehat{X}_Z) 
\end{align}
is representable by finitely connected (in the sense of Definition~\ref{def:admissible-rpas-connected}), quasi-compact, lci Artin derived stacks (it is, in fact, a linear stack as in \cite[Proposition~3.6]{Porta_Sala_Hall}), and that the map
\begin{align}
	\partial_1 \colon \dstackCohext(\widehat{X}_Z) \longrightarrow \dstackCoh(\widehat{X}_Z) 
\end{align}
is locally rpas (in the sense of Definition~\ref{def:admissible-rpas-connected}). Notice that both source and target of these stacks are indgeometric, so it is a priori unclear that these maps are representable by Artin stacks.

\medskip

Finally, in order to have a more streamlined treatment of the Hall product, we introduce the class of \textit{admissible indgeometric derived stacks}, which is essentially characterized by the property of having a geometric reduced stack (see Theorem~\ref{thm:checking_admissibility_on_reduced}). We then extend Khan's theory of motivic Borel-Moore homology to this class of stacks in \S\ref{sec:BM-homology}. We can summarize the main result obtained there as follows:
\begin{theoremintroduction}[Theorem~\ref{thm:functoriality_of_BM_homology_admissible}]
	There exists a lax symmetric monoidal functor
	\begin{align}
		\HBM_\ast(-; \Q) \colon \Corr^\times(\indGeomadm_k)_{\mathsf{rep.lci},\mathsf{rpas}} \longrightarrow \Pro(\catMod_\Q^\heartsuit) 
	\end{align}
	that sends $\scrX$ to the pro-object
	\begin{align}
		\HBM_\ast(\scrX;\Q) \coloneqq \flim_{\scrU \Subset \scrX} \HBM_\ast(\Red{\scrU}; \Q) \ , 
	\end{align}
	where the limit is taken over the quasi-compact admissible open substacks of $\scrX$. This assignment is covariant in locally rpas morphisms s and contravariant in morphisms that are representable by finitely connected, quasi-compact, and lci Artin derived stacks.
\end{theoremintroduction}
The proof is essentially a routine extension of the six operation formalism, up to a certain renormalization procedure that allows to have pushforward functoriality for \textit{locally} proper morphism instead of only for the proper ones.
This is done in \S\ref{sec:BM-homology}. We also prove a much more general result (in the spirit of \cite{Khan_VFC}), dealing with any oriented motivic Borel-Moore homology theory.

\medskip

Finally, we combine all the ingredients discussed so far to prove the following.
\begin{theoremintroduction}[Theorem~\ref{thm:nilpotent-COHA}]\label{thmintro:existscoha} 
	Let $X$ be a smooth surface over a field $k$ of characteristic zero and let $j \colon Z \hookrightarrow X$ be the inclusion of a closed subscheme. Assume that $X$ admits a projective compactification $\overline{X}$ that contains $Z$ as a closed subscheme. Then, there exists a unital associative topological algebra structure on
	\begin{align}
		\coha_{\widehat{X}_Z} \coloneqq \HBM_\ast( \dstackCohpsnil(\widehat{X}_Z); \Q ) 
	\end{align}
	with the multiplication $p_\ast q^!$, where the map $p$ and $q$ are those in the diagram 
	\begin{align}
		\begin{tikzcd}[ampersand replacement=\&]			
			\dstackCohpsnil(\widehat{X}_Z) \times \dstackCohpsnil(\widehat{X}_Z) \& \calS_2\dstackCohpsnil(\widehat{X}_Z) \ar[swap]{l}{q} \ar{r}{p} \& \dstackCohpsnil(\widehat{X}_Z)
		\end{tikzcd}\ .
	\end{align}
	Moreover:
	\begin{enumerate}\itemsep=0.2cm
		\item \label{item:COHA-ii} if $Z' \hookrightarrow X'$ is a second closed immersion satisfying the above assumptions, an abstract isomorphism $\widehat{X}_Z \simeq \widehat{X'}_{Z'}$ of formal schemes induces an isomorphism of topological algebras
		\begin{align}
			\coha_{\widehat{X}_Z} \simeq \coha_{\widehat{X'}_{Z'}} \ . 
		\end{align}
		
		\item if $i \colon Z' \hookrightarrow Z$ is a nested closed subscheme of $X$, then the direct image $i_\ast$ gives rise to a continuous algebra morphism
		\begin{align}
			i_\ast \colon \coha_{\widehat{X}_{Z'}} \longrightarrow \coha_{\widehat{X}_Z} \ . 
		\end{align}
	\end{enumerate}
\end{theoremintroduction}

Point~\eqref{item:COHA-ii} of the above theorem plays a crucial role in the explicit computations done in \cite[\S\ref*{COHA-Yangian-sec:explicit-computations}]{DPSSV-3}, as it allows to reduce global computations to `local' ones. Its proof follows from the fact that \textit{by construction} $\dstackCohpsnil(\widehat{X}_Z)$ only depends on $\widehat{X}_Z$ as a formal scheme, and not on the full pair $(X,Z)$.

If there is an algebraic torus $T$ acting on $X$ such that $Z$ is $T$-invariant, the above theorem extends verbatim to the equivariant setting. 

We state and prove Theorem~\ref{thmintro:existscoha} in the more general context of an arbitrary motivic formalism $\bfD$, and the construction of $\coha^{\bfD}_{X,Z}$ is functorial in $\bfD$: any transformation $\bfD \to \bfD'$ induces a continuous algebra morphism $\coha^{\bfD}_{X,{Z}} \to \coha^{\bfD'}_{X,Z}$. When $\bfD$ is taken to be sheaves on the $\C$-analytification, one recovers the Borel-Moore homology discussed above, but this framework allows to deal simultaneously with motivic Chow groups and $G$-theory. As in \cite{Porta_Sala_Hall}, our approach hinges on the notion of \textit{2-Segal structures} considered by Dyckerhoff and Kapranov \cite{Dyckerhoff_Kapranov_Higher_Segal}.

\begin{remark}
	Let $Z \hookrightarrow X$ be a closed immersion as in Theorem~\ref{thmintro:existscoha}.
	Assume that $Z$ is a smooth projective curve and that the formal completion $\widehat{X}_Z$ coincide with the formal completion $\widehat{\sfT^\ast Z}_Z$; note that the latter is partially understood thanks \cite[Theorem 5.1]{Sala_Schiffmann}.
	In the general case it is possible to consider the deformation to the normal bundle of $Z$ inside $X$.
	This should give rise to a deformation of the nilpotent COHA of the pair $(X,Z)$ to the nilpotent COHA of the pair $(\mathsf N_{Z/X}, X)$.
	We will address this in the future (and we expect that the results of \cite{Zhao_Deformation_normal_cone} will be useful in this direction).
\end{remark}

Finally, in \S\ref{sec:0-dimensional}, we describe the nilpotent cohomological Hall algebra of $0$-dimensional sheaves under a very mild assumption, which is satisfied in the cases of primary interest for us in \cite{DPSSV-2, DPSSV-3}: namely, the minimal resolutions of Kleinian singularities and elliptic surfaces with singular fibers of affine type \textit{ADE}. In both situations, we are able to reduce our analysis to the well-understood case of the usual COHAs of $0$-dimensional sheaves, following \cite{MMSV}. 

\subsection*{Structure of the paper}

In \S\ref{sec:admissible_indgeometric}, we introduce the notion of \textit{admissible} indgeometric stacks. For this class of indgeometric stacks, we define in \S\ref{sec:BM-homology} their (motivic, equivariant) Borel–Moore homology, adapting to our setting the framework developed in \cite[\S\ref*{torsion-pairs-sec:homological-invariants}]{DPS_Torsion-pairs}.

In \S\ref{sec:nilpotent_coherent_sheaves}, we introduce the main moduli stack considered in this paper: the moduli stack of nilpotent coherent sheaves on the formal completion of a smooth surface along a closed subscheme. Under suitable assumptions, we associate to this stack a cohomological Hall algebra in \S\ref{sec:COHA}. Finally, in \S\ref{sec:0-dimensional}, we describe the nilpotent cohomological Hall algebra of $0$-dimensional sheaves.

The paper also contains three appendices. Appendix~\ref{appendix:ind-objects} collects general results on ind-objects used in \S\ref{sec:admissible_indgeometric}, while a theory of \textit{relative} indcoherent sheaves à la Gaitsgory is developed in Appendix~\ref{appendix:relative_ind_coherent}. Finally, Appendix~\ref{appendix:set-theoreticity-pure} provides a characterization of coherent sheaves on smooth surfaces that are set-theoretically supported on a closed subscheme; this is used in \S\ref{sec:nilpotent_coherent_sheaves}. 

\subsection*{Notation}

We set $\N\coloneqq \Z_{\geqslant 0}$.

We denote by $\Spc$ the $\infty$-category of spaces and by $\PSh$ the $\infty$-category of presheaves with coefficients in $\Spc$.

We fix a base commutative noetherian ring $k$. We denote by $\mathsf{Poly}_k$ the category of polynomial $k$-algebras in finitely many variables and we write $\mathsf{dCAlg}_k$ for the sifted completion of $\mathsf{Poly}_k$. We refer to $\mathsf{dCAlg}_k$ as the $\infty$-category of \textit{derived commutative $k$-algebras}. We also set $\dAff_k \coloneqq \mathsf{dCAlg}_k\op$, and we refer to it as the $\infty$-category of affine derived schemes. We set
\begin{align}
	\PreSt_k \coloneqq \PSh(\dAff_k)\ .
\end{align}
We shall use the attribute \textit{geometric} as a synonym of \textit{algebraic} or \textit{Artin}. 

We refer to the notation introduced in \cite[\S1.6]{Porta_Sala_Hall} for the stable $\infty$-category of quasi-coherent sheaves on derived stacks and its subcategories. Furthermore, we use the \textit{implicitly derived convention}: given a morphism of derived stacks $f \colon X \to Y$, we let 
\begin{align}
	f^\ast \colon \catQCoh(Y) \to \catQCoh(X)
\end{align}
be the \textit{derived} pullback functor, and we let 
\begin{align}
	f_\ast \colon \catQCoh(X) \to \catQCoh(Y)
\end{align}
be the \textit{derived} pushforward. Similarly, all fiber products, Hom sheaves and spaces, and tensor products will be understood in the derived sense, unless otherwise stated.

\medskip

For a smooth projective complex surface $X$, let $K_0(X)$ be its \textit{Grothendieck group} and let $N(X)$ be its \textit{numerical Grothendieck group}, where the latter is defined by:
\begin{align}
	N(X)\coloneqq K_0(X)/\equiv \ ,
\end{align}
where $F_1 \equiv F_2$ if $\ch(F_1) = \ch(F_2)$ for $F_1, F_2 \in K_0(X)$. Then, $N(X)$ is a finitely generated free abelian group. In addition, we denote by $N_1(X)$ the subgroup of numerical equivalence classes of divisors on $X$.

We denote by $\mathsf{NS}(X)$ the \textit{Neron-Severi group} of $X$. For a coherent sheaf $E$ on $X$ whose support has dimension less than or equal to one, we denote by $\ell(E)\in \mathsf{NS}(X)$ the fundamental one cycle of $E$.

Let $\catCoh_{\leqslant 1}(X) \subset \catCoh(X)$ be the subcategory of sheaves $\calE$ with $\dim \mathsf{Supp}(\calE) \leq 1$. We define the subgroup $N_{\leqslant 1}(X) \subset N(X)$ to be
\begin{align}
	N_{\leqslant 1}(X)\coloneqq \mathsf{Im}(K_0(\catCoh_{\leqslant 1}(X)) \longrightarrow N(X) ) \ .
\end{align}
Note that we have an isomorphism $N_{\leqslant 1}(X) \simeq \mathsf{NS}(X)\oplus \Z$ sending $E$ to the pair $(\ell(E), \chi(E))$. We shall identify an element $v\in N_{\leqslant 1}(X)$ with $(\beta, n) \in \mathsf{NS}(X) \oplus \Z$ by the above isomorphism.

\subsection*{Acknowledgments}

The first-named author would like to thank Yan Soibelman for enlightening discussions on cohomological Hall algebras.

The second-named author would like to thank Giuseppe Ancona, Federico Binda, Fréderic Déglise, Niels Feld, Drago\v{s} Fr\u{a}til\u{a}, Lie Fu, Benjamin Hennion, Marc Hoyois, Fangzhou Jin, Adeel Khan, Massimo Pippi, Marco Robalo, and Tony Yue Yu for many discussions on motivic Borel-Moore homology and their extension to admissible indgeometric stacks.

The third-named author would like to thank Andrei Neguţ for very helpful discussions, which took place under the MIT-UNIPI Project (XI call). Moreover, he acknowledges the MIUR Excellence Department Project awarded to the Department of Mathematics, University of Pisa, CUP I57G22000700001. Finally, he is a member of GNSAGA of INDAM.

A part of the paper was finalized during a research visit of the second-named and fourth-named authors at the Department of Mathematics of the University of Pisa under GNSAGA – INDAM research visits program and the 2023 Visiting Fellow program of the University of Pisa, respectively. Another part of the paper was finalized during research visits of the third-named author at Kavli IPMU, the University of Tokyo, in 2023 and 2024. Finally, the fourth-named author would like to thank the Simion Stoilow Institute of Mathematics for the great working conditions provided during the preparation of this paper.

\bigskip\section{Admissible indgeometric stacks}\label{sec:admissible_indgeometric}

In this section, we introduce the class of \textit{admissible indgeometric stacks}, that are essentially characterized by the property that their reduced stack is geometric, see Theorem~\ref{thm:checking_admissibility_on_reduced}. We study their main properties, in particular establishing the existence of canonical ind-presentations, see Theorem~\ref{thm:indization}. This foundational work will be later needed to define Borel-Moore homology for these stacks, together with its natural functorialities, in \S\ref{sec:BM-homology}. The main source of examples for us would be the derived stack of coherent sheaves on a quasi-projective scheme $X$ set-theoretically supported on a closed subscheme $Z$, that will be discussed in detail in \S\ref{sec:nilpotent_coherent_sheaves}.

\subsection{Generalities on derived stacks}

We review, following \cite[Chapter 2]{Gaitsgory_Rozenblyum_Study_I} and \cite[\S17.4]{Lurie_SAG}, the main deformation-theoretic properties of derived stacks needed throughout the paper. A similar review can be found in \cite[Appendix B]{Brantner_Magidson_Nuiten}.

\begin{notation}
	\hfill
	\begin{itemize}\itemsep=0.2cm
		\item Fix an integer $n \geqslant 0$. We denote by $\leqdAff{n}{k}$ be the full subcategory spanned by $n$-coconnective affine derived schemes, i.e., the objects of the form $\Spec(A)$ where $A$ is an animated commutative ring satisfying $\pi_i(A) \simeq 0 \quad \text{ for } i > n$. We formally allow the case $n = \infty$, which we refer to as the \textit{eventually coconnective} subcategory of $\dAff_k$, denoted
		\begin{align}
			\evccndAff_k \coloneqq \bigcup_{n \geqslant 0} \leqdAff{n}{k} \ . 
		\end{align}
		We let $\evccnj \colon \evccndAff_k \hookrightarrow \dAff_k$ be the natural inclusion.
		
		\item For $S = \Spec(A) \in \dAff_k$ and $\calF \in \catQCoh(S)_{\geqslant 0}$. We let $A \oplus \calF$ be the \textit{split square-zero extension of $A$ by the connective $A$-module $\calF$}, and we set $S[\calF] \coloneqq \Spec(A \oplus \calF) \in \dAff_k$.
		For $S \in \dAff_k$ and $\calF \in \catQCoh(S)_{\geqslant 0}$, a \textit{derivation of $S$ in $\calF$ is a retraction of the canonical inclusion $S \to S[\calF]$}. Given a derivation $d$, we denote the \textit{square-zero extension associated to $d$} the affine derived scheme $S_d[\calF]$ given as the pushout in $\dAff_k$
		\begin{align}
			\begin{tikzcd}[ampersand replacement=\&]
				S[\calF] \arrow{r}{d} \arrow{d}{d_0} \& S \arrow{d} \\
				S \arrow{r} \& S_d[\calF] 
			\end{tikzcd} \ ,
		\end{align}
		where $d_0$ is the zero derivation.
		
		\item Given $S = \Spec(A)$ we set $\Red{S} \coloneqq \Spec(\pi_0(A)_{\mathsf{red}})$. We say that $S$ is reduced if the canonical map $\Red{S} \to S$ is an equivalence, and we denote by $\reddAff_k$ the full subcategory of $\dAff_k$ spanned by reduced derived schemes.
		We let $\Red{j} \colon \reddAff_k \longrightarrow \dAff_k$ be the canonical inclusion. Given $F \in \dSt$, we set
		\begin{align}
			\Red{F} \coloneqq \Red{j}_! \Red{j}^\ast(F) \ , 
		\end{align}
		where $\Red{j}_!$ denotes the left Kan extension along $\Red{j}$. We refer to $\Red{F}$ as the \textit{reduced stack of $F$} (see \cite[\S2.1]{CPTVV}). \qedhere
	\end{itemize}
\end{notation}

\begin{definition}\label{def:laft-convergent}
	Let $f \colon F \to G$ be a morphism of derived prestacks. We say that $f$ is:
	\begin{enumerate}\itemsep=0.2cm
		\item \textit{locally almost of finite type} (laft) if for every integer $n \geqslant 0$ and every cofiltered diagram $S_\bullet \colon I \to \tensor*[^{\leqslant n}]{\dAff}{_k}$ with limit $S$, the square
		\begin{align}
			\begin{tikzcd}[ampersand replacement=\&,cells={font=\everymath\expandafter{\the\everymath\displaystyle}}]
				\colim_{i \in I\op} F(S_i) \arrow{r} \arrow{d} \& F(S) \arrow{d} \\
				\colim_{i \in I\op} G(S_i) \arrow{r} \& G(S)
			\end{tikzcd}
		\end{align}
		is a pullback square;
		
		\item \textit{convergent} if for every affine derived scheme $S \in \dAff_k$ the square
		\begin{align}
			\begin{tikzcd}[ampersand replacement=\&,cells={font=\everymath\expandafter{\the\everymath\displaystyle}}]
				F(S) \arrow{r} \arrow{d} \& \lim_{n \geqslant 0} F(\sft_{\leqslant n}(S)) \arrow{d} \\
				G(S) \arrow{r} \& \lim_{n \geqslant 0} G(\sft_{\leqslant n}(S))
			\end{tikzcd}
		\end{align}
		is a pullback.
		
		\item \textit{infinitesimally cohesive} if for every $S \in \dAff_k$, every $\calF \in \catQCoh(S)_{\geqslant 1}$ and every derivation $d \colon S[\calF] \to S$, the square
		\begin{align}
			\begin{tikzcd}[ampersand replacement=\&,cells={font=\everymath\expandafter{\the\everymath\displaystyle}}]
				F(S_d[\calF]) \arrow{r} \arrow{d} \& F(S) \times_{F(S[\calF])} F(S) \arrow{d} \\
				G(S_d[\calF]) \arrow{r} \& G(S) \times_{G(S[\calF])} G(S)
			\end{tikzcd}
		\end{align}
		is a pullback;
		
		\item a \textit{nil-equivalence} if the induced morphism $\Red{f} \colon \Red{X} \to \Red{Y}$ is an equivalence.
	\end{enumerate}
	When $G = \Spec(k)$, we say that $F$ is \textit{laft} (resp.\ \textit{convergent}, \textit{infinitesimally cohesive}).
\end{definition}

\begin{notation}
	We denote by $\conv{\PreSt}_k$ the full subcategory of $\PreSt_k$ of convergent prestacks.
\end{notation}

\begin{recollection}[Convergence]\label{recollection:convergence}
	Set $\evccnPreSt_k \coloneqq \PSh(\evccndAff_k)$. Left Kan extension along $\evccnj$ allows to see $\evccnPreSt_k$ as a full subcategory of $\PreSt_k$. Then a derived prestack $F$ is convergent if and only if the canonical morphism
	\begin{align}
		F \longrightarrow \evccnj_\ast \evccnj^\ast(F) 
	\end{align}
	is an equivalence, where $\evccnj_\ast$ denotes the right Kan extension along $\evccnj$. We set
	\begin{align}
		\conv{(-)} \coloneqq \evccnj_\ast \evccnj^\ast \colon \PreSt_k \longrightarrow \PreSt_k \ , 
	\end{align}
	and we refer to it as the \textit{convergent completion functor}. It provides a left adjoint for the inclusion $\evccnPreSt_k \hookrightarrow \PreSt_k$. Concretely, for every derived prestack $F$ and every affine derived scheme $S$, one has
	\begin{align}
		\conv{F}(S) \coloneqq \lim_{n \geqslant 0} F(\sft_{\leqslant n}(S)) \ . 
	\end{align}
	In particular, a derived prestack and its convergent completion always share the same classical truncation. Notice that by construction $\conv{(-)}$ commutes with arbitrary limits in $\PreSt_k$.
\end{recollection}

Finally, we record the following elementary but useful. Recall that a map of derived prestack $F \to G$ is \textit{formally étale} if it admits $0$ as cotangent complex (see \cite[Definition~17.2.4.2]{Lurie_SAG}).
\begin{lemma}\label{lem:testing_equivalences_on_reduced}
	Let $f \colon F \to G$ be a formally étale morphism between laft derived prestacks. Assume that $f$ is convergent and infinitesimally cohesive. Then, $f$ is an equivalence if and only if it is a nil-equivalence.
\end{lemma}

\begin{proof}
	If $f$ is an equivalence, the same goes for $\Red{f}$. As for the converse, it suffices to show that for every affine derived scheme $S = \Spec(A)$, the map $f$ induces an equivalence on $S$-points. Since both $F$ and $G$ are laft, it is enough to treat the case where $S$ is almost of finite presentation over $k$. Proceeding by induction on the Postnikov tower, we further reduce to treat the case where $S$ is underived. In this case, the map $A = \pi_0(A) \to \pi_0(A)_{\red}$ can be factored as a finite composition of square-zero extensions. Proceeding by induction on the length of this presentation, the conclusion follows combining convergence, infinitesimal cohesiveness and \cite[Theorem~7.4.1.23]{Lurie_HA}.
	\end{proof}

\subsection{Indgeometric stacks and their basic properties}\label{subsec:indgeometricity}

Mimicking \cite[\S2.1]{Gaitsgory_Rozenblyum_Study_II}, we set:
\begin{definition}\label{def:indgeometric_stack}
	A derived prestack $F$ is said to be an \textit{$n$-indgeometric derived stack} if it satisfies:
	\begin{enumerate}\itemsep=0.2cm
		\item \label{item:indgeometric_stack-1} $F$ is convergent;
		
		\item  \label{item:indgeometric_stack-2} there exists a filtered diagram $X_\bullet \colon I \to \PreSt_k$ together with a map
		\begin{align}
			\phi \colon X \coloneqq \colim_{\alpha \in I} X_\alpha \longrightarrow F 
		\end{align}
		such that:
		\begin{enumerate}\itemsep=0.2cm
			\item \label{item:indgeometric_stack-2-a} for every $\alpha \in I$, $X_\alpha$ is an $n$-geometric derived stack;
			
			\item \label{item:indgeometric_stack-2-b} for every morphism $\alpha \to \beta$ in $I$, the transition map $X_\alpha \to X_\beta$ is a closed immersion almost of finite type;
			
			\item \label{item:indgeometric_stack-2-c} the restriction
			\begin{align}
				\evccnj^\ast \colon X \vert_{\evccndAff_k\op} \longrightarrow F \vert_{\evccndAff_k\op} 
			\end{align}
			is an equivalence in $\evccnPreSt_k$.
		\end{enumerate}
		We refer to the diagram $\{X_\alpha\}_{\alpha \in I}$ as a \textit{presentation for $F$}.
	\end{enumerate}
	We say that $F$ is an \textit{indgeometric derived stack} if it is an $n$-indgeometric derived stack for some $n \geqslant 0$, and we denote by $\indGeom_k$ the full subcategory of $\PreSt_k$ spanned by indgeometric stacks.
		\end{definition}

\begin{warning}
	Note that in \cite[Definition~4.1]{CW-Ind-Geometric}, the authors introduce a notion of ind-geometric stacks, which is more general than ours. In particular, the definition of \textit{loc.\ cit.} is the natural generalization to stacks of the definition of ind-schemes in \cite{Gaitsgory_Rozenblyum_DG_indschemes}. Some of the results proved in this section overlap with similar ones proved in \cite[\S4.2 \& \S4.3]{CW-Ind-Geometric}.
\end{warning}

\begin{remark}\label{rem:presentation_in_prestacks}
	Let $\{X_\alpha\}_{\alpha \in I}$ be a presentation for a indgeometric stack $F$ and let
	\begin{align}
		X \coloneqq \colim_{\alpha \in I} X_\alpha \ , 
	\end{align}
	the colimit being computed in $\PreSt_k$. Then, the structural morphism $\phi \colon X \to F$ is not necessarily an equivalence, but it induces a canonical equivalence $\overline{\phi} \colon \conv{X} \stackrel{\sim}{\to} F$.
\end{remark}

\begin{warning}\label{rem:indization}
	Notice that in the above definition the presentation $\{X_\alpha\}_{\alpha \in I}$ is not part of the structure defining $F$. Notice also that, despite the terminology, $F$ is not defined as an ind-object in some category. Of course, if we \textit{fix} some presentation $\{X_\alpha\}_{\alpha \in I}$ for $F$, then we can define an ind-object
	\begin{align}
		\fcolim_{\alpha \in I} X_\alpha \in \Ind(\dGeom_k) \ , 
	\end{align}
	which, a priori, \textit{depends on the presentation} of $F$. Our first step is to show that under some mild conditions, it is possible to attach to $F$ a \textit{canonical} ind-object.
\end{warning}

Before starting working towards the goal sketched in Warning~\ref{rem:indization}, let us collect a couple of basic properties of indgeometric derived stacks. We start by the following simple lemma, asserting that asking for convergence is a harmless assumption (compare with Remark~\ref{rem:presentation_in_prestacks}):
\begin{lemma}\label{lem:convergent_completion_produces_indstack}
	Let $X_\bullet \colon I \to \PreSt_k$ be a filtered diagram and assume that:
	\begin{enumerate}[label=(\alph*)] \itemsep0.2cm
		\item there exists an integer $n$ such that for every $i \in I$, $X_i$ is an $n$-geometric derived stack;
		
		\item for every morphism $i \to j$ in $I$, the morphism $X_i \to X_j$ is a closed immersion almost of finite type.
	\end{enumerate}
	Set
	\begin{align}
		X \coloneqq \colim_{i \in I} X_i \ , 
	\end{align}
	the colimit being computed in $\PreSt_k$. Then $\conv{X}$ is an $n$-indgeometric derived stack.
\end{lemma}

\begin{proof}
	Set $F \coloneqq \conv{X}$. Then $F$ is convergent by assumption, and since $\conv{(-)}$ is a left adjoint to the inclusion $\conv{\PreSt}_k \hookrightarrow \PreSt_k$, there is a canonical map $X \to F$. Since $\evccnj \colon \evccndAff_k \hookrightarrow \dAff_k$ is fully faithful, the counit of $\evccnj^\ast \dashv \evccnj_\ast$ is an equivalence, whence the conclusion.
\end{proof}

Another easy but handy fact is that one can ask for a more restrictive condition on the presentation of an indgeometric derived stack at no cost. Indeed, we have the following.
\begin{lemma}\label{lem:indstack_eventually_coconnective_presentation}
	For a convergent derived stack $F \colon \dAff_k\op \to \Spc$, the following conditions are equivalent:
	\begin{enumerate}\itemsep=0.2cm
		\item \label{item:indstack_eventually_coconnective_presentation-1} there exists a presentation $\{X_\alpha\}_{\alpha \in I}$ for $F$ satisfying Definition~\ref{def:indgeometric_stack}--\eqref{item:indgeometric_stack-2};
		
		\item \label{item:indstack_eventually_coconnective_presentation-2} there exists a presentation $\{X_\alpha\}_{\alpha \in I}$ for $F$ satisfying Definition~\ref{def:indgeometric_stack}--\eqref{item:indgeometric_stack-2} and moreover for every $\alpha \in I$, $X_\alpha \in \evccnPreSt_k$.
	\end{enumerate}
\end{lemma}

\begin{proof}
	Obviously \eqref{item:indstack_eventually_coconnective_presentation-2} implies \eqref{item:indstack_eventually_coconnective_presentation-1}.
	
	For the converse, fix a presentation $\{X_\alpha\}_{\alpha \in I}$ for $F$ satisfying Definition~\ref{def:indgeometric_stack}--\eqref{item:indgeometric_stack-2}. Considering $\N$ as a poset in the natural way, we see that $I \times \N$ is obviously filtered. We define a functor
	\begin{align}
		Y_{\bullet,\ast} \colon  \longrightarrow \PreSt_k 
	\end{align}
	by the rule $Y_{\alpha,n} \coloneqq \sft_{\leqslant n}(X_\alpha)$. There are canonical maps
	\begin{align}
		\colim_{n \in \N} \sft_{\leqslant n}(X_\alpha) \longrightarrow X_\alpha \ , 
	\end{align}
	which induce a canonical morphism
	\begin{align}
		\colim_{(\alpha,n) \in I \times \N} Y_{\alpha,n} \longrightarrow F \ . 
	\end{align}
	Conditions \eqref{item:indgeometric_stack-2-a} and \eqref{item:indgeometric_stack-2-b} in Definition~\ref{def:indgeometric_stack} are trivially satisfied.
	As for condition \eqref{item:indgeometric_stack-2-c}, our assumption is that for every $S \in \evccndAff_k$ the canonical morphism
	\begin{align}
		\colim_{\alpha \in I} \Map_{\PreSt_k}(S, X_\alpha) \longrightarrow \Map_{\PreSt_k}(S, F) 
	\end{align}
	is an equivalence. Thus, to conclude it is enough to argue that the canonical morphism
	\begin{align}
		\colim_{n \in \N} \Map_{\PreSt_k}( S, \sft_{\leqslant n}(X_\alpha) ) \longrightarrow \Map_{\PreSt_k}(S, X_\alpha) 
	\end{align}
	is an equivalence. Choose an integer $m$ such that $S \in \tensor[^{\leqslant m}]{\dAff}{_k}$. Then, we have
	\begin{align}
		\colim_{n \in \N} \Map_{\PreSt_k}( S, \sft_{\leqslant n}(X_\alpha) ) & \simeq \colim_{n \in \N} \Map_{\PreSt_k}( S, \sft_{\leqslant m}(\sft_{\leqslant n}(X_\alpha)) ) \simeq \Map_{\PreSt_k}(S, \sft_{\leqslant m}(X_\alpha)) \\
		& \simeq \Map_{\PreSt_k}(S, X_\alpha) \ ,
	\end{align}
	whence the conclusion.
\end{proof}

Next, let us observe that indgeometric stacks automatically satisfy flat hyperdescent, and hence define objects in $\dSt_k$. To see this, we need the following mild generalization of \cite[Chapter~2, Proposition~1.2.2]{Gaitsgory_Rozenblyum_Study_II}.
\begin{lemma}\label{lem:truncatedness_mapping_spaces}
	Let $F$ be an $n$-indgeometric stack. Then for every $m$-truncated affine derived scheme $S \in \tensor*[^{\leqslant m}]{\dAff}{_k}$, $\Map(S,F)$ is $(n+m)$-truncated. 
\end{lemma}

\begin{proof}
	Since filtered colimits in $\Spc$ commute with taking homotopy groups and since colimits in $\PreSt_k$ are computed objectwise, we can replace $F$ by an $n$-geometric derived stack. In this case, the statement is well known, and follows for instance from \cite[Lemma~2.1.1.2]{TV_HAG-II}.
\end{proof}

\begin{proposition}\label{prop:indgeometric_descent}
	Every indgeometric stack satisfies flat hyperdescent, and in particular defines an object in $\dSt_k$.
\end{proposition}

\begin{proof}
	Let $S_\bullet \to S$ be a flat hypercover in $\dAff$ and let $F$ be a indgeometric stack. Consider the commutative square
	\begin{align}
		\begin{tikzcd}[ampersand replacement=\&]
			\Map(S,F) \arrow{r} \arrow{d} \& \lim_{[n] \in \bfDelta} \Map(S_n, F) \arrow{d} \\
			\lim_{m \geqslant 0} \Map(\sft_{\leqslant m}(S), F) \arrow{r} \& \lim_{[n] \in \bfDelta} \Map(\sft_{\leqslant m}(S_n), F)
		\end{tikzcd}  \ .
	\end{align}
	We have to prove that the top horizontal map is an equivalence. Since $F$ is convergent, the vertical maps are equivalences. Therefore, it is enough to prove the same statement assuming in addition that $S$ (and hence each $S_n$) is $m$-truncated for some $m \geqslant 0$.
	
	Choose a presentation $\{X_\alpha\}_\alpha$ for $F$ by $\ell$-geometric stacks. Then, we are called to check that the canonical map
	\begin{align}
		\colim_\alpha \Map(S, X_\alpha) \longrightarrow \lim_{[n] \in \bfDelta} \colim_\alpha \Map(S_n, X_\alpha) 
	\end{align}
	is an equivalence. However, the previous lemma guarantees that $\Map(S_n, X_\alpha)$ is $(m + \ell)$-truncated for every $[n] \in \bfDelta$ and every integer $\ell \geqslant 0$. Therefore, the right hand side is the limit of a diagram with values in $\Spc_{\leqslant m + \ell}$, and therefore the canonical map
	\begin{align}
		\lim_{[n] \in \bfDelta} \colim_\alpha \Map(S_n, X_\alpha) \longrightarrow \lim_{[n] \in \bfDelta_{\leqslant m + \ell + 2}} \colim_\alpha \Map(S_n, X_\alpha) 
	\end{align}
	is an equivalence. At this point, the conclusion follows from the fact that filtered colimits commute with finite limits.
\end{proof}

\subsection{Admissibility}\label{subsec:admissible}

For our applications, we only need a special class of indgeometric stacks, that we call \textit{admissible indgeometric stacks}. They are akin to the notion of inf-scheme (and not of ind-inf-scheme) introduced in \cite[Definition~2.3.1.2]{Gaitsgory_Rozenblyum_Study_II}.

\begin{definition}\label{def:admissible_indgeometric_stack}
	Let $F$ be an indgeometric derived stack. We say that:
	\begin{enumerate}\itemsep=0.2cm
		\item $F$ is \textit{ind-qcqs} if it admits a presentation $\{X_\alpha\}$ where each $X_\alpha$ is a quasi-compact and quasi-separated geometric derived stack;
		
		\item $F$ is \textit{qcqs} if it is ind-qcqs with a presentation $\{X_\alpha\}_{\alpha \in I}$ where all the transition maps are nil-equivalences;
		
		\item $F$ is \textit{admissible} if there exists a (possibly transfinite) sequence
		\begin{align}
			\emptyset = U_0 \hookrightarrow U_1 \hookrightarrow \cdots U_\alpha \hookrightarrow U_{\alpha + 1} \hookrightarrow \cdots 
		\end{align}
		of open Zariski immersions between qcqs indgeometric stacks, whose colimit in $\PreSt_k$ is $F$. We refer to the collection $\{U_\alpha\}$ as an \textit{admissible open exhaustion of $F$}.
	\end{enumerate}
	We let $\indGeomindqcqs_k$, $\indGeomqcqs_k$, and $\indGeomadm_k$ denote the full subcategories of $\indGeom_k$ spanned by ind-qcqs, qcqs, and admissible objects, respectively.
\end{definition}

\begin{remark}\label{rem:cofinal_presentations}
	Let $F$ be an indgeometric derived stack and let $\{X_\alpha\}_{\alpha \in I}$ be a presentation. If $f \colon J \to I$ is a cofinal map, then the restricted diagram $\{X_{f(\beta)}\}_{\beta \in J}$ is again a presentation for $F$. In particular, $F$ is qcqs if and only if there exists a cofinal map $f \colon J \to I$ such that all the transition maps in $\{X_{f(\beta)}\}_{\beta \in J}$ are nil-equivalences. Notice that since $I$ is filtered, for every $\alpha_0 \in I$ the forgetful map $I_{\alpha_0 /} \to I$ is cofinal. In other words, $F$ is qcqs if and only if it admits a presentation where the transition maps become nil-equivalences for ``sufficiently large indexes''.
\end{remark}

A very important fact concerning the notion of admissibility is that it can be tested at the reduced level, in particular completely disregarding the derived structure. To prove this fact, we need the following two technical results.
\begin{lemma}\label{lem:compactness}
	\hfill
	\begin{enumerate}\itemsep=0.2cm
		\item \label{item:compactness-1} Let $F$ be an $n$-indgeometric stack and let $S \in \tensor*[^{<\infty}]{\dGeom}{^{\qcqs}_k}$ be an eventually coconnective qcqs geometric derived stack. Then, for every presentation $\{X_\alpha\}_\alpha$ of $F$, the canonical map
		\begin{align}
			\colim_\alpha \Map_{\dSt}(S,X_\alpha) \longrightarrow \Map_{\dSt}(S,F) 
		\end{align}
		is an equivalence.
		
		\item \label{item:compactness-2} Let $S \in \indGeomqcqs_k$ and let $F \in \indGeomadm_k$ be an admissible indgeometric derived stack. Then, for any quasi-compact open exhaustion $\{U_\alpha\}$ of $F$, the canonical map
		\begin{align}\label{eq:qc_exhaustion}
			\colim_\alpha \Map_{\dSt_k}(S, U_\alpha) \longrightarrow \Map_{\dSt_k}(S, F)
		\end{align}
		is an equivalence.
	\end{enumerate}
\end{lemma}

\begin{proof}
	We first prove \eqref{item:compactness-1}. Since $S$ is geometric and qcqs, we can find a smooth hypercover $S_\bullet$ of $S$ with the property that for every $[a] \in \bfDelta$, $S_a$ is a \textit{finite} disjoint union of derived affines. Furthermore, if $S \in \tensor*[^{<m}]{\dGeom}{^{\mathsf{qc}}_k}$, then we can equally assume that $S_a \in \tensor*[^{<m}]{\dAff}{_k}$. At this point, consider the following commutative diagram:
	\begin{align}
		\begin{tikzcd}[ampersand replacement=\&]
			\Map_{\dSt_k}(S,F) \arrow{r} \arrow{d} \& \colim_\alpha \Map_{\dSt_k}(S, X_\alpha)  \arrow{d} \\
			\lim_{[a] \in \bfDelta} \Map_{\dSt_k}(S_a, F) \arrow{r} \& \colim_\alpha \lim_{[a] \in \bfDelta} \Map_{\dSt_k}(S_a, X_\alpha) 
		\end{tikzcd}\ .
	\end{align}	
	Proposition~\ref{prop:indgeometric_descent} guarantees that the vertical maps are equivalences. It is therefore enough to prove that the bottom horizontal arrow is an equivalence. Applying Lemma~\ref{lem:truncatedness_mapping_spaces}, we see that all the mapping spaces are $(n+m)$-truncated. We can therefore replace $\bfDelta$ by $\bfDelta_{\leqslant n+m+2}$, whence the conclusion.
	
	We now prove \eqref{item:compactness-2}. For every index $\alpha$, the map $i_\alpha \colon U_\alpha \to F$ is representable by open Zariski immersions. It immediately follows that the canonical map
	\begin{align}
		U_\alpha \longrightarrow U_\alpha \times_F U_\alpha 
	\end{align}
	is an equivalence. In other words, the map $i_\alpha$ is $(-1)$-truncated. It follows that the induced map $\Map_{\dSt_k}(S,U_\alpha) \to \Map_{\dSt_k}(S,F)$ is $(-1)$-truncated as well in $\Spc$. Since the colimit is filtered, it follows that the map \eqref{eq:qc_exhaustion} is also $(-1)$-truncated. To complete the proof, it is enough to check that it is surjective on connected components. Let therefore $S \to F$ be a morphism in $\dSt_k$. Since $i_\alpha$ is representable by open Zariski immersions, we immediately see that the canonical map
	\begin{align}
		\Map_{\dSt_k}(S, U_\alpha) \longrightarrow \Map_{\dSt_k}(\Red{S}, U_\alpha) \times_{\Map_{\dSt_k}(\Red{S},F)} \Map_{\dSt_k}(S,F)
	\end{align}
	is an equivalence.
	In other words: a map $f \colon S \to F$ factors through $U_\alpha$ if and only if its restriction $\Red{f} \colon \Red{S} \to F$ factors through $U_\alpha$.
	We can therefore replace $S$ by $\Red{S}$; in other words, we can assume from the very beginning that $S \in \dGeomqc_k$. In this case, set $V_\alpha \coloneqq S \times_F U_\alpha$. This gives an increasing open exhaustion of $S$ by the open $V_\alpha$. Since $S$ is quasi-compact, there must exist an index $\alpha$ such that $S = V_\alpha$, so that $f$ must factor through $U_\alpha$. The proof is thus complete.
\end{proof}

\begin{lemma}\label{lem:topological_invariance}
	Let $F$ be a derived stack and let $G \to \Red{F}$ be a morphism representable by open Zariski immersions. Then, there exists an essentially unique derived stack $\overline{G}$ and dashed arrows making the square
	\begin{align}
		\begin{tikzcd}[ampersand replacement=\&]
			G \arrow[dashed]{r} \arrow{d} \& \overline{G} \arrow[dashed]{d} \\
			\Red{F} \arrow{r} \& F
		\end{tikzcd} 
	\end{align}
	a pullback. Moreover, the map $\overline{G} \to F$ is representable by open Zariski immersions.
\end{lemma}

\begin{proof}
	Recall that the canonical map
	\begin{align}
		(\Red{F})_{\mathsf{dR}} \longrightarrow F_{\mathsf{dR}} 
	\end{align}
	is an equivalence.\footnote{Here $(-)_{\mathsf{dR}}$ is the de Rham stack. Concretely, for any derived stack $F$ and every affine scheme $S$, one has $F_{\mathsf{dR}}(S) \coloneqq F(\Red{S})$.}
	We define $\overline{G}$ as the fiber product
	\begin{align}
		\begin{tikzcd}[ampersand replacement=\&]
			G \arrow[dashed]{r} \arrow{d} \& \overline{G} \arrow{r} \arrow{d} \& G_{\mathsf{dR}} \arrow{d} \\
			\Red{F} \arrow{r} \& F \arrow{r} \& F_{\mathsf{dR}} 
		\end{tikzcd} \ .
	\end{align}
	The canonical map $G \to G_{\mathsf{dR}}$ makes the outer rectangle commutative, so the existence of the dashed arrow follows immediately. Furthermore, unraveling the definitions we see that the outer rectangle is in fact a pullback, so the same holds for the square on the left. We are thus reduced to check that the map $G_{\mathsf{dR}} \to F_{\mathsf{dR}}$ is representable by open Zariski immersions. Fix $S \in \dAff_k$ and a morphism $x \colon S \to F_{\mathsf{dR}}$, corresponding to a morphism $\Red{S} \to F$. By assumption, the map $U \coloneqq S \times_{\Red{F}} G \to S$ is an open Zariski immersion. It is thus enough to check that the square
	\begin{align}
		\begin{tikzcd}[ampersand replacement=\&]
			U \arrow{r} \arrow{d} \& G_{\mathsf{dR}} \arrow{d} \\
			S \arrow{r} \& F_{\mathsf{dR}}
		\end{tikzcd}
	\end{align}
	is a pullback. Notice that both vertical maps are $(-1)$-truncated. Then, the conclusion follows by observing that for every $T \in \dAff_k$, a map $T \to S$ factors through $U$ if and only if $\Red{T} \to \Red{S}$ factors through $\Red{U}$, hence if and only if the composite $\Red{T} \to \Red{F}$ factors through $G$.
\end{proof}

\begin{theorem}\label{thm:checking_admissibility_on_reduced}
	Let $\scrX$ be an indgeometric derived stack.
	\begin{enumerate}\itemsep=0.2cm
		\item \label{item:checking_admissibility_on_reduced-1} $\scrX$ is a qcqs indgeometric derived stack if and only if $\Red{\scrX}$ is a qcqs geometric stack.
		
		\item \label{item:checking_admissibility_on_reduced-2} $\scrX$ is an admissible indgeometric derived stack if and only if $\Red{\scrX}$ is a quasi-separated geometric stack.
	\end{enumerate}
\end{theorem}

\begin{proof}
	We first establish the ``only if'' implications. Fix a presentation $\{X_\alpha\}_{\alpha \in I}$ for $F$. Then it follows that
	\begin{align}
		\Red{F} \simeq \colim_{\alpha \in I} \Red{X}_\alpha \ ,
	\end{align}
	the colimit being computed in $\PreSt_k$. If $F$ is qcqs, the above colimit becomes eventually constant, so both statements follow automatically.
	
	We now complete the proof of point \eqref{item:checking_admissibility_on_reduced-1}. Assume that $\Red{\scrX}$ is a qcqs geometric derived stack. Pick any qcqs presentation $\{X_\alpha\}_{\alpha \in I}$ for $\scrX$. The canonical map
	\begin{align}
		\colim_{\alpha \in I} \Red{X}_\alpha \longrightarrow \Red{\scrX} 
	\end{align}
	is then an isomorphism. Lemma~\ref{lem:compactness} guarantees that $\Red{\scrX}$ is a retract of $\Red{X}_\alpha$ for some $\alpha \in I$. Thus, for every $\beta \in I_{\alpha /}$, $\Red{\scrX}$ is also a retract of $\Red{X}_\beta$. Moreover, since we required the transition maps to be closed embeddings, we see that they are $(-1)$-truncated morphisms in the category of \textit{reduced} geometric derived stacks. Thus, the map $\Red{X}_\beta \to \Red{\scrX}$ is also $(-1)$-truncated, and therefore, having a retract, it is an equivalence. It follows that for every $\beta \in I_{\alpha /}$, the map $X_\beta \to \scrX$ is a nil-equivalence, and therefore that the transition maps in $\{X_\beta\}_{\beta \in I_{\alpha/}}$ are nil-equivalences. Thus, the conclusion follows from Remark~\ref{rem:cofinal_presentations}.
	
	Now complete the proof of point \eqref{item:checking_admissibility_on_reduced-2}. Assume therefore that $\Red{\scrX}$ is a quasi-separated geometric derived stack. Then \cite[Lemma~A.1]{Porta_Sala_Hall} allows to pick a quasi-compact open exhaustion $\{U_\alpha\}_{\alpha \in I}$ of $\Red{\scrX}$. Since $\Red{\scrX}$ is quasi-separated, every the same goes for every $U_\alpha$. Applying Lemma~\ref{lem:topological_invariance}, we see that there exists a uniquely defined derived stack $V_\alpha$ equipped with a Zariski open immersion $V_\alpha \to \scrX$ such that the square
	\begin{align}
		\begin{tikzcd}[ampersand replacement=\&]
			U_\alpha \arrow{r} \arrow{d} \& V_\alpha \arrow{d} \\
			\Red{\scrX} \arrow{r} \& \scrX
		\end{tikzcd}
	\end{align}
	a pullback. It follows that $V_\alpha$ is an indgeometric derived stack, and that $\Red{V}_\alpha \simeq U_\alpha$. The uniqueness of $V_\alpha$ immediately implies that if $U_\alpha \subseteq U_\beta$ then $V_\alpha \subseteq V_\beta$. Thus, the first half of the proof guarantees that $V_\alpha$ is a qcqs geometric derived stack. Since $\{V_\alpha\}_{\alpha \in I}$ obviously is an open exhaustion of $\scrX$, it follows that $\scrX$ is an admissible indgeometric derived stack.
\end{proof}

\begin{corollary}\label{cor:qc_open_in_admissible_is_qcqs}
	Let $\scrX$ be an admissible indgeometric derived stack and let $\scrU \to \scrX$ be a morphism representable by quasi-compact open Zariski immersions. Then, $\scrU$ is a qcqs indgeometric derived stack.
\end{corollary}

\begin{proof}
	The representability assumption implies that $\scrU \to \scrX$ is a convergent morphism. Since $\scrX$ is convergent, the same goes for $\scrU$. Besides, if $\{X_\alpha\}_{\alpha \in I}$ is a presentation for $\scrX$, then it follows formally that $\{\scrU \times_\scrX X_\alpha\}_{\alpha \in I}$ is a presentation for $\scrU$. Therefore, $\scrU$ is an indgeometric derived stack. Besides, there is a canonical morphism
	\begin{align}
		\Red{\scrU} \longrightarrow \scrU \times_{\scrX} \Red{\scrX} \ , 
	\end{align}
	and since $\scrU \to \scrX$ is representable by open Zariski immersions, we immediately see that this morphism is an equivalence. Since $\Red{\scrX}$ is a quasi-separated geometric stack and $\Red{\scrU}$ is a quasi-compact open inside $\Red{\scrX}$, we obtain that $\Red{\scrU}$ is qcqs geometric stack and the conclusion follows from Theorem~\ref{thm:checking_admissibility_on_reduced}--\eqref{item:checking_admissibility_on_reduced-1}.
\end{proof}

\begin{corollary}\label{cor:qcqs_indgeometric_nilequivalences_in_every_presentation}
	Let $\scrX$ be a qcqs indgeometric derived stack and let $\{X_\alpha\}_{\alpha \in I}$ be any presentation. Define $I'$ to be the full subcategory of $I$ spanned by those indexes $\alpha$ for which the structural morphism $X_\alpha \to \scrX$ is a nil-equivalence. Then $I'$ is a cofinal subcategory of $I$, and in particular $\{X_\alpha\}_{\alpha \in I'}$ is again a presentation for $\scrX$.
\end{corollary}

\begin{proof}
	Consider the structural morphism $\Red{\scrX} \to \scrX$. Since $\scrX$ is a qcqs indgeometric derived stack, Theorem~\ref{thm:checking_admissibility_on_reduced} implies that $\Red{\scrX}$ is a qcqs geometric stack. Besides,
	\begin{align}
		\Red{\scrX} \simeq \colim_{\alpha \in I} \Red{X}_\alpha \ . 
	\end{align}
	Thus, Lemma~\ref{lem:compactness} implies that $\Red{\scrX}$ is a retract of $\Red{X}_\alpha$ for some $\alpha \in I$. Since $\Red{X}_\alpha \to \Red{\scrX}$ is a monomorphism, it follows that $\Red{\scrX} \simeq \Red{X}_\alpha$. This shows that $I'$ is non-empty. Replacing $I$ by $I_{\beta/}$ for any $\beta \in I$, we conclude that $I'$ is a cofinal subcategory of $I$.
\end{proof}

In order to define Borel-Moore homology of admissible indgeometric derived stacks in \S\ref{sec:BM-homology}, we need to make the ind-structure on such objects explicit. Following \cite[\S\ref*{torsion-pairs-appendix:ind_objects}]{DPS_Torsion-pairs}, we consider the restricted Yoneda embeddings
\begin{align}
	\Phi_{\ind\textrm{-}\qcqs} \colon \indGeomindqcqs_k \longrightarrow \PSh(\tensor*[^{<\infty}]{\dGeom}{^{\qcqs}_k})  \ , 
\end{align}
and
\begin{align}
	\Phi_{\qcqs} \colon \indGeomadm_k \longrightarrow \PSh(\indGeomqcqs_k) \ . 
\end{align}
The choices of $\tensor*[^{<\infty}]{\dGeom}{^{\qcqs}_k}$ and $\indGeomqcqs_k$ are justified by the following result.

\begin{theorem}\label{thm:indization}
	\hfill
	\begin{enumerate}\itemsep=0.2cm
		\item \label{item:indization-1} The functor $\Phi_{\ind\textrm{-}\qcqs}$ commutes with limits, is fully faithful and factors through $\Ind(\tensor*[^{<\infty}]{\dGeom}{^{\qcqs}_k})$.
		Furthermore, if $\{X_\alpha\}_{\alpha \in I}$ is a qcqs presentation for $\scrX \in \indGeomindqcqs_k$, then
		\begin{align}
			\Phi_{\ind\textrm{-}\qcqs}(\scrX) \simeq \fcolim_{(n, \alpha) \in \N \times I} \sft_{\leqslant n} X_\alpha \ , 
		\end{align}
		where $\sft_{\leqslant n}$ denotes the $n$-th Postnikov truncation functor.
		
		\item \label{item:indization-2} The functor $\Phi_{\qcqs}$ commutes with limits, is fully faithful and factors through $\Ind(\indGeomqcqs_k)$.
		Furthermore, if $\{U_\alpha\}_{\alpha \in I}$ is a quasi-compact open exhaustion of $F$, then
		\begin{align}
			\Phi_{\qcqs}(\scrX) \simeq \fcolim_{\alpha \in I} X_\alpha \ . 
		\end{align}
	\end{enumerate}
\end{theorem}

\begin{proof}
	Choose an ind-qcqs presentation $\{X_\alpha\}_{\alpha \in I}$ for $\scrX$. Using Lemma~\ref{lem:indstack_eventually_coconnective_presentation}, we see that
	\begin{align}
		\{\sft_{\leqslant n}(X_\alpha)\}_{(\alpha, n) \in I \times \N} 
	\end{align}
	is again a presentation for $\scrX$, whose pieces are eventually coconnective. At this point, the conclusion follows combining \cite[Proposition~\ref*{torsion-pairs-prop:indization}]{DPS_Torsion-pairs} and Lemma~\ref{lem:compactness}.
\end{proof}

\begin{notation}\label{notation:indization}
	The two parts of Theorem~\ref{thm:indization} provide us with a limit-preserving and fully faithful functor
	\begin{align}\label{eq:indization}
		(-)_{\ind} \colon \indGeomadm_k \to \Ind(\indGeomqcqs) \to \Ind(\Ind(\tensor*[^{<\infty}]{\dGeom}{^{\qcqs}_k})) \ ,
	\end{align}
	which we refer to as the \textit{indization functor}.
\end{notation}

\begin{corollary}\label{cor:indgeometric_finite_limits}
	The full subcategory $\indGeomindqcqs_k \subset \PreSt_k$ is closed under finite limits.
\end{corollary}

\begin{proof}
	It is enough to consider the case of fiber products. For this, apply first Corollary~\ref{cor:indP_objects_fiber_products} to $\Ind( \tensor*[^{<\infty}]{\dGeom}{^{\qcqs}_k} )$ with $P$ the collection of closed immersions, and then combine it with Theorem~\ref{thm:indization}--\eqref{item:indization-1}.
\end{proof}

\begin{corollary}\label{cor:admissible_stacks_finite_limits}
	The category $\indGeomadm_k$ is closed under finite limits in $\PreSt_k$.
\end{corollary}

\begin{proof}
	It suffices to show that $\indGeomadm_k$ is closed under fiber products. Consider therefore a fiber product $\scrW \coloneqq \scrX \times_{\scrY} \scrZ$ of admissible indgeometric stacks. Corollary~\ref{cor:indgeometric_finite_limits} implies that $\scrW$ is again indgeometric. On the other hand,
	\begin{align}
		\Red{\scrW} \simeq \Red{\big(}\Red{\scrX} \times_{\Red{\scrY}} \Red{\scrZ} \big) \ . 
	\end{align}
	Then the conclusion follows combining Theorem~\ref{thm:checking_admissibility_on_reduced} and the fact that quasi-separated geometric derived stacks are closed under fiber products.
\end{proof}

We will need one more criterion to establish admissibility. To state it, we first prove:
\begin{lemma}\label{lem:smooth_implies_relatively_reduced}
	Let $f \colon \scrX \to \scrY$ be a morphism of derived indgeometric stacks. Assume that $f$ is representable by smooth geometric stacks. Then the square
	\begin{align}
		\begin{tikzcd}[ampersand replacement=\&]
			\Red{\scrX} \arrow{r} \arrow{d} \& \scrX \arrow{d} \\
			\Red{\scrY} \arrow{r} \& \scrY
		\end{tikzcd} 
	\end{align}
	is a pullback.
\end{lemma}

\begin{proof}
	Choose a presentation $\{Y_\alpha\}_{\alpha \in I}$ for $\scrY$ and set $X_\alpha \coloneqq Y_\alpha \times_{\scrY} \scrX$. Since $f$ is representable by geometric stacks, it follows formally that $\{X_\alpha\}_{\alpha \in I}$ is a presentation for $\scrX$. Then Proposition~\ref{prop:indgeometric_descent} guarantees that
	\begin{align}
		\scrY \simeq \colim_{\alpha \in I} Y_\alpha \qquad \text{and} \qquad \scrX \simeq \colim_{\alpha \in I} X_\alpha 
	\end{align}
	hold in $\dSt_k$. In particular, since $\Red{j} \colon \reddAff_k \to \dAff_k$ is cocontinuous, it follows from \cite[Lemma~2.20]{Porta_Yu_Higher_Analytic_Stacks} that both $\Red{j}_!$ and $\Red{j}^\ast$ commute with colimits and hence that the canonical morphisms
	\begin{align}
		\colim_{\alpha \in I} \Red{Y}_\alpha \longrightarrow \Red{\scrY} \quad \text{and} \quad \colim_{\alpha \in I} \Red{X}_\alpha \longrightarrow \Red{\scrX} 
	\end{align}
	are equivalences. Since these colimit are filtered, it follows that it suffices to argue that for each $\alpha \in I$ the square
	\begin{align}
		\begin{tikzcd}[ampersand replacement=\&]
			\Red{X}_\alpha \arrow{r} \arrow{d} \& X_\alpha \arrow{d}{f_\alpha} \\
			\Red{Y}_\alpha \arrow{r} \& Y_\alpha
		\end{tikzcd} 
	\end{align}
	is a pullback. Notice that after applying $\Red{j}^\ast$, the horizontal arrows become equivalences and therefore the square becomes a pullback. It suffices then to argue that the fiber product $\Red{Y}_\alpha \times_{Y_\alpha} X_\alpha$ is automatically reduced. Since $f_\alpha$ is smooth by assumption, this is automatic.
\end{proof}

\begin{corollary}\label{cor:testing_admissibility_on_an_atlas}
	Let
	 \begin{align}
	 	\begin{tikzcd}[ampersand replacement=\&]
			\scrX_U \arrow{r} \arrow{d} \& \scrX \arrow{d} \\
			U \arrow{r}{u} \& Y
		\end{tikzcd} 
	 \end{align}
	be a pullback square of derived stacks. Assume that $Y$ is a geometric derived stack, that $u$ is a smooth epimorphism and that $\scrX$ is indgeometric. Then $\scrX$ is admissible if and only if $\scrX_U$ is admissible.
\end{corollary}

\begin{proof}
	Both $U$ and $\scrY$ are admissible, so if $\scrX$ is admissible Corollary~\ref{cor:admissible_stacks_finite_limits} implies that $\scrX_U$ is admissible.
	
	For the vice-versa, we observe that since $u$ is smooth the square
	\begin{align}
		\begin{tikzcd}[ampersand replacement=\&]
			\Red{U} \arrow{r} \arrow{d} \& \Red{Y} \arrow{d} \\
			U \arrow{r} \& Y
		\end{tikzcd} 
	\end{align}
	is a pullback square.
	Let
	\begin{align}
		\scrX' \coloneqq \Red{U} \times_{\Red{Y}} \Red{\scrX} \ . 
	\end{align}
	The canonical morphism $\scrX' \to \scrX_U$ induces an equivalence $\Red{\scrX}' \simeq \Red{\scrX}_U$. Thus, if $\scrX_U$ is admissible, $\Red{\scrX}_U$ is a quasi-separated by Theorem~\ref{thm:checking_admissibility_on_reduced}. By assumption, the canonical morphism $\scrX' \to \Red{\scrX}$ is representable by smooth derived stacks. In particular, Lemma~\ref{lem:smooth_implies_relatively_reduced} guarantees that $\scrX'$ is itself reduced, and therefore it is a quasi-separated geometric stack. Since $\Red{U} \to \Red{Y}$ is again a smooth epimorphism, we deduce from \cite[Corollary~1.3.4.5]{TV_HAG-II} that $\Red{\scrX}$ is itself a quasi-separated geometric stack, so the conclusion follows from Theorem~\ref{thm:checking_admissibility_on_reduced}.
\end{proof}

\bigskip\section{Borel-Moore homology for admissible indgeometric stacks}\label{sec:BM-homology}

In this section we define Borel-Moore homology for admissible indgeometric stacks, and establish its basic functorialities. We use in an essential way the motivic theory developed by A.\ Khan and its collaborators (see e.g. \cite{Khan_VFC,Khan_Voevodsky_criterion}). The extension from geometric to admissible indgeometric is a formal procedure, essentially made possible thanks to the theory developed in \cite[\S\ref*{torsion-pairs-sec:homological-invariants}]{DPS_Torsion-pairs}. We refer to \textit{loc. cit.} for all the relevant notation.

\medskip

Fix a \textit{motivic formalism} $\bfD^\ast$ (in the sense of \cite[Definition~\ref*{torsion-pairs-def:motivic_formalism}]{DPS_Torsion-pairs}) and its associated six-functors formalism
\begin{align}
	\proD^\ast_! \colon \Corr\big( \dGeom^{\mathsf{qs}}_k \big)_{\mathsf{lft},\mathsf{all}} \longrightarrow \CAlg(\Pro(\Cat_\infty)) \ .
\end{align}
Let $S \in \dGeomqcqs_k$ be a qcqs geometric derived stack. For an oriented $\calA \in \CAlg(\bfD(S))$ and an abelian subgroup $\Gamma \subset \Pic(\bfD(S))$ closed under Thom twists.
In \cite[\S\ref*{torsion-pairs-subsec:beyond_qc}]{DPS_Torsion-pairs}, the first three-named authors defined for a quasi-separated geometric derived stack locally of finite type $\scrX$ over $S$ its \textit{relative Borel-Moore homology groups with coefficients in $\calA$}:
\begin{align}
	\HBMDGamma_0(\scrX / S; \calA ) \ , 
\end{align}
as a $\Gamma$-graded \textit{topological} abelian group.
This is an extension of Khan's framework \cite{Khan_VFC}, and the topology keeps track of the non-quasi-compactness of $\scrX$.

\begin{remark}[Underlying topology]\label{rem:quasi-compact_topology}
	Let us take $S = \Spec(k)$.
	It follows combining \cite[Remarks~\ref*{torsion-pairs-rem:quasi-compact_topology} and \ref*{torsion-pairs-rem:underlying_topology_II}]{DPS_Torsion-pairs} that $\HBMDGamma_0(\scrX/S;\calA)$ is explicitly given by the formula
	\begin{align}
		\HBMDGamma_0(\scrX/S;\calA) \simeq \bigoplus_{\alpha \in \pi_0(\scrX)} \lim_{\scrU \subset_{\mathsf{qc}} \scrX_\alpha} \HBMDGamma_0(\scrU/S;\calA) \ , 
	\end{align}
	where:
	\begin{enumerate}\itemsep=0.2cm
		\item the direct sum is taken over the set of connected components of $\scrX$, and if $\alpha \in \pi_0(\scrX)$, we denote by $\scrX_\alpha$ the corresponding connected component;
		
		\item the limit is taken over the set of quasi-compact open substacks of $\scrX_\alpha$;
		
		\item each group $\HBMDGamma_0(\scrU/S;\calA)$ can be expressed in terms of the six-operations, and essentially coincide with Khan's motivic Borel-Moore homology (except for the fact that we are summing together all twists by elements in $\Gamma$, see \cite[Remark~\ref*{torsion-pairs-rem:explicit_bigrading}]{DPS_Torsion-pairs} for more on this point);
		
		\item both the inverse limit and the coproduct are topologized in the natural way, endowing the building pieces $\HBMDGamma_0(\scrU/S;\calA)$ with the discrete topology.
	\end{enumerate}
	Technically speaking, the results obtained in \cite{DPS_Torsion-pairs} are strictly stronger, in that we do not work with topological abelian groups, but with free categorical (co)completions.
	However, in applications only this underlying topological structure will be taken into account.
\end{remark}

Contrary to \cite{Khan_VFC}, the above homology groups are not obtained directly out of the six-functors formalism, as there is a certain \textit{renormalization} procedure that needs to be carried out and that allows to have larger functoriality properties (see \cite[Remarks~\ref*{torsion-pairs-rem:difference_with_extended} and \ref*{torsion-pairs-rem:Gaitsgory_renormalization}]{DPS_Torsion-pairs}).

\begin{definition}\label{def:rpas-connected}
	A morphism $f \colon \scrX \to \scrY$ is said to be \textit{locally rpas} if for every connected component $\scrX_0\subset \scrX$ the induced map $\scrX_0 \to \scrY$ is representable by proper algebraic spaces; and that is said to be \textit{finitely connected} if for every connected component $\scrY_0 \subset \scrY$ the preimage $\scrX\times_\scrY \scrY_0$ has finitely many connected components.
\end{definition}
With this terminology, we establish in \cite[\S\ref*{torsion-pairs-subsec:beyond_qc}]{DPS_Torsion-pairs} the following functoriality properties for $\HBMDGamma_\ast(-/S;\scrA)$:
\begin{itemize}\itemsep=0.2cm
	\item for derived lci, quasi-compact and finitely connected morphisms $f \colon \scrX \to \scrY$, a \textit{Gysin pullback}
	\begin{align}
		f^! \colon \HBMDGamma_0(\scrY / S; \calA) \longrightarrow \HBMDGamma_0(\scrX/S;\calA) \ , 
	\end{align}
	which is continuous for the natural topology on both sides;
	\item for locally rpas morphisms $f \colon \scrX \to \scrY$, a proper pushforward
	\begin{align}
		f_\ast \colon \HBMDGamma_0(\scrX / S ; \calA) \longrightarrow \HBMDGamma_0(\scrY/S;\calA) \ . 
	\end{align}
\end{itemize}
Together with the exterior products (essentially induced by the algebra structure on $\calA$), these properties produce a lax-monoidal functor
\begin{align}\label{eq:BM_quasi_separated_stacks}
	\HBMDGamma_0(-/S; \calA) \colon \Corr^\times(\dGeomqs_S)_{\mathsf{qc.lci}\:\cap\:\mathsf{fconn},\:\mathsf{lrpas}} \longrightarrow \Pro^\sqcup(\Modd_R^\heartsuit) \ , 
\end{align}
see \cite[Theorem~\ref*{torsion-pairs-thm:functoriality_of_genuine_BM_homology_admissible}]{DPS_Torsion-pairs}.

We now extend this framework to admissible indgeometric derived stacks. In light of Theorem~\ref{thm:checking_admissibility_on_reduced}, we have a well defined functor
\begin{align}
	\Red{(-)} \colon \indGeomadm_S \longrightarrow \mathsf{Geom}^{\mathsf{qs}}_S \ .
\end{align}

We extend the terminology above:
\begin{definition}\label{def:admissible-rpas-connected}
	A morphism $f \colon \scrX \to \scrY$ of admissible indgeometric derived stacks is \textit{finitely connected} (resp.\ \textit{locally rpas}) if the underlying morphism $\Red{f} \colon \Red{\scrX} \to \Red{\scrY}$ has the same property. On the other hand, we say that $f$ is \textit{quasi-compact} (resp.\ \textit{derived lci}) if it is representable by quasi-compact (resp.\ derived lci) geometric derived stacks. 	
\end{definition}

We now define 
\begin{align}
	\HBMDGamma_0(-/S; \calA) \colon \Corr^\times(\indGeomadm_S)_{\mathsf{qc.lci}\:\cap\:\mathsf{fconn},\:\mathsf{lrpas}} \longrightarrow \Pro^\sqcup(\Modd_R^\heartsuit) \ , 
\end{align}
as the composition of \eqref{eq:BM_quasi_separated_stacks} with the functor $\Red{(-)}$.
The following three remarks concern operations for the resulting homology theory:

\begin{remark}[Proper pushforward]
	Fix a qcqs geometric derived stack $S$ and a morphism $f \colon \scrX \to \scrY$ of admissible indgeometric derived stacks over $S$.
	If $f$ is locally rpas, the same goes, by definition, for $\Red{f}$.
	In particular the proper pushforward of \cite[Theorem~\ref*{torsion-pairs-thm:functoriality_of_genuine_BM_homology_admissible}]{DPS_Torsion-pairs} induces a well defined, continuous map
	\begin{align}
		f_\ast \colon \HBMDGamma_0(\scrX / S; \calA) \longrightarrow \HBMDGamma_0(\scrY / S; \calA) \ . \tag*{\qedhere} 
	\end{align}
\end{remark}

\begin{remark}[Gysin pullback]\label{rem:extension-functoriality}
	Fix a qcqs geometric derived stack $S$ and a morphism $f \colon \scrX \to \scrY$ of admissible indgeometric derived stacks over $S$. Assume that $f$ is quasi-compact, derived lci and finitely connected. In this situation, the morphism $f' \colon \scrX \times_{\scrY} \Red{\scrY} \to \Red{\scrY}$ has the same properties, and it is nil-equivalent to $f$. Besides, $\Red{\scrY}$ is a quasi-separated geometric derived stack by Theorem~\ref{thm:checking_admissibility_on_reduced}, and the same goes for
	\begin{align}
		\Red{\left( \scrX \times_{\scrY} \Red{\scrY} \right)} \simeq \Red{\scrX} \ . 
	\end{align}
	Thus, the Gysin pullback of \cite[Theorem~\ref*{torsion-pairs-thm:functoriality_of_genuine_BM_homology_admissible}]{DPS_Torsion-pairs} yields a well defined map
	\begin{align}
		f^! \colon \HBMDGamma_0(\scrY / S; \calA) \longrightarrow \HBMDGamma_0(\scrX / S; \calA) \ . \tag*{\qedhere} 
	\end{align}
\end{remark}

\begin{remark}[Exterior products]\label{rem:exterior-product}
	Fix a qcqs geometric derived stack $S$ and two admissible indgeometric derived stacks $\scrX$ and $\scrY$ over $S$. Via the canonical equivalence
	\begin{align}
		\Red{ ( {\Red{\scrX} \times_{\Red{S}} \Red{\scrY}} ) } \simeq \Red{(\scrX \times_S \scrY)} \ , 
	\end{align}
	the exterior product of \cite[Theorem~\ref*{torsion-pairs-thm:functoriality_of_genuine_BM_homology_admissible}]{DPS_Torsion-pairs} yields a well defined, continuous map
	\begin{align}
		\boxtimes \colon \HBMDGamma_0(\scrX / S; \calA) \otimes \HBMDGamma_0(\scrY / S; \calA) \longrightarrow \HBMDGamma_0(\scrX \times_S \scrY / S; \calA) \ . \tag*{\qedhere} 
	\end{align}
\end{remark}

We can summarize this discussion stating the following:
\begin{theorem}\label{thm:functoriality_of_BM_homology_admissible}
	Let $S \in \dGeomqcqs_k$, $\calA \in \CAlg(\bfD(S))$ and let $\Gamma \subseteq \Pic(\bfD^\ast(S))$ be an abelian subgroup.
	Assume that $\calA$ is oriented and that $\Gamma$ is closed under Thom twists.
	Then the construction
	\begin{align}
		\HBMDGamma_0(-/S; \calA) \colon \Corr^\times(\indGeomadm_S)_{\mathsf{qc.lci}\:\cap\:\mathsf{fconn},\:\mathsf{lrpas}} \longrightarrow \Pro^\sqcup(\catMod_R^\heartsuit)
	\end{align}
	that sends $\scrX \to S$ to the topological $R$-module $\HBMDGamma_0(\scrX/S;\calA)$ (disregarding the extra $\Gamma$-grading) and a correspondence
	\begin{align}
		\begin{tikzcd}[column sep=small,ampersand replacement=\&]
			{} \& \scrZ \arrow{dr}{p} \arrow{dl}[swap]{f} \\
			\scrX \& \& \scrY
		\end{tikzcd} \ ,
	\end{align}
	where $f$ is representable by finitely connected and quasi-compact derived lci geometric derived stacks while $p$ is locally rpas, to the composition
	\begin{align}
		p_\ast \circ f^! \colon \HBMDGamma_0(\scrX / S; \calA) \longrightarrow \HBMDGamma_0(\scrY / S; \calA)
	\end{align}
	defines a lax-monoidal functor.
	
	Moreover, if $(\bfD', \calA', \Gamma')$ is a second motivic formalism with a choice of an oriented ring of coefficients $\calA'$ and an abelian subgroup $\Gamma'$ closed under Thom twist, then a morphism
	\begin{align}
		(s,\phi) \colon (\bfD, \calA, \Gamma) \longrightarrow (\bfD', \calA', \Gamma') 
	\end{align}
	induces a lax symmetric monoidal transformation
	\begin{align}
		\HBMDGamma_0(-/S;\calA) \longrightarrow \sfH^{\bfD', \Gamma'}_0(-/S;\calA') \ . 
	\end{align}
\end{theorem}

\begin{example}\label{ex:three-H}
	Let $S \in \dGeomqcqs_k$ and $\scrX\in \indGeomadm_S$.
	\begin{enumerate}\itemsep=0.2cm
		\item when $\bfD^\ast \coloneqq \mathbf{DM}_\Q^\ast$ is the \textit{rational Voevodsky's formalism} (cf.\ \cite[Example~\ref*{torsion-pairs-ex:genuine_motivic_formalism}]{DPS_Torsion-pairs}) and $\calA \coloneqq \sfH \Q$ is the motivic Eilenberg-MacLane $\E_\infty$-ring spectrum, and $\Gamma \coloneqq \Z\langle 1 \rangle$, we write
		\begin{align}
			\Hmotbullet(\scrX/S; 0)\coloneqq \HBMDGamma_0(\scrX/S; \sfH \Q) \ .
		\end{align}
		This recovers (motivically defined) Chow groups of $\scrX$. 
		
		\item When $\bfD^\ast \coloneqq \mathbf{DM}_\Q^\ast$ is the rational Voevodsky's formalism and $\calA \coloneqq \mathsf{KGL}^{\mathsf{et}}$ is the étale hypersheafification of the algebraic $K$-theory spectrum, and $\Gamma \coloneqq \Z\langle 1 \rangle$, we write
		\begin{align}
			\bfG(\scrX/S)\coloneqq \HBMDGamma_0(\scrX/S; \mathsf{KGL}^{\mathsf{et}}) \ .
		\end{align}
		This recovers the étale sheafification of the algebraic $G$-theory of $\scrX$.
		
		\item When $\bfD^\ast \coloneqq \BettiD$ is the \textit{topological formalism} (cf.\ \cite[Example~\ref*{torsion-pairs-ex:topological_formalism}]{DPS_Torsion-pairs}), $\calA \coloneqq \Q$, and $\Gamma \coloneqq \Z\langle 1/2 \rangle$ (see \cite[Remark~\ref*{torsion-pairs-rem:explicit_bigrading}]{DPS_Torsion-pairs}), we simply write
		\begin{align}
			\HBMbullet(\scrX/S)\coloneqq \HBMGamma_0(\scrX/S; \Q) \ .
		\end{align}
		
		We will also be interested in taking $\Gamma \coloneqq \Z \langle 1 \rangle$, in which case we write
		\begin{align}
			 \HBM_{\mathsf{even}}(\scrX/S)\coloneqq \HBMDGamma_0(\scrX/S; \Q)\ .
		\end{align}
		The natural inclusion $\Z\langle 1 \rangle \subset \Z\langle 1/2 \rangle$ induces a continuous morphism of algebras
		\begin{align}
			\HBM_{\mathsf{even}}(\scrX/S)\longrightarrow \HBMbullet(\scrX/S)\ .
		\end{align}
		Also, the natural transformation $\bfDM_\Q \to \BettiD_\Q$ (see \cite[ Remark~\ref*{torsion-pairs-rem:comparison_motivic_to_topological}]{DPS_Torsion-pairs}) induces a continuous map
		\begin{align}
			\Hmotbullet(\scrX/S) \longrightarrow \HBM_{\mathsf{even}}(\scrX/S) \ , 
		\end{align}
		given by the cycle class map. \qedhere
	\end{enumerate}
\end{example}

\begin{remark}\label{rem:equivariant-H}
	Let $S \in \dGeomqcqs_k$ and $\scrX\in \indGeomadm_S$. Assume that there is an action by a torus $T$ on $\scrX$. We consider $S$ equipped with the trivial $T$-action. By construction
	\begin{align}
		S \times_{[S/T]} [\scrX / T]  \simeq \scrX \ . 
	\end{align}
	Therefore, Corollary~\ref{cor:testing_admissibility_on_an_atlas} implies that each $[\scrX / T]$ is admissible. Thus, for fixed $\calA \in \CAlg(\bfD(S))$ and a fixed abelian group $\Gamma \subseteq \Pic(\bfD^\ast(S))$, we set
	\begin{align}
		\mathsf{H}^{\bfD,\Gamma, T}_0(\scrX; \calA) \coloneqq\mathsf{H}^{\mathsf{BM},\pi^\ast\Gamma}_0([\scrX/T]/[S/T]; \pi^\ast \calA) \ ,
	\end{align}
	where $\pi\colon [S/T]\to S$ is the coarse moduli space map. In particular, we denote by $\HBMbulletT(\scrX)$, $\bfG^T(\scrX)$, $\HmotbulletT(\scrX)$, $\sfH_{\mathsf{even}}^T(\scrX)$ the equivariant versions of the examples discussed in Example~\ref{ex:three-H}.
\end{remark}

We fix an abelian group $(\Lambda,+)$. We define the category of \textit{$\Lambda$-graded derived stacks} as:
\begin{align}
	\LambdadSt_k \coloneqq \Fun(\Lambda, \dSt_k) \ , 
\end{align}
and we consider it with the symmetric monoidal structure induced by $\Lambda$ via \textit{Day's convolution} (see \cite[Recollection~\ref*{torsion-pairs-recollection:Day_convolution}]{DPS_Torsion-pairs}). This symmetric monoidal structure propagates to the $\infty$-category of correspondences $\Corr(\LambdadSt_k)$.

\begin{definition}
	Let $P$ be a property of derived stacks (resp.\ of morphisms of derived stacks). A $\Lambda$-graded derived stack $F$ (resp.\ a morphism $F \to G$) in $\Fun(\Lambda, \dSt_k)$ is \textit{said to have the property $P$} if for every $\bfv \in \Lambda$ the derived stack $F(\bfv)$ (resp.\ the morphism $F(\bfv) \to G(\bfv)$) has the property $P$.
\end{definition}

We denote by $\LambdaindGeomadm_k$ the full subcategory of $\LambdadSt_k$ spanned by admissible indgeometric stacks.
Concretely,
\begin{align}
	\LambdaindGeomadm_k \coloneqq \Fun(\Lambda, \indGeomadm_k) \ . 
\end{align}
In the same way, we have a well defined symmetric monoidal $\infty$-category of correspondences
\begin{align}
	\Corr^\times(\LambdaindGeomadm_k)_{\mathsf{qc.lci}\:\cap\:\mathsf{fconn},\:\mathsf{lrpas}} \ . 
\end{align}
Fix now a coefficient ring $R$ of characteristic zero. Once again, Day's convolution endows the category
\begin{align}
	\Lambda\textrm{-}\Pro^\sqcup(\catMod_R^\heartsuit) \coloneqq \Fun(\Lambda, \Pro^\sqcup(\catMod_R^\heartsuit)) 
\end{align}
with a symmetric monoidal structure. Combining the formal properties of Day's convolution (see \cite[Recollection~\ref*{torsion-pairs-recollection:Day_convolution}]{DPS_Torsion-pairs}) with Theorem~\ref{thm:functoriality_of_BM_homology_admissible}, we obtain:
\begin{theorem}\label{thm:BM_Lambda_graded_functoriality}
	Let $\bfD^\ast$ be a motivic formalism. Let $S \in \dGeomqcqs$, $\calA \in \CAlg(\bfD(S))$ and let $\Gamma \subseteq \Pic(\bfD^\ast(S))$ be an abelian subgroup. Assume that $\calA$ is oriented and that $\Gamma$ is closed under Thom twists. Then, the construction
	\begin{align}
		\HBMDGamma_0(-/S; \calA) \colon \Corr^\times(\LambdaindGeomadm_k)_{\mathsf{qc.lci}\:\cap\:\mathsf{fconn},\:\mathsf{lrpas}} \longrightarrow \Lambda\textrm{-}\Pro^\sqcup(\catMod_R^\heartsuit) 
	\end{align}
	that sends $\scrX \to S$ to
	\begin{align}
		\HBMDGamma_0(\scrX/S;\calA) \coloneqq \bigoplus_{\bfv \in \Lambda} \HBMDGamma_0(\scrX(\bfv)/S;\calA)
	\end{align}
	and whose functoriality is given as in Theorem~\ref{thm:functoriality_of_BM_homology_admissible} defines a lax symmetric monoidal functor.
\end{theorem}

\begin{remark}
	All the constructions described in Example~\ref{ex:three-H} and Remark~\ref{rem:equivariant-H} carry over to the $\Lambda$-graded setting.
\end{remark}

\bigskip\section{Derived moduli of coherent sheaves on formal schemes}\label{sec:nilpotent_coherent_sheaves}

In this section we study families of nilpotent sheaves on formal schemes. All throughout this section, we fix a field $k$ of characteristic zero. Note that all affine schemes are assumed to be \textit{laft} in the sense of Definition~\ref{def:laft-convergent}.

\subsection{Derived formal schemes}

We start by discussing \textit{formal completions}:
\begin{definition}[{\cite[Definition~2.1.3]{CPTVV}}]\label{def:formal_completion}
	Let $X$ be a derived scheme locally almost of finite presentation and let $j \colon Z \to X$ be a closed immersion. The \textit{formal completion $\widehat{X}_Z$ of $X$ along $Z$} is defined as the following pullback:
	\begin{align}
		\begin{tikzcd}[ampersand replacement = \&]
			\widehat{X}_Z \arrow{r} \arrow{d} \& X \arrow{d}{\lambda_X} \\
			Z_{\mathsf{dR}} \arrow{r}{j_{\mathsf{dR}}} \& X_{\mathsf{dR}} 
		\end{tikzcd} \ . 
	\end{align}
\end{definition}

\begin{definition}\label{def:formal_scheme}
	A derived stack $\calX \in \dSt_k$ is a \textit{derived formal scheme} if there exists a jointly epimorphic collection of morphism $\{ \calU_i \to \calX \}$ which are representable by open Zariski immersions and such that each $\calU_i$ is a formal completion in the sense of Definition~\ref{def:formal_completion}.
	
	We write $\formaldSch_k$ for the full subcategory of $\dSt_k$ spanned by derived formal schemes.
\end{definition}

\begin{remark}\label{rem:formal_scheme_reduced}
	If $\calX$ is a derived formal scheme, then $\Red{\calX}$ is an ordinary scheme. 
\end{remark}

\begin{example}\label{eg:formal_scheme_times_scheme}
	Let $S \in \dSch_k$ be a derived scheme. Then, the induced morphism
	\begin{align}
		(\trunc{S})_{\mathsf{dR}} \longrightarrow S_{\mathsf{dR}} 
	\end{align}
	is an equivalence, and in particular $S$ coincides with the formal completion of $S$ along $\trunc{S}$. More generally, if $Z$ is a closed subscheme of $X$ and $S \in \dSch_k$ is a derived scheme, then $\widehat{X}_Z \times S$ coincides with the formal completion of $X \times S$ along $Z \times \trunc{S_0}$. It follows that if $\calX$ is a derived formal scheme, the same goes for $\calX \times S$ for every derived scheme $S$.
\end{example}

\begin{definition}
	We say that $\calX$ is \textit{quasi-compact} (resp.\ \textit{quasi-separated}) if $\Red{\calX}$ is quasi-compact (resp.\ quasi-separated).
	
	We denote by $\formaldSch_k^{\mathsf{qc}}$ the corresponding full subcategory of $\formaldSch_k$.
\end{definition}

\begin{recollection}\label{recollection:presentation_formal_completion}
	Let $\calX$ be a derived formal scheme and set $Z \coloneqq \Red{\calX}$. A \textit{laft thickening} of $Z$ inside $\calX$ is a factorization of $j \colon \Red{\calX} \to \calX$ as
	\begin{align}
		\begin{tikzcd}[ampersand replacement=\&]
			\Red{\calX} \arrow{r}{j_{Z,W}} \& W \arrow{r}{j_W} \& X
		\end{tikzcd} \ ,
	\end{align}
	where $j_{Z,W}$ is a nil-isomorphism\footnote{i.e., $\Red{j}_{Z,W}$ is an isomorphism.}, $W$ is a laft derived scheme, and $j_W$ is a representable by a closed immersion. We denote by $\scrT_{Z \overunder \calX}$ the full subcategory of $\dSt_{Z \overunder X}$ spanned by laft thickenings of $Z$ inside $X$. It follows from \cite[Chapter~2, Proposition~1.8.6]{Gaitsgory_Rozenblyum_Study_II} that $\scrT_{Z \overunder \calX}$ is filtered and that the canonical map
	\begin{align}
		\colim_{W \in \scrT_{Z \overunder \calX}} W \longrightarrow \calX 
	\end{align}
	is an equivalence, where the colimit is taken in $\PreSt_k$. Finally, notice that all the morphisms in $\scrT_{Z \overunder \calX}$ are closed embeddings.
\end{recollection}

When $\calX = \widehat{X}_Z$, this provides an explicit formula of the formal completion in terms of infinitesimal thickenings. Notice that \textit{a priori} one needs all derived thickenings to compute $\widehat{X}_Z$. Nevertheless, when $X$ is underived then the same goes for $\widehat{X}_Z$:
\begin{proposition}\label{prop:derived_completion_classical}
	Let $j \colon Z \to X$ be a closed immersion of classical schemes almost of finite presentation over $k$. Then, $\widehat{X}_Z$ is underived and it coincides with the classical formal completion of $X$ along $Z$.
\end{proposition}

\begin{proof}
	The question is easily seen to be local on $X$, and in the affine case it follows from \cite[Corollary~7.3.6.9]{Lurie_SAG}. See also \cite[Proposition~2.1.4]{HL_Categorical_properness}.
	\end{proof}

\begin{remark}\label{rem:infinitesimal_neighborhood}
	In \cite[Chapter~9, \S5]{Gaitsgory_Rozenblyum_Study_II}, the authors construct a derived version of the infinitesimal neighborhood of order $n$ of $Z$ inside $X$. We denote this construction by $Z^{(n)}$. When $X$ is underived, its truncation $\trunc{Z^{(n)}}$ coincides with the usual $n$-th order thickening of $Z$ inside $X$. If $X$ is underived and $j$ is a regular closed immersion, $Z^{(n)}$ is underived and it coincides with the usual infinitesimal neighborhood of order $n$, $\Spec_X(\scrO_X/\calJ_Z^{n+1})$, where $\calJ_Z$ is the ideal of definition of $Z$. Finally, the natural functor
	\begin{align}
		\N \longrightarrow \scrT_{Z \overunder X} 
	\end{align}
	sending $n$ to $Z^{(n)}$ is cofinal and therefore there is a canonical equivalence
	\begin{align}
		\colim_n Z^{(n)} \simeq \widehat{X}_Z \ . 
	\end{align}
	Notice that when $X$ is underived, Proposition~\ref{prop:derived_completion_classical} implies that the canonical comparison map
	\begin{align}
		\colim_n \trunc{Z^{(n)}} \longrightarrow \colim_n Z^{(n)} 
	\end{align}
	is an equivalence, even when the individual maps $\trunc{Z^{(n)}} \to Z^{(n)}$ are not.
\end{remark}

Derived formal schemes have representable diagonal. Since this is somehow important in what follows, we include a short proof.
\begin{lemma}\label{lem:formal_completion_to_ambient_-1_truncated}
	Let $X$ be a quasi-compact and quasi-separated derived scheme locally almost of finite presentation and let $j \colon Z \hookrightarrow X$ be a closed immersion. Then, the natural morphism
	\begin{align}
		\widehat{X}_Z \longrightarrow X 
	\end{align}
	is $(-1)$-truncated.
\end{lemma}

\begin{proof}
	Since $(-1)$-truncated morphisms are stable under pullbacks, it is enough to verify that $Z_{\mathsf{dr}} \to X_{\mathsf{dr}}$ is $(-1)$-truncated. Unraveling the definitions, we see that we have to check that for every derived affine scheme $S$ locally almost of finite type the fibers of the morphism
	\begin{align}
		\Map_{\dSt}(\Red{S}, Z) \longrightarrow \Map_{\dSt}(\Red{S}, X) 
	\end{align}
	are empty or contractible. This is true because $Z \to X$ is a closed immersion.
\end{proof}

\begin{proposition}\label{prop:formal_scheme_diagonal_representable}
	Let $\calX$ be a derived formal scheme. Then, the diagonal
	\begin{align}
		\Delta_\calX \colon \calX \longrightarrow \calX \times \calX 
	\end{align}
	is representable by derived schemes.
\end{proposition}

\begin{proof}
	The question is Zariski-local on $\calX$, so we can assume from the very beginning that $\calX \simeq \widehat{X}_Z$ for some quasi-compact and quasi-separated derived scheme locally almost of finite presentation and some closed immersion $j \colon Z \hookrightarrow X$. Let $S = \Spec(A)$ be an affine derived scheme and let $x \colon S \to \widehat{X}_Z \times \widehat{X}_Z$. It follows from Recollection~\ref{recollection:presentation_formal_completion} that we can find a laft thickening $W$ of $Z$ inside $X$ and a factorization of $x$ as
	\begin{align}
		\begin{tikzcd}[column sep=small, ampersand replacement=\&]
			S \arrow{r}{x'} \& W \times W \arrow{r} \& \widehat{X}_Z \times \widehat{X}_Z 
		\end{tikzcd} \ .
	\end{align}
	Since
	\begin{align}
		(W \times W) \times_{\widehat{X}_Z \times \widehat{X}_Z} \widehat{X}_Z \simeq W \times_{\widehat{X}_Z} W \ , 
	\end{align}
	it suffices to prove that this is a derived scheme. Using Lemma~\ref{lem:formal_completion_to_ambient_-1_truncated} we see that the canonical map
	\begin{align}
		W \times_{\widehat{X}_Z} W \longrightarrow W \times_X W 
	\end{align}
	is an equivalence, whence the conclusion.
\end{proof}

The above Proposition immediately implies the following.
\begin{corollary}\label{cor:formal_scheme_diagonal_representable}
	Let $\calX$ be a derived formal scheme and let $S$ be a derived scheme. Then, any morphism $S \to \calX$ is representable by derived schemes.
\end{corollary}

We finish with the following simple recognition criterion for formal completions.
\begin{lemma}\label{lem:recognition_criterion_formal_completions}
	Let $f \colon F \to G$ be a formally étale morphism between laft, convergent and infinitesimally cohesive derived stacks. Then, $f$ exhibits $F$ as the formal completion of $G$ along $\Red{F} \to G$.
\end{lemma}

\begin{proof}
	First of all, observe that the canonical map
	\begin{align}
		(\Red{F})_{\mathsf{dR}} \longrightarrow F_{\mathsf{dR}} 
	\end{align}
	is an equivalence. It is therefore sufficient to argue that the square
	\begin{align}
		\begin{tikzcd}[ampersand replacement=\&]
			F \arrow{r}{f} \arrow{d} \& G \arrow{d} \\
			F_{\mathsf{dR}} \arrow{r} \& G_{\mathsf{dR}}
		\end{tikzcd} 
	\end{align}
	is a pullback. Since both $F_{\mathsf{dR}}$ and $G_{\mathsf{dR}}$ are formally étale over $\Spec(k)$, the bottom horizontal map is formally étale as well.
	Since by assumption $f$ is formally étale, the same goes for the induced morphism
	\begin{align}
		F \longrightarrow \widehat{G}_F \ . 
	\end{align}
	Besides, since $F$ and $G$ are laft, convergent and infinitesimally cohesive, the same goes for $\widehat{G}_F$ and therefore for the above morphism. In particular, Lemma~\ref{lem:testing_equivalences_on_reduced} shows that it suffices to check that the induced morphism $\Red{F} \to \Red{\widehat{G}}_F$ is an equivalence. For this, it suffices to argue that the above square is a pullback when evaluated on $S \in \reddAff_k$. Since for such a choice of $S$ the vertical arrows in the above square become equivalences, the conclusion follows.
\end{proof}

\subsection{Nilpotent sheaves on a formal scheme}

In order to define the moduli of coherent sheaves on a formal scheme, we use \textit{ind-coherent sheaves} \cite{Gaitsgory_IndCoh}, or rather their relative version discussed briefly in Appendix~\ref{appendix:relative_ind_coherent}.

Fix a base derived scheme $S \in \evccndAff_k$ and a morphism $\calX \to S$ which is representable by qcqs derived formal schemes. We set $Z \coloneqq \Red{\calX}$. Notice that $\calX$ is itself (i.e., forgetting the map to $S$) a qcqs derived formal scheme and therefore Recollection~\ref{recollection:presentation_formal_completion} implies that
\begin{align}
	\calX \simeq \colim_{W \in \scrT_{Z \overunder \calX}} W \ , 
\end{align}
the colimit being computed in $\PreSt_{/S}$.
It follows that
\begin{align}\label{eq:IndCoh_colimit}
	\catIndCoh(\calX/S) \simeq \colim_{W \in \scrT_{Z \overunder \calX}} \catIndCoh(W/S) \ ,
\end{align}
where the transition maps are computed with respect to the $(-)_\ast$-functoriality of $\catIndCoh$ as in Construction~\ref{construction:relative_IndCoh_Kan_extension} and the colimit is computed in $\PrL$ (see Remark~\ref{rem:relative_IndCoh_on_formal_schemes}). Notice moreover that each $W$ is qcqs, so by definition
\begin{align}
	\catIndCoh(W/S) \simeq \Ind(\catCohb(W/S)) 
\end{align}
is compactly generated and all the transition maps preserve compact objects. Moreover, Corollary~\ref{cor:formal_scheme_diagonal_representable} implies that each morphism $j_{W,\calX} \colon W \to \calX$ is representable by closed immersions; it follows from Construction~\ref{construction:relative_IndCoh_proper_uppershriek} and Remark~\ref{rem:relative_IndCoh_on_formal_schemes} that we have an adjunction $j_{W,\calX,\ast} \dashv j_{W,\calX}^!$ in $\PrL$ and that in particular $j_{W,\calX,\ast}$ commutes with compact objects. In particular, we obtain a well defined morphism
\begin{align}\label{eq:IndCoh_colimit_PrLomega}
	\colim_{W \in \scrT_{Z \overunder \calX}} \catCohb(W/S) \longrightarrow \catIndCoh(\calX/S)^\omega \ ,
\end{align}
We prove below that this is an equivalence, a statement that is well-known when $S = \Spec(k)$.

We first focus on the case of a formal completion, that is $\calX \simeq \widehat{X}_Z$ for a qcqs derived scheme $X$ over $S$ equipped with a closed subscheme $Z \hookrightarrow X$. Write $i \colon U \coloneqq X \smallsetminus Z \hookrightarrow X$ for the inclusion of the open complementary and
\begin{align}
	\hat{\jmath} \colon \widehat{X}_Z \longrightarrow X 
\end{align}
for the induced morphism.
The presentation of $\catIndCoh(\widehat{X}_Z)$ given in \eqref{eq:IndCoh_colimit} provides a functor
\begin{align}
	\hat{\jmath}_\ast \colon \catIndCoh(\widehat{X}_Z / S) \longrightarrow \catIndCoh(X/S) \ . 
\end{align}
At least when $S = \Spec(k)$, the following result is well-known:

\begin{theorem}\label{thm:nilpotent_vs_set_theoretic}
	The functor $\hat{\jmath}_\ast$ is fully faithful and identifies $\catIndCoh(\widehat{X}_Z/S)$ with
	\begin{align}
		\catIndCoh(X/S)_Z \coloneqq \ker( i^\ast \colon \catIndCoh(X/S) \to \catIndCoh(U/S) ) \ . 
	\end{align}
\end{theorem}

\begin{remark}
	When $S = \Spec(k)$, this result is essentially the content of \cite[Proposition~7.4.5]{Gaitsgory_Rozenblyum_DG_indschemes}. In that case, the proof of the essential surjectivity part can be seen as an adaptation of \cite[Lemma~7.41]{Rouquier_Dimensions_of_triangulated}.
\end{remark}

Before giving the proof, we need the following two technical lemmas:
\begin{lemma}\label{lem:nilpotent_vs_theoretic_full_faithfulness}
	Let $\calX \to S$ be a relative qcqs derived formal scheme and set $Z \coloneqq \Red{\calX}$. The functor \eqref{eq:IndCoh_colimit_PrLomega} is fully faithful. Moreover, if $\calX \simeq \widehat{X}_Z$, the same statement holds with $X$ in place of $\calX$.
\end{lemma}

\begin{proof}
	We first prove the statement for $\calX$. Since $\scrT_{Z \overunder \calX}$ is filtered, it suffices to show that for $W \in \scrT_{Z \overunder \calX}$ and every pair $\calF, \calG \in \catCohb(W/S)$, the map
	\begin{align}
		\colim_{W' \in \scrT_{W \overunder \calX}} \Map_{\catIndCoh(W/S)}( j_{W,W',\ast}(\calF), j_{W,W',\ast}(\calG)) \longrightarrow \Map_{\catIndCoh(\calX/S)}( j_{W,\calX,\ast}(\calF), j_{W,\calX,\ast}(\calG) ) 
	\end{align}
	is an equivalence. Here $j_{W,W'} \colon W \to W'$ is the given transition morphism, which is a closed immersion. Using the adjunctions $j_{W,W',\ast} \dashv j_{W,W'}^!$ of Construction~\ref{construction:relative_IndCoh_operations}--\eqref{item:construction-relative_IndCoh_operations-1} and the fact that $\calF \in \catCohb(W/S) \simeq \catIndCoh(W/S)^\omega$, we can rewrite the source of this morphism as
	\begin{align}
		\Map_{\catIndCoh(W/S)}\Big( \calF, \colim_{W' \in \scrT_{W \overunder \calX}} j_{W,W'}^! j_{W,W',\ast}(\calG) \Big) \ . 
	\end{align}
	On the other hand, since the map $W \to \calX$ is representable by closed immersions by Corollary~\ref{cor:formal_scheme_diagonal_representable}, we also have the adjunction $j_{W,\calX,\ast} \dashv j_{W,\calX}^!$, which allows to rewrite the target of the above morphism as
	\begin{align}
		\Map_{\catIndCoh(W/S)}\big(\calF, j_{W,\calX}^! j_{W,\calX,\ast}(\calG)\big) \ . 
	\end{align}
	It therefore suffices to verify that the canonical comparison map
	\begin{align}
		\colim_{W' \in \scrT_{W \overunder \calX}} j_{W,W'}^! j_{W,W',\ast}(\calG) \longrightarrow j_{W,\calX}^! j_{W,\calX,\ast}(\calG)) 
	\end{align}
	is an equivalence in $\catIndCoh(W/S)$. Using base-change of Proposition~\ref{prop:relative_IndCoh_proper_basechange} and the formula of \cite[\S8.2]{Porta_Yu_Higher_Analytic_Stacks}, this reduces to check that the canonical map
	\begin{align}
		\colim_{W' \in \scrT_{W \overunder \calX}} W \times_{W'} W \longrightarrow W \times_{\calX} W 
	\end{align}
	is an equivalence. This follows from the descent property of the $\infty$-topos $\PreSt_{/S}$.
	
	When $X$ is in place of $\calX$ we see that the equivalence
	\begin{align}
		W \times_{\calX} W \simeq W \times_X W 
	\end{align}
	supplied by Lemma~\ref{lem:formal_completion_to_ambient_-1_truncated} implies the above proof still applies \textit{verbatim}.
\end{proof}

\begin{lemma}\label{lem:nilpotent_vs_theoretic_ess_surj}
	In the setting of Theorem~\ref{thm:nilpotent_vs_set_theoretic}, every object of
	\begin{align}
		\catCohb(X)_U \coloneqq \ker( i^\ast \colon \catCohb(X) \to \catCohb(U) ) 
	\end{align}
	is in the essential image of $\hat{\jmath}_\ast$.
\end{lemma}

\begin{proof}
	Given $\calF \in \ker(i^\ast)$ it suffices to construct $W \in \scrT_{Z \overunder \calX}$, $\widetilde{\calF} \in \catCohb(W/S)$ and an equivalence
	\begin{align}
		j_{W,\calX,\ast}(\widetilde{\calF}) \simeq \calF \ . 
	\end{align}
	Let $m$ be the number of non-vanishing homotopy groups of $\calF$. When $m = 0$, up to a shift we can assume that $\calF \in \catCoh^\heartsuit(\calX) \simeq \catCoh^\heartsuit(\trunc{\calX})$, in which case it follows immediately that a power of the sheaf of ideals defining $Z$ annihilates $\calF$ - and this is equivalent to find $W \in \scrT_{Z \overunder \calX}$, $\widetilde{\calF} \in \catCohb(W)$ and an equivalence
	\begin{align}
		j_{W,\calX,\ast}(\widetilde{\calF}) \simeq \calF \ . 
	\end{align}
	Recall once more from Corollary~\ref{cor:formal_scheme_diagonal_representable} that $j_{W,\calX}$ is representable by closed immersions.
	This implies that $j_{W,\calX,\ast} \colon \catQCoh(W) \to \catQCoh(\calX)$ satisfies the projection formula, and in turn this implies that since $\calF$ has finite tor-amplitude relative to $S$, the same goes for $\widetilde{\calF}$.
	
	As for the inductive step, we can assume, up to a shift, that $\calF$ is connective and that $m$ is the biggest integer such that $\pi_{m+1}(\calF) \ne 0$. Consider the fiber sequence
	\begin{align}
		\pi_{m+1}(\calF)[m+1] \longrightarrow \calF \longrightarrow \tau_{\leqslant m}(\calF) \ . 
	\end{align}
	Since $i^\ast$ is $t$-exact, the inductive hypothesis guarantees the existence of two thickenings $W_1$ and $W_2$ together with $\calF_i \in \catCohb(W_i/S)$ and identifications
	\begin{align}
		j_{W_1,\calX,\ast}(\calF_1) \simeq \pi_{m+1}(\calF)[m+1] \ , \qquad j_{W_2,\calX,\ast}(\calF_2) \simeq \tau_{\leqslant m}(\calF) \ . 
	\end{align}
	As before, this immediately implies that both $\calF_i$ have finite tor-amplitude relative to $S$. Since $\scrT_{Z \overunder \calX}$ is filtered, we can moreover assume without loss of generality that $W_1 = W_2$. Write $W$ instead of $W_1$ or $W_2$. The extension $\calF$ determines a canonical element in
	\begin{align}
		\pi_0 \Map_{\catIndCoh(\calX/S)}\big( \tau_{\leqslant m}(\calF), \pi_{m+1}(\calF)[m+2]\big) \simeq \pi_0 \Map_{\catIndCoh(\calX/S)}\big( j_{W,\calX,\ast}(\calF_1), j_{W,\calX,\ast}(\calF_2) \big) \ . 
	\end{align}
	On the other hand,  by the full faithfulness part, we have an identification
	\begin{align}
		\pi_0 \Map_{\catIndCoh(\calX)}\big( j_{W,\calX,\ast}(\calF_1), j_{W,\calX,\ast}(\calF_2) \big) \simeq \colim_{W' \in \scrT_{W \overunder \calX}} \pi_0 \Map_{\catIndCoh(W')}\big( j_{W, W',\ast}(\calF_1), j_{W,W',\ast}(\calF_2) \big) \ , 
	\end{align}
	where we used the fact that $\pi_0$ commutes with filtered colimits in $\Spc$. Thus, there exists a thickening $W'$ of $W$ such that the extension determined by $\calF$ is defined at the $W'$-level. The conclusion follows.
\end{proof}

\begin{proof}[Proof of Theorem~\ref{thm:nilpotent_vs_set_theoretic}]
	Combining Lemma~\ref{lem:nilpotent_vs_theoretic_ess_surj} and Proposition~\ref{prop:relative_IndCoh_localization}, we deduce that $\hat{\jmath}_\ast$ is essentially surjective. Moreover, Lemma~\ref{lem:nilpotent_vs_theoretic_full_faithfulness} implies that $\hat{\jmath}_\ast$ is fully faithful on the full subcategory of $\catIndCoh(\widehat{X}_Z)$ spanned by objects of the form
	\begin{align}
		\Big\{ j_{W,\widehat{X}_Z,\ast}(\calF) \: \Big| \: W \in \scrT_{Z \overunder \widehat{X}_Z}, \: \calF \in \catCohb(W/S)\Big\} \ . 
	\end{align}
	Notice that all these objects are compact in $\catIndCoh(\widehat{X}_Z / S)$ and that the same goes for their images in $\catIndCoh(X/S)$.
	Since $\hat{\jmath}_\ast$ is in $\PrL$ by construction, it suffices to argue that this collection of objects generates $\catIndCoh(\widehat{X}_Z)$.
	Let $\calF \in \catIndCoh(\widehat{X}_Z/S)$. Once again, to keep notations to a minimum, we will write $\calX \coloneqq \widehat{X}_Z$.
	Since
	\begin{align}
		\catIndCoh(\calX/S) \simeq \lim_{W\in \scrT_{Z \overunder \calX}\op} \catIndCoh(W/S) \ , 
	\end{align}
	we can represent $\calF$ as the compatible system $\{j_{W,\calX}^!(\calF)\}$. It is thus enough to argue that if
	\begin{align}
		\Map_{\catIndCoh(\calX/S)}( j_{Z,\calX,\ast}(\calG), \calF ) \simeq 0 
	\end{align}
	for all $\calG \in \catCohb(Z/S)$, then $\calF \simeq 0$. However,
	\begin{align}
		\Map_{\catIndCoh(\calX/S)}( j_{Z,\calX,\ast}(\calG), \calF ) \simeq \Map_{\catIndCoh(Z/S)}( \calG, j_{Z,\calX}^!( \calF ) ) \ , 
	\end{align}
	and since $\catIndCoh(Z / S) \simeq \Ind(\catCohb((Z/S))$, the conclusion follows.
\end{proof}

\begin{corollary}\label{cor:computing_Cohnil}
	In Formula~\eqref{eq:IndCoh_colimit}, the colimit can equally be computed in $\PrLomega$. Equivalently, the functor \eqref{eq:IndCoh_colimit_PrLomega} is an equivalence.
\end{corollary}

\begin{proof}
	To begin with, we observe that since $\scrT_{Z \overunder \calX}$ is filtered, the source of \eqref{eq:IndCoh_colimit_PrLomega} satisfies Zariski descent in $\calX$. The same goes for $\catIndCoh$, as shown in Proposition~\ref{prop:relative_IndCoh_Zariski_descent}. Applying \cite[Lemma~2.6]{Binda_Porta}, we see that the same goes for the target of \eqref{eq:IndCoh_colimit_PrLomega}. Moreover, this also shows that the statement is also Zariski local in $S$. Thus, it is enough to consider the case where $\calX$ is of the form $\widehat{X}_Z$. In this case, full faithfulness follows from Lemma~\ref{lem:nilpotent_vs_theoretic_full_faithfulness} and essential surjectivity from Lemma~\ref{lem:nilpotent_vs_theoretic_ess_surj}.
\end{proof}

\begin{definition}
	Let $S \in \evccndAff_k$ and let $\calX \to S$ be a relative qcqs derived formal scheme. A \textit{nilpotent bounded coherent complex} on $\calX$ is a compact object in $\mathsf{IndCoh}(\calX/S)$. We denote the resulting $\infty$-category as
	\begin{align}
		\catCohbnil(\calX/S) \coloneqq \catIndCoh(\calX/S)^\omega \ . \tag*{\qedhere}  
	\end{align}
\end{definition}

\begin{warning}	Let $\calX$ be a derived formal scheme. Then $\calX$ is in particular a derived stack, so $\catAPerf(\calX)$ is well defined and has a $t$-structure, but bounded objects in $\catAPerf(\calX)$ typically differs from $\catCohbnil(\calX)$. Consider for instance the following example: take $X = \A^1_k$, $Z = \Spec(k)$ embedded as $0$ inside $X$ and $\calX \coloneqq \widehat{X}_Z$.
	Then
	\begin{align}
		\catAPerf(\calX) \simeq \lim_{n} \catAPerf(k[T] / (T^n)) \simeq \catAPerf(k\llbracket T \rrbracket) \ ,
	\end{align}
	where the limit is computed with respect to pullbacks. It follows that bounded objects coincide with $\catCohb(k \llbracket T \rrbracket )$.
	On the other hand, we have:
	\begin{align}
		\catCohbnil(\calX) \simeq \colim_n \catCohbnil(k[T] / (T^n)) \ . 
	\end{align}
	Here the limit is computed with respect to the pullback, while the colimit is computed with respect to the pushforward. The ring of formal power series $k\llbracket T \rrbracket$ defines an object in $\catCohb(\calX)$ but not in $\catCohbnil(\calX)$.
\end{warning}

\begin{variant}\label{variant:nilpotent_qcoh_equivariant_case}
	Let $S$ be a derived stack and let $\formaldSch_S$ be the full subcategory of $\dSt_{/S}$ spanned by relative formal derived schemes.
	Using Lemma~\ref{lem:relative_IndCoh_base_change} and Proposition~\ref{prop:relative_IndCoh_Zariski_descent} we can define
	\[ \catIndCoh(-/S) \colon (\formaldSch_S)\op_{\mathsf{sch}} \longrightarrow \Cat_\infty \]
	by
	\begin{align}
		\catIndCoh(\calX/S) \coloneqq \lim_{T \in \dAff_{/S}} \catIndCoh(\calX \times_S T / T) \ . 
	\end{align}
	This applies in particular to the following situation: let $G$ be an algebraic group acting on a quasi-compact scheme $X$ and let $Z \hookrightarrow X$ be a closed subscheme which is $G$-invariant.
	Then the maps
	\begin{align}
		[X/G] \longrightarrow \sfB G \qquad \text{and} \qquad [Z/G] \simeq Z \times \sfB G \to \sfB G 
	\end{align}
	are schematic. The formal completion of $[X/G]$ along $[Z/G]$ coincides with $[\widehat{X}_Z/G]$, and the structural map to $\sfB G$ is representable by formal derived schemes. In particular, the $\infty$-category $\catIndCoh( [\widehat{X}_Z / G] / \sfB G )$ is well defined.
\end{variant}

\subsection{Derived moduli of nilpotent sheaves}

We now adapt the definition of the derived stack of coherent sheaves in \cite[\S2.1]{Porta_Sala_Hall} to the formal setting.
\begin{definition}[Flat nilpotent families]
	Let $\calX$ be a qcqs derived formal scheme and fix an eventually coconnective derived affine scheme $S \in \evccndAff_k$.
	\begin{itemize}\itemsep=0.2cm
		\item An \textit{$S$-family of nilpotent coherent complexes} is an object $\calF \in \catCohbnil(\calX \times S / S)$.
		\item We say that an $S$-family of nilpotent coherent complexes $\calF$ is \textit{$S$-flat} (resp.\ \textit{$S$-properly supported}) if there exists a laft thickening $W$ of $\Red{\calX}$ inside $\calX$, $\overline{\calF} \in \catCohb(W \times S / S)$ an isomorphism $j_{W \times S,\calX \times S,\ast}(\overline{\calF}) \simeq \calF$ and moreover $\widetilde{\calF}$ is flat over $S$ in the sense of \cite[Definition~2.1]{Porta_Sala_Hall} (resp.\ proper over $S$ in the sense of \cite[Definition~2.27]{Porta_Sala_Hall})
		\item An \textit{$S$-family of nilpotent coherent sheaves} is an $S$-family of nilpotent coherent complexes which is furthermore $S$-flat and $S$-properly supported.
	\end{itemize}
	We write $\catCohnil_{S,\ps}( \calX \times S )$ for the full subcategory of $\catCohbnil(\calX \times S / S)$ spanned by $S$-families of nilpotent coherent sheaves.
\end{definition}

\begin{lemma}\label{lem:pushforward_coherent}
	Let $S \in \evccndAff_k$ and let $f \colon Y_1 \to Y_2$ be a finite morphism between quasi-compact derived schemes over $S$. Then, an object $\calF \in \catCohb(Y_1)$ belongs to $\catCoh_{S,\ps}(Y_1)$ if and only if $f_{S,\ast}(\calF)$ belongs to $\catCoh_{S,\ps}(Y_2)$. Thus, the functor $f_{S,\ast} \colon \catIndCoh(Y_1 / S) \to \catIndCoh(Y_2 / S)$ induces a well defined functor
	\begin{align}
		\catCoh_{S,\ps}(Y_1) \longrightarrow \catCoh_{S,\ps}(Y_2) \ . 
	\end{align}
\end{lemma}

\begin{proof}
	Since $f$ is finite, in particular proper, it takes $\catAPerf(Y_1 \times S)$ to $\catAPerf(Y_2 \times S)$, see \cite[Theorem~5.6.0.2]{Lurie_SAG}. Besides, the image of a subset of $Y_1 \times S$ proper over $S$ via $f_S$ will still be proper over $S$. Therefore, $f_{S,\ast}$ preserves the properly supported condition. It also reflects it, because since $f$ is proper, the preimage via $f$ of a proper subset of $Y_2 \times S$ will be proper in $Y_1 \times S$. We are left to show that it preserves and reflects $S$-flatness. Let $\calF \in \catAPerf(Y_1 \times S)$ be $S$-flat and let $\calG\in \catQCoh^\heartsuit(Y_2 \times S)$. Write
	\begin{align}
		q_S \colon Y_1 \times S \longrightarrow S \qquad \text{and} \qquad p_S \colon Y_2 \times S \longrightarrow S 
	\end{align}
	for the canonical projections. Then, the projection formula yields
	\begin{align}
		p_S^\ast(\calG) \otimes_{\scrO_{X_S}} f_{S,\ast}(\calF) \simeq f_{S,\ast}( f^\ast p_S^\ast(\calG) \otimes_{\scrO_{Y_S}} \calF ) \simeq f_\ast( q_S^\ast(\calG) \otimes_{\scrO_{Y_S}} \calF ) \ . 
	\end{align}
	Recall that $f$ being finite guarantees that it is affine as well. Then, the conclusion follows from the fact that $f_{S,\ast}$ is $t$-exact and conservative.
\end{proof}

\begin{remark}
	Assume that $\calX \simeq \widehat{X}_Z$ and fix $S \in \evccndAff_k$. Given $\calF \in \catCohbnil(\widehat{X}_Z \times S / S)$, choose a representative $\widetilde{\calF}$ defined over a thickening $W \times S$ of $Z \times S$ inside $X \times Z$. Write
	\begin{align}
		 j_{W,X} \colon W \times S \longrightarrow X \times S \qquad \text{and} \qquad \hat{\jmath}_{W,\calX} \colon W \times S \longrightarrow \calX \times S 
	\end{align}
	for the corresponding closed immersions (we omit $S$ in the notation for an increased readability). Then $\calF$ is $S$-flat if and only if $j_{W,X,\ast}(\widetilde{\calF})$ is $S$-flat. Moreover, flatness can also be formulated intrinsically in terms of the $t$-structure on $\catIndCoh(\calX/S)$; at least when $X$ is smooth, the two definitions coincide thanks to Lemma~\ref{lem:relative_IndCoh_relative_Serre_regularity}.
\end{remark}

\begin{construction}
	Fix a quasi-compact derived formal scheme $\calX$. The functoriality of $\catCohb(-/S)$ in $S$ discussed in Proposition~\ref{lem:relative_IndCoh_base_change} together with \cite[Lemma 2.4]{Porta_Sala_Hall} imply that families of nilpotent coherent sheaves assemble into a well-defined functor
	\begin{align}
		\widetilde{\dstackCoh}_{\ps}^{\mathsf{nil}}(\calX) \colon \evccndAff_k\op \longrightarrow \Spc 
	\end{align}
	by the rule
	\begin{align}
		\widetilde{\dstackCoh}_{\ps}^{\mathsf{nil}}(\calX)(S) \coloneqq \catCohnil_{S,\ps}(\calX \times S)^\simeq \ , 
	\end{align}
	where $(-)^\simeq$ denotes the maximal $\infty$-groupoid. Finally, we let
	\begin{align}
		\dstackCohpsnil(\calX) \coloneqq \evccnj_\ast( \widetilde{\dstackCoh}_{\ps}^{\mathsf{nil}}(\calX) )
	\end{align}
	be the right Kan extension of $\widetilde{\dstackCoh}^{\mathsf{nil}}_{\ps}(\calX)$ along $\evccnj \colon \evccndAff_k \hookrightarrow \dAff_k$. We write $\stackCohpsnil(\calX)$ for its truncation.
\end{construction}

\begin{remark}\label{rem:Cohnil_convergent}
	It follows from Recollection~\ref{recollection:convergence} that $\dstackCohnil_{\ps}(\calX)$ is convergent.
\end{remark}

\begin{example}\label{eg:Cohnil_indgeometric_presentation_I}
	Let $\calX \simeq \widehat{X}_Z$ for $X$ a quasi-compact derived scheme and let $S \in \dAff_k$ be an affine derived scheme. Then combining Recollection~\ref{recollection:convergence} and Corollary~\ref{cor:computing_Cohnil}, we find
	\begin{align}
		\dstackCohpsnil(\calX)(S) & \simeq \lim_n \colim_{W \in \scrT_{Z \overunder X}} \dstackCohps(W)(\sft_{\leqslant n}(S)) \\
		& \simeq \lim_{n} \colim_{W \in \scrT_{Z \overunder X}} \catCoh_{S,\ps}(W \times \sft_{\leqslant n}(S))^\simeq \ .
	\end{align}
	Notice that when computing the truncation $\stackCohpsnil(\calX)$, the exterior limit disappears.
\end{example}

\subsection{Indgeometricity}

Our first goal is to show that for formal derived schemes of the form $\widehat{X}_Z$, the derived prestack $\dstackCohpsnil(\widehat{X}_Z)$ is indgeometric under some mild assumptions on $X$ (see Theorem~\ref{thm:indgeometricity} below). Before giving the precise statement, we need to introduce the following:

\begin{definition}\label{def:perfect_resolution_property}
	Let $X$ be a derived scheme. We say that $X$ \textit{satisfies the perfect resolution property} if for every $\calF \in \catCoh^\heartsuit(X)$ there exists $\calE \in \catPerf(X)_{\geqslant 0}$ and a morphism
	\begin{align}
		\calE \longrightarrow \calF 
	\end{align}
	inducing a surjection on $\pi_0$. We say that $X$ \textit{universally satisfies the perfect resolution property} if for every derived affine scheme $S$, $X \times S$ satisfies the perfect resolution property.
\end{definition}

\begin{example}\label{eg:quasi_projective_implies_perfect_resolution}
	If $X$ is underived and it has the resolution property (i.e., for every classical coherent sheaf $\calF$ there exists an epimorphism $\calE \to \calF$, with $\calE$ being a vector bundle), then it also satisfies the perfect resolution property. In particular, if $X$ is quasi-projective, then it universally satisfies the perfect resolution property.
\end{example}

\begin{lemma}\label{lem:perfect_resolution_property_hereditary}
	Let $j \colon X \hookrightarrow \overline{X}$ be a closed immersion of derived schemes. If $\overline{X}$ (universally) satisfies the perfect resolution property, then the same goes for $X$.
\end{lemma}

\begin{proof}
	We treat the universal case; the other follows by the same proof below taking $S = \Spec(k)$.
	
	\medskip
	
	Fix an affine derived scheme $S$ and let $\calF \in \catCoh^\heartsuit(X \times S)$. Then $j_S \colon X \times S \to \overline{X} \times S$ is again a closed immersion and $\overline{X} \times S$ satisfies the perfect resolution property. Let $\calF \in \catCoh^\heartsuit(X \times S)$.
	Then $j_{S,\ast}(\calF) \in \catCoh^\heartsuit(\overline{X} \times S)$. It follows that there exists $\calE\in \catPerf(\overline{X} \times S)_{\geqslant 0}$ and a morphism
	\begin{align}
		\calE \longrightarrow j_{S,\ast}(\calF) 
	\end{align}
	being surjective on $\pi_0$. Since $j_{S}$ is a closed immersion, it follows that the composite
	\begin{align}
		j_S^\ast(\calE) \longrightarrow j_S^\ast(j_{S,\ast}(\calF)) \longrightarrow \calF 
	\end{align}
	is also surjective on $\pi_0$.
\end{proof}

\begin{corollary}\label{cor:smooth_implies_universal_perfect_resolution}
	Smooth schemes universally satisfy the perfect resolution property.
\end{corollary}

\begin{proof}
	Let $Y$ be a smooth scheme. Note that every object in $\catCoh^\heartsuit(Y)$ is perfect. Thus, $Y$ satisfies the perfect resolution property. 
	
	Let $S$ be a derived affine scheme and let $X$ be a smooth scheme. Since $S$ is affine, there exists a closed immersion in $\A^n_k$ for some $n \geqslant 0$. Thus, the induced morphism $X \times S \to X \times \A^n_k$ is a closed immersion. Since $Y\coloneqq X \times \A^n_k$ is again smooth, $Y$ satisfies the perfect resolution property. Therefore, the conclusion follows from Lemma~\ref{lem:perfect_resolution_property_hereditary}.
\end{proof}

We can now formulate the main result of this section.
\begin{definition}
	We say that a qcqs derived formal scheme $\calX$ \textit{universally satisfies the perfect resolution property} if the full subcategory of $\scrT_{\Red{\calX} \overunder \calX}$ spanned by those thickenings $W$ that universally satisfy the perfect resolution property is cofinal.
\end{definition}

\begin{example}
	If $\calX \simeq \widehat{X}_Z$, with $X$ either a quasi-projective or a smooth scheme, then Lemma~\ref{lem:perfect_resolution_property_hereditary} shows that $\calX$ universally satisfies the perfect resolution property.
\end{example}

\begin{theorem}\label{thm:indgeometricity}
	Let $\calX$ be a qcqs derived formal scheme universally satisfying the perfect resolution property.
	Then, the derived stack $\dstackCohpsnil(\calX)$ is indgeometric.
\end{theorem}

Before giving the proof of Theorem~\ref{thm:indgeometricity}, we need some technical preliminaries. The key statement is Proposition~\ref{prop:perfect_resolution_implies_closed_immersion_of_stacks} below. We start with two general constructions.

\begin{construction}\label{construction:modular_pushforward}
	Let $f \colon Y \to Y$ be a closed immersion between quasi-compact and quasi-separated schemes locally almost of finite presentation over $k$. Notice that given a morphism $T \to S$ of derived affine schemes, the square
	\begin{align}
		\begin{tikzcd}[ampersand replacement=\&]
			Y \times T \arrow{r} \arrow{d}{f_{T}} \& Y \times S \arrow{d}{f_S} \\
			X \times T \arrow{r} \& X \times S
		\end{tikzcd}
	\end{align}
	is a derived pullback. In particular, derived base-change (see e.g.\ \cite[Proposition~2.5.4.5]{Lurie_SAG}) and Lemma~\ref{lem:pushforward_coherent} imply that $f$ induces a well-defined morphism
	\begin{align}\label{eq:pushforward_derived_stacks_coherent}
		\bff_\ast \colon \dstackCohps(Y) \longrightarrow \dstackCohps(X) \ .
	\end{align}
	Our goal is to establish that, under mild assumptions on $X$, $\bff_\ast$ is itself a closed immersion.
\end{construction}

\begin{construction}\label{construction:zero_locus_morphism_of_sheaves}
	Let $p \colon X \to S$ be a flat morphism of derived schemes, with $S$ being affine. Given $f \colon T \to S$, we write $X_T \coloneqq T \times_S X$ and we let $f_X \colon X_T \to X$ be the naturally induced morphism. Let $\phi \colon \calG \to \calF$ be a morphism in $\catAPerf(X)$.	We define the functor of points
	\begin{align}
		\sfZ_\phi \colon \dAff_{/S}\op \longrightarrow \Spc
	\end{align}
	sending $f \colon T \to S$ to
	\begin{align}
		\Map_{\Map_{\catQCoh(X_T)}(f_X^\ast(\calG), f_X^\ast(\calF))}( f_X^\ast(\phi), 0 ) \ , 
	\end{align}
	i.e., the space of nullhomotopies of $f_X^\ast(\phi) \colon f_X^\ast(\calG) \to f_X^\ast(\calF)$ in $\catQCoh(X_T)$. Notice that there is a canonical natural transformation $i \colon \sfZ_\phi \to S$.
\end{construction}

\begin{remark}\label{rem:truncation_of_zero_locus_of_a_morphism}
	In the setting of Construction~\ref{construction:zero_locus_morphism_of_sheaves}, assume that both $\calF$ and $\calG$ are connective and that $\calF$ is $S$-flat. For every \textit{underived} affine scheme $T$, $f_X^\ast(\calF)$ is $T$-flat and hence it belongs to $\catQCoh^\heartsuit(X_T)$. Since $f_X^\ast(\calG)$ is connective, it follows that
	\begin{align}
		\Map_{\catQCoh(X_T)}\big(f_X^\ast(\calG), f_X^\ast(\calF)\big) \simeq \Map_{\catQCoh^\heartsuit(X_T)}\big( \pi_0(f_X^\ast(\calG)), f_X^\ast(\calF) \big) \ . 
	\end{align}
	Under this equivalence, we see that $\Map_{\catQCoh^\heartsuit(X_T)}\big( \pi_0(f_X^\ast(\calG)), f_X^\ast(\calF) \big)$ is a \textit{set}, and $\sfZ_\phi(T)$ can be explicitly characterized as follows:
	\begin{align}
		\sfZ_\phi(T) = \begin{cases}
			\ast & \text{if } \pi_0( f_X^\ast(\phi) ) = 0 \ ,\\
			\emptyset & \text{otherwise}\ .
		\end{cases} 
	\end{align}
	Said differently, $\trunc{\sfZ_\phi}$ is identified with the subfunctor of $\Map_{\dSt_k}(-,S)$ that sends an underived affine scheme $T$ to the collection of morphisms $f \colon T \to S$ such that $\pi_0( f_X^\ast(\phi) ) = 0$.
\end{remark}

The following lemma is a generalization of \cite[Theorem~5.8 and Remark 5.9]{FGA_Explained} to the context of derived geometry.
\begin{lemma}\label{lem:zero_locus_morphism_of_sheaves}
	In the setting of Construction~\ref{construction:zero_locus_morphism_of_sheaves}, assume that $\calG$ is connective and that $\calF$ is $S$-flat.
	Then:
	\begin{enumerate}\itemsep=0.2cm
		\item \label{item:zero_locus_morphism_of_sheaves-1} let $\alpha \colon \calE \to \calG$ be a morphism with $\calE$ connective and let
		\begin{align}
			\psi \colon \calE \longrightarrow \calF 
		\end{align}
		be the morphism obtained by composition with $\phi$. If $\alpha$ induces an epimorphism on $\pi_0$, then the induced natural transformation
		\begin{align}
			\sfZ_\phi \longrightarrow \sfZ_\psi 
		\end{align}
		is an equivalence after passing to the truncations.
		
		\item \label{item:zero_locus_morphism_of_sheaves-2} Assume furthermore that $\calG$ is perfect and that $\calF$ is $S$-properly supported. Then, $\sfZ_\phi$ is a derived scheme and the morphism $i \colon \sfZ_\phi \to S$ is a closed immersion.
	\end{enumerate}
\end{lemma}

\begin{proof}
	Statement~\eqref{item:zero_locus_morphism_of_sheaves-1} follows directly from the definitions and from Remark~\ref{rem:truncation_of_zero_locus_of_a_morphism}.
	
	Let us prove Statement~\eqref{item:zero_locus_morphism_of_sheaves-2}. Consider $\calHom_X(\calG, \calF) \simeq \calG^\vee \otimes_{\scrO_X} \calF$. Notice that for every $f \colon T \to S$ and every $M \in \catQCoh^\heartsuit(T)$, the complex
	\begin{align}
		M \otimes_{\scrO_T} p_{T,\ast} \calHom_{X_T}( f_X^\ast(\calG), f_X^\ast(\calF) ) \simeq p_{T,\ast}\big( \calHom_{X_T}( f_X^\ast(\calG), p_T^\ast(M) \otimes f_X^\ast(\calF) ) \big) 
	\end{align}
	is in negative homological degrees (where $p_T \colon X_T \coloneqq X \times_S T \to T$ and $f_X \colon X_T \to X$ are the naturally induced morphisms). Thus,
	\begin{align}
		\calH \coloneqq p_\ast \calHom_X(\calG, \calF) 
	\end{align}
	has tor-amplitude $\leqslant 0$ on $S$. Since $\calF$ is $S$-properly supported, we have $\calH \in \catAPerf(S)$ and therefore $\calH$ is perfect. Its dual $\calH^\vee$ is connective, so it follows that
	\begin{align}
		V \coloneqq \Spec_S( \Sym_{\scrO_S}( \calH^\vee ) ) 
	\end{align}
	is a derived scheme. For every $f \colon T = \Spec(A) \to S$, we have
	\begin{align}
		\Map_{/S}(T, V) & \simeq \Map_{\catQCoh(T)}( f^\ast( \calH^\vee ), A ) \\
		& \simeq \Map_{\catQCoh(T)}( A, f^\ast(\calH) ) \\
		& \simeq \Map_{\catQCoh(T)}( A, p_{T,\ast}(\calHom_{X_T}(f_X^\ast(\calG), f_X^\ast(\calF))) ) \ .
	\end{align}
	It follows that the morphism $\phi \colon \calG \to \calF$ is classified by a morphism $s_\phi \colon S \to V$. Write $s_0 \colon S \to V$ for the zero section, which is clearly a closed immersion. Unraveling the definitions, we see that
	\begin{align}
		\begin{tikzcd}[ampersand replacement=\&]
			Z_\phi \arrow{r}{i} \arrow{d}{i} \& S \arrow{d}{s_\phi} \\
			S \arrow{r}{s_0} \& V
		\end{tikzcd} 
	\end{align}
	is a derived pullback square, whence the conclusion.
\end{proof}

\begin{proposition}\label{prop:perfect_resolution_implies_closed_immersion_of_stacks}
	Assume that $f \colon Y \to X$ is a closed immersion and assume that $X$ universally satisfies the perfect resolution property. Then
	\begin{align}
		\bff_\ast \colon \dstackCohps(Y) \longrightarrow \dstackCohps(X) 
	\end{align}
	is a closed immersion as well.
\end{proposition}

\begin{proof}
	First observe that \cite[\S2.3.1]{Porta_Sala_Hall} (or \cite[Theorem~5.2.2]{HL_Categorical_properness}) shows that both source and target of $\bff_\ast$ are derived $1$-geometric stacks, locally almost of finite presentation (see also \cite[Tag~0DLX]{stacks-project} for the analogous statement on truncations). Thus, it is enough to check the statement for the induced morphism between the classical truncations, which we still denote
	\begin{align}
		\bff_\ast \colon \stackCohps(Y) \to \stackCohps(X) \ . 
	\end{align}
	Let $S$ be an underived affine scheme and let $x \colon S \to \stackCohps(X)$ be a point classifying a family of $S$-properly supported, $S$-flat coherent sheaves $\calF \in \catAPerf(X \times S)$. Since $S$ is underived, $S$-flatness implies that $\calF \in \catQCoh^\heartsuit(X \times S)$. Set
	\begin{align}
		\calI_S \coloneqq \fib\big( \scrO_{X \times S} \to f_{S,\ast}(\scrO_{Y \times S}) \big) \ . 
	\end{align}
	Then $x$ belongs to the image of $f_\ast$ if and only if the map
	\begin{align}
		\phi \colon \pi_0( \calI \otimes \calF ) \longrightarrow \calF 
	\end{align}
	is zero. It follows that
	\begin{align}
		S \times_{\stackCohps(X)} \stackCohps(Y) \simeq \sfZ_\phi \ . 
	\end{align}
	Since $X$ universally satisfies the perfect resolution property, we can find a connective perfect complex $\calE$ and a morphism
	\begin{align}
		\calE \longrightarrow \pi_0( \calI \otimes \calF ) 
	\end{align}
	which induces an epimorphism on $\pi_0$. At this point, the conclusion follows from Lemma~\ref{lem:zero_locus_morphism_of_sheaves}.
\end{proof}

\begin{proof}[Proof of Theorem~\ref{thm:indgeometricity}]
	The derived stack $\dstackCohpsnil(\calX)$ is convergent by Remark~\ref{rem:Cohnil_convergent}. It is therefore enough to construct a presentation. Notice that Example~\ref{eg:Cohnil_indgeometric_presentation_I} implies that for every $S \in \evccndAff_k$ there is an equivalence
	\begin{align}
		\colim_{W \in \scrT_{Z \overunder X}} \dstackCohps(W)(S) \stackrel{\sim}{\longrightarrow} \dstackCohpsnil(\widehat{X}_Z)(S) \ . 
	\end{align}
	Besides, it follows from \cite[\S2.3.1]{Porta_Sala_Hall} or \cite[Theorem~5.2.2]{HL_Categorical_properness} that $\dstackCohps(W)$ is a derived $1$-geometric stack, locally almost of finite presentation (see also \cite[Tag~0DLX]{stacks-project} for the analogous statement on truncations). This shows that condition \eqref{item:indgeometric_stack-2-c} in Definition~\ref{def:indgeometric_stack} is satisfied. We are therefore left to check condition \eqref{item:indgeometric_stack-2-b}; in other words, we have to check that for a morphism $W \to W'$ in $\scrT_{Z\overunder X}$ the induced morphism
	\begin{align}
		\dstackCohpsnil(W) \longrightarrow \dstackCohps(W') 
	\end{align}
	is a closed immersion almost of finite type. Since both source and target are locally almost of finite presentation, the above morphism is automatically almost of finite type. By assumption, we can assume that $W'$ universally satisfies the perfect resolution property; thus Lemma~\ref{lem:perfect_resolution_property_hereditary} guarantees that the same goes for $W'$. At this point, the conclusion follows from Proposition~\ref{prop:perfect_resolution_implies_closed_immersion_of_stacks}.
\end{proof}

Theorem~\ref{thm:indgeometricity} and Proposition~\ref{prop:indgeometric_descent} yield the following.
\begin{corollary}
	Under the assumptions of Theorem~\ref{thm:indgeometricity}, the derived prestack $\dstackCohpsnil(\calX)$ is a derived stack, i.e., it satisfies \'etale hyperdescent.
\end{corollary}

\subsection{Nilpotent versus set-theoretic support}\label{subsec:nilpotent-vs-set-theoretic}

We keep fixing a field of characteristic zero $k$. We also fix a qcqs derived $k$-scheme $X$ together with a closed immersion $j \colon Z \to X$.
The goal of this section is to provide an alternative description of $\dstackCohpsnil(\widehat{X}_Z)$ in terms of $\dstackCohps(X)$.
This alternative description is particularly well suited for explicit computations, but it is not intrinsic in the formal completion $\widehat{X}_Z$ and relies on the extra knowledge of the ambient space $X$.

Denote by $i \colon U \to X$ the inclusion of the open complementary of $Z$.
\begin{definition}\label{def:families_of_nilpotent_sheaves}
	Let $S \in \dAff_k$ be an affine derived scheme. The $\infty$-category of \textit{$S$-families of quasi-coherent sheaves on $X$ set-theoretically supported on $Z$} is the kernel
	\begin{align}
		\catQCoh_Z(X \times S) \coloneqq \ker\big( i_S^\ast \colon \catQCoh(X \times S) \to \catQCoh(U \times S) \big) \ . 
	\end{align}
	Moreover, given $\calF \in \catQCoh_Z(X \times S)$ we say that $\calF$ is \textit{almost perfect} (resp.\ \textit{$S$-flat}, \textit{$S$-properly supported}) if it is almost perfect (resp.\ $S$-flat, $S$-properly supported) as an object in $\catQCoh(X \times S)$.
\end{definition}

It is straightforward to check that the above definition cuts a derived substack
\begin{align}
	\dstackCoh_{Z,\ps}(X) \hookrightarrow \dstackCohps(X) \ , 
\end{align}
that parametrizes $S$-flat families of $S$-properly supported almost perfect sheaves on $X$ set-theoretically supported on $Z$.
By restricting and right Kan extending along $\evccnj \colon \evccndAff_k \hookrightarrow \dAff_k$, we see that Corollary~\ref{cor:computing_Cohnil}, together with Lemma~\ref{lem:pushforward_coherent}, provides a well defined morphism of derived stacks
\begin{align}
	\bfjmathhat_\ast \colon \dstackCohpsnil(\widehat{X}_Z) \longrightarrow \dstackCoh_{Z,\ps}(X) \ . 
\end{align}
We can now state the main result of this section.

\begin{theorem}\label{thm:nilpotent_into_set_theoretic_moduli}
	The morphism $\bfjmathhat_\ast$ defined above is an equivalence of derived stacks.
\end{theorem}

\begin{proof}
	We first notice that the derived stack $\dstackCoh_{Z,\ps}(X)$ fits in the following pullback square:
	\begin{align}\label{eq:set_theoretic_pullback_description}
		\begin{tikzcd}[ampersand replacement=\&]
			\dstackCoh_{Z,\ps}(X) \arrow{r} \arrow{d} \& \dstackCohps(X) \arrow{d} \\
			\Spec(k) \arrow{r}{0} \& \dstackCoh(U) \ ,
		\end{tikzcd}
	\end{align}
	where the morphism $0 \colon \Spec(k) \to \dstackCoh(U)$ selects the zero sheaf. It immediately follows that $\dstackCoh_{Z,\ps}(X)$ is convergent. Since $\dstackCohpsnil(\widehat{X}_Z)$ is convergent as well, it is enough to prove that for every $S \in \evccndAff_k$, the morphism $\bfjmathhat_\ast$ induces an equivalence
	\begin{align}
		\bfjmathhat_\ast \colon \dstackCohpsnil(\widehat{X}_Z)(S) \longrightarrow \dstackCoh_{Z,\ps}(X)(S) \ . 
	\end{align}
	For this, it is enough to prove that this morphism is fully faithful and essentially surjective. Full faithfulness follows from Theorem~\ref{thm:nilpotent_vs_set_theoretic}, while essential surjectivity follows combining the same theorem with Lemma~\ref{lem:pushforward_coherent}.
\end{proof}

Composing the top horizontal arrow in the diagram \eqref{eq:set_theoretic_pullback_description} with $\bfjmathhat_\ast$ we obtain a canonical morphism
\begin{align}\label{eq:nilpotent_into_ambient}
	\dstackCohpsnil(\widehat{X}_Z) \longrightarrow \dstackCohps(X) \ .
\end{align}
We have the following.
\begin{corollary}\label{cor:Cohnil_as_formal_completion-1}
	The morphism \eqref{eq:nilpotent_into_ambient} is formally étale and it exhibits $\dstackCohpsnil(\widehat{X}_Z)$ as the formal completion of $\dstackCohps(X)$ along $\Red{\dstackCohpsnil(\widehat{X}_Z)}$. In other words, the square
	\begin{align}
		\begin{tikzcd}[ampersand replacement=\&]
			\dstackCohpsnil(\widehat{X}_Z) \arrow{r} \arrow{r} \arrow{d} \& \dstackCohps(X) \arrow{d} \\
			\big(\Red{\dstackCohpsnil(\widehat{X}_Z)}\big)_{\mathsf{dR}} \arrow{r} \& \dstackCohps(X)_{\mathsf{dR}}
		\end{tikzcd}
	\end{align}
	is a pullback.
\end{corollary}

\begin{proof}
	We start by the first statement. Since $\bfjmathhat_\ast$ is an equivalence by Theorem~\ref{thm:nilpotent_into_set_theoretic_moduli}, it is enough to prove that the top horizontal morphism in the square \eqref{eq:set_theoretic_pullback_description} is formally étale. Since that square is a pullback, it suffices to show that the bottom horizontal arrow is formally étale, which is indeed true as a consequence of the Nakayama lemma . This prove the first statement. As for the second one, in virtue of Lemma~\ref{lem:recognition_criterion_formal_completions} and what we just proved, it suffices to argue that both $\dstackCohpsnil(\widehat{X}_Z)$ and $\dstackCoh(X)$ are convergent and laft. For the latter, both properties are well known; on the other hand the former is convergent by definition. Concerning it being laft, observe that the pullback description \eqref{eq:set_theoretic_pullback_description} implies that $\dstackCoh_{Z,\ps}(X)$ is the pullback of three laft stack, and that it is therefore laft itself. At this point, the conclusion follows from Theorem~\ref{thm:nilpotent_into_set_theoretic_moduli}.
\end{proof}

\begin{remark}
	It is natural to wonder whether $\Red{\dstackCohpsnil(\widehat{X}_Z)}$ is actually a closed substack of $\dstackCohps(X)$.
	This is is not entirely obvious, but it is indeed true, at least under some mild additional assumptions on $X$.
	See Corollary~\ref{cor:Cohnil_as_formal_completion-2}.
\end{remark}

\subsection{Admissibility}\label{subsec:admissibility}

Let $X$ be a quasi-compact scheme and let $Z \hookrightarrow X$ be a closed immersion. We saw in Theorem~\ref{thm:indgeometricity} that mild assumptions on $X$ imply that $\dstackCohpsnil(\widehat{X}_Z)$ is indgeometric. The goal of this section is to prove that under slightly stronger assumptions, $\dstackCohpsnil(\widehat{X}_Z)$ is admissible (see Definition~\ref{def:admissible_indgeometric_stack}). We will achieve this by constructing an explicit admissible open exhaustion of $\dstackCohpsnil(\widehat{X}_Z)$ given by \textit{Harder-Narasimhan strata}.

\medskip

To begin with, fix a thickening $W \in \scrT_{Z \overunder X}$. Construction~\ref{construction:modular_pushforward} yields a well-defined morphism
\begin{align}
	\mathbf{j}_{W,\ast} \colon \dstackCohps(W) \longrightarrow \dstackCohps(X) \ , 
\end{align}
and the universal property of colimits and of the convergent completion (see Recollection~\ref{recollection:convergence}) produce a canonical morphism
\begin{align}\label{eq:nilpotent_embedding_into_ambient}
	\dstackCohpsnil(\widehat{X}_Z) \longrightarrow \dstackCohps(X) \ .
\end{align}
We will use it to induce the desired admissible open exhaustion of $\dstackCohpsnil(\widehat{X}_Z)$, out of the geometry of $\dstackCohps(X)$.

To begin our construction, we start assuming that $X$ is a projective scheme over a field $k$. This condition will be relaxed later on, see Theorem~\ref{thm:Cohnil_admissible}.

We denote by $H$ an ample divisor on $X$. Given any coherent sheaf $\calE$ on $X$, we define as usual the \textit{Hilbert polynomial $P_H(\calE, t)$ of $\calE$} as
\begin{align}
	P_H(\calE, t)\coloneqq \chi(E\otimes\scrO_X(tH))\in \Q[m]\ .
\end{align}
By \cite[Lemma~1.2.1]{HL_Moduli}, $P_H(\calE, t)$ can be uniquely written in the form
\begin{align}
	P_H(\calE, t)=\sum_{i=0}^{\dim(\calE)}c_i(\calE) \frac{t^i}{i!}\ .
\end{align}
Following \cite[\S1]{Moduli_pi_one}, when $\dim(\calE)\geq 1$, the \textit{reduced Hilbert polynomial of $\calE$} is
\begin{align}
	p_H(\calE, t)\coloneqq \frac{P_H(\calE, t)}{c_{\dim(\calE)}(\calE)}\ .
\end{align}
Recall moreover that there is a natural ordering of polynomials given by the lexicographic order of their coefficients. Explicitly, $P (\leq) P'$ if and only if $P(m) (\leq) P'(m)$ for $m \gg 0$.\footnote{Here, we are using the convention in \cite[Notation~1.2.5]{HL_Moduli}.} We can now introduce a notion of semistability.
\begin{definition}
	A coherent sheaf $\calE$ of dimension $m\geq 1$ on $X$ is (semi)stable if it is pure and if for all proper non-trivial subsheaves $\calE'\subsetneq \calE$ one has $p_H(\calE')(\leq)p_H(\calE)$.
\end{definition}
If $\calE$ is a coherent sheaf of dimension $m\geq 1$, it admits a unique filtration in (the \textit{Harder-Narasimhan filtration})
\begin{align}
	0\subset \calE_0 \subset \calE_1 \subset \cdots \subset \calE_\ell \coloneqq \calE 
\end{align}
where $\calE_0\subset \calE$ is the maximal $(m-1)$-dimensional subsheaf of $\calE$ and the factors $\calE_i/\calE_{i-1}$ are semistable and their reduced Hilbert polynomial satisfy
\begin{align}
	p_H(\calE_1/\calE_0) > \cdots > p_H(\calE/\calE_{\ell-1})\ .
\end{align}
Employing standard terminology, the \textit{minimal Hilbert polynomial of $\calE$} is defined by
\begin{align}
	p_{H\textrm{-}\min}(\calE) \coloneqq p_H(\calE/\calE_{\ell-1})\ . 
\end{align}
If $\calE_0=0$, \textit{the maximal Hilbert polynomial of $\calE$} is defined by 
\begin{align}
	p_{H\textrm{-}\max}(\calE) \coloneqq p_H(\calE_1)\ . 
\end{align}
For any fixed polynomial $P(t) \in \Q[t]$, we let $\dstackCohps(X;P(t))$ be the open and closed substack of $\dstackCohps(X)$ parametrizing flat and properly supported families of coherent sheaves having $P(t)$ as their Hilbert polynomial. We also set
\begin{align}
	\dstackCohpsnil(\widehat{X}_Z;P(t)) \coloneqq \dstackCohpsnil(\widehat{X}_Z) \times_{\dstackCohps(X)} \dstackCohps(X;P(t)) \ . 
\end{align}
Similarly, for every thickening $W \in \scrT_{Z \overunder X}$, we set
\begin{align}
	\dstackCohps(W;P(t)) \coloneqq \dstackCohps(W) \times_{\dstackCohps(X)} \dstackCohps(X;P(t)) \ . 
\end{align}
It follows that
\begin{align}
	\dstackCohpsnil(\widehat{X}_Z)=\bigsqcup_{P(t)\in \Q[t]}\dstackCohpsnil(\widehat{X}_Z; P(t)) \quad \text{and} \quad \dstackCohps(W) = \bigsqcup_{P(t)\in\Q[t]}\dstackCohps(W;P(t)) \ .
\end{align}
Furthermore, there is a canonical morphism in $\PreSt_k$
\begin{align}
	\colim_{W \in \scrT_{Z \overunder X}} \dstackCohps(W;P(t)) \longrightarrow \dstackCohpsnil(\widehat{X}_Z;P(t)) 
\end{align}
which exhibits $\dstackCohpsnil(\widehat{X}_Z;P(t))$ as the convergent completion of the colimit.

This analysis reduces us to construct admissible open covers for each $\dstackCohpsnil(\widehat{X}_Z;P(t))$. This will be achieved thanks to the properties of Harder-Narasimhan filtrations. We therefore fix a polynomial $P(t)$ of degree $m$, as well as an auxiliary monic polynomial of the same degree $\alpha(t) = \sum_{i=0}^m \alpha_i t^i/i! \in \Q[t]$. When $m \geqslant 1$, we define
\begin{align}
	\frakU_\alpha(X;P(t)) \subset \dstackCohps(X;P(t)) 
\end{align}
as the open substack parametrizing $S$-flat and properly supported families of coherent sheaves $\calE$ on $X$ for which on every geometric point $s\in S$ the inequality
\begin{align}
	p_{H\textrm{-}\min}(\calE_s, t) > \alpha(t)
\end{align}
holds. When instead $m = 0$, i.e., $P(t)$ is constant, we conventionally set
\begin{align}
	\frakU_\alpha(X;P(t)) \coloneqq \dstackCohps(X;P(t)) \ , 
\end{align}
for any (constant) $\alpha$. In either case, $\frakU_\alpha(X;P(t))$ is a quasi-compact and quasi-separated open substack of $\dstackCohps(X;P(t))$.
For a fixed thickening $W \in \scrT_{Z \overunder X}$ we set
\begin{align}
	\frakU_\alpha(W;P(t)) \coloneqq \dstackCohps(W;P(t)) \times_{\dstackCohps(X;P(t))} \frakU_\alpha(X;P(t)) \ . 
\end{align}
Similarly, we introduce
\begin{align}\label{eq:U_alpha}
	\frakU_\alpha(\widehat{X}_Z;P(t)) \coloneqq \dstackCohps(\widehat{X}_Z;P(t)) \times_{\dstackCohps(X;P(t))} \frakU_\alpha(X;P(t)) \ . 
\end{align}
We have:
\begin{lemma}
	For every choice of polynomials $P(t), \alpha(t) \in \Q[t]$ as above, $\frakU_\alpha(\widehat{X}_Z;P(t))$ is an ind-qcqs indgeometric derived stack.
\end{lemma}

\begin{proof}
	There is a canonical map in $\PreSt_k$
	\begin{align}
		\colim_{W \in \scrT_{Z \overunder X}} \frakU_\alpha(W;P(t)) \longrightarrow \frakU_\alpha(\widehat{X}_Z;P(t)) \ , 
	\end{align}
	which exhibits $\frakU_\alpha(\widehat{X}_Z;P(t))$ as the convergent completion of the colimit. Besides, for every morphism $j_{W,W'} \colon W \to W'$ in $\scrT_{Z \overunder X}$, there is an induced square
	\begin{align}
		\begin{tikzcd}[column sep=large,ampersand replacement=\&]
			\frakU_\alpha(W;P(t)) \arrow{d} \arrow{r} \& \frakU_\alpha(W';P(t)) \arrow{d} \\
			\dstackCohps(W;P(t)) \arrow{r}{\mathbf j_{W,W',\ast}} \& \dstackCohps(W';P(t)) \ ,
		\end{tikzcd} 
	\end{align}
	which furthermore is a pullback. Since $X$ is projective, Lemma~\ref{lem:perfect_resolution_property_hereditary} and Corollary~\ref{cor:smooth_implies_universal_perfect_resolution} imply that $W'$ universally satisfies the perfect resolution property. Therefore, Proposition~\ref{prop:perfect_resolution_implies_closed_immersion_of_stacks} implies that $\mathbf j_{W,W',\ast}$ is a closed immersion. Thus, the same goes for the top horizontal map in the above square. Finally, using Proposition~\ref{prop:perfect_resolution_implies_closed_immersion_of_stacks} once more, we see that the natural morphism
	\begin{align}
		\mathbf j_{W,\ast} \colon \dstackCohps(W;P(t)) \longrightarrow \dstackCohps(X;P(t)) 
	\end{align}
	is a closed immersion. In particular, it is quasi-compact and it thus follows that each $\frakU_\alpha(W;P(t))$ is quasi-compact and quasi-separated.
\end{proof}

For any $k\in \N$, $k\geq 1$, let $Z^{(k)}$ denote Gaitsgory-Rozenblyum's $k$-thickening of $Z$ along $X$ (cf.\ Remark~\ref{rem:infinitesimal_neighborhood}). We shall show the following.
\begin{proposition}\label{prop:HN-admissibility}
	Fix polynomials $P(t), \alpha(t)\in \Q[t]$ of the same degree $m$ such that $\alpha(t)$ is monic. Then, there exists $k\in \N$, with $k\geq 1$, depending on $P(t)$ and $\alpha(t)$, such that the canonical map
	\begin{align}
		\frakU_\alpha(\trunc{Z^{(k_1)}};P(t))\longrightarrow \frakU_\alpha(\trunc{Z^{(k_2)}};P(t))
	\end{align}
	is a nil-equivalence for all $k\leq k_1\leq k_2$.
	In particular, $\frakU_\alpha(\widehat{X}_Z;P(t))$ is a qcqs indgeometric derived stack.
\end{proposition}

We need some preliminary results. The following proposition will be proved in Appendix~\ref{appendix:set-theoreticity-pure}.
\begin{proposition}\label{prop:filtration-scheme-theoretically}
	Let $\calF$ be a pure coherent sheaf on $X$ of dimension $m\geq 1$, set-theoretically supported on $Z$. Then, there exists a filtration 
	\begin{align}
		0\eqqcolon\calF_{\ell+1}\subset \calF_\ell\subset \calF_{\ell-1} \subset \cdots \subset \calF_0 \coloneqq \calF \ ,
	\end{align}
	for $\ell\geq 1$, so that each subquotient is a pure $m$-dimensional sheaf with scheme-theoretic support contained in $Z$.
\end{proposition}

\begin{lemma}\label{lem:bound}
	Fix a positive integer $c$. Then, any coherent sheaf $\calF$ on $X$ of dimension $m$, set-theoretically supported on $Z$, with fixed Hilbert polynomial having $c$ as leading coefficient, which is also pure if $m\geq 1$, is scheme-theoretically supported on $\trunc{Z^{(c)}}$.
\end{lemma}

\begin{proof}
	Let $\calF$ be a coherent sheaf on $X$ of dimension $m$, set-theoretically supported on $Z$, with fixed Hilbert polynomial $P(t)$. 
	
	Let us first consider the case $m=0$. In this case, if $c=1$, i.e., if $\calF$ has length one, it is scheme-theoretically supported on $Z$. When $c\geq 2$, the sheaf $\calF$ admits a Jordan-Hölder filtration where all subquotients are zero-dimensional sheaves of length one. Thus, all subquotients are scheme-theoretically supported on $Z$. Hence, the scheme-theoretic support of $\calE$ is contained in $\trunc{Z^{(c)}}$ as a closed subscheme.\footnote{Here and it what follows, we make use of the following fact: if a coherent sheaf $\calF$ fits into a short exact sequence $0\to \calF_1\to \calF\to \calF_2\to0$, its scheme-theoretic support is contained in the union of the scheme-theoretic supports of $\calF_1$ and $\calF_2$.}
	
	Now, let us assume that $m\geq 1$. By Proposition~\ref{prop:filtration-scheme-theoretically}, there exists a filtration 
	\begin{align}
		0\eqqcolon\calF_{\ell+1}\subset \calF_\ell\subset \calF_{\ell-1} \subset \cdots \subset \calF_0 \coloneqq \calF \ ,
	\end{align}
	for $\ell\geq 1$, so that each subquotient is a pure $m$-dimensional sheaf with scheme-theoretic support contained in $Z$. In particular, the scheme-theoretic support of $\calF$ is contained in $\trunc{Z^{(\ell)}}$. Since $1\leq \ell\leq c$, we have that the scheme-theoretic support of $\calF$ is also contained in $\trunc{Z^{(c)}}$ as a closed subscheme.
\end{proof}

\begin{proposition}\label{prop:alpha-bound}
	Fix polynomials $P(t)=\sum_{i=0}^m c_i \, t^i/i!\in \Q[t]$ and $\alpha(t)=\sum_{i=0}^m \alpha_i t^i/i!\in \Q[t]$ of the same degree $m\geq 1$ such that $\alpha(t)$ is monic. Then, there exists an integer $k \geq 1$, depending on $P(t)$ and $\alpha(t)$, so that any coherent sheaf $\calE$ on $X$, set-theoretically supported on $Z$, with Hilbert polynomial $P(t)$, such that $p_{H\textrm{-}\min}(\calE) > \alpha$ is scheme-theoretically supported on $\trunc{Z^{(k)}}$. 
\end{proposition} 

\begin{proof} 
	Fix a coherent sheaf $\calE$ on $X$, set-theoretically supported on $Z$, with Hilbert polynomial $P(t)$. Consider the torsion filtration of $\calE$:
	\begin{align}
		0\subseteq \calT_0 \subseteq \cdots \subseteq \calT_{m-1}\subseteq \calE\ ,
	\end{align}
	where $\calT_i$ is the maximal $i$-dimensional subsheaf of $\calE$ for $0\leq i \leq m-1$. 
	
	If $\calE$ is pure, i.e., $\calT_{m-1}=0$, the claim follows from Lemma~\ref{lem:bound}. Otherwise, set $\calF_m\coloneqq \calE/\calT_{m-1}$. Then, $\calF_m$ is pure $m$-dimensional and set-theoretically supported on $Z$, with $c_m(\calF_m)=c_m$. By Lemma~\ref{lem:bound}, $\calF_m$ is scheme-theoretically supported on $\trunc{Z^{(c_m)}}$. 
	
	Since $p_{H\textrm{-}\min}(\calE)>\alpha$, we get $p_H(\calF_m)\geq \alpha$. In particular, we have
	\begin{align}
		P_H(\calF_m)\geq c_m\cdot \alpha \ ,
	\end{align}
	which implies
	\begin{align}
		P_H(\calT_{m-1}) = P-P_H(\calF_m)\leq P-c_m\cdot \alpha\ .
	\end{align}
	If $\calT_{m-1}$ is pure, i.e., if $\calT_{m-2}=0$, by Lemma~\ref{lem:bound}, we get that $\calT_{m-1}$ is supported on $Z^{(d_{m, m-1})}$, where $k_{m, m-1}\coloneqq \lfloor c_{m-1}-c_m\cdot \alpha_{m-1}  \rfloor$. Thus, the scheme-theoretic support of $\calE$ is contained in $\trunc{Z^{(k)}}$ as a closed subscheme, where $k \coloneqq c_m + k_{m, m-1}$. 
	
	Otherwise, consider the short exact sequence
	\begin{align}
		0\longrightarrow\calT_{m-1}/\calT_{m-2}\longrightarrow \calE/\calT_{m-2}\longrightarrow\calF_m\longrightarrow 0\ .
	\end{align}
	We have
	\begin{align}
		c_{m-1}(\calT_{m-1}/\calT_{m-2})&=c_{m-1}(\calE/\calT_{m-2})-c_{m-1}(\calF_m)\leq c_{m-1}-c_m\cdot \alpha_{m-1}\\
		&\leq  c_{m-1}-c_m\cdot \alpha_{m-1}\ .
	\end{align}
	Now, $\calT_{m-1}/\calT_{m-2}$ is a pure $(m-1)$-dimensional sheaf. Hence, by Lemma~\ref{lem:bound}, $\calT_{m-1}/\calT_{m-2}$ is scheme-theoretically supported on $\trunc{Z^{(k_{m, m-1})}}$.
	
	Furthermore, also $p_H(\calE/\calT_{m-2})\geq \alpha$. Thus,
	\begin{align}
		c_{m-2}(\calT_{m-2})=c_{m-2}-c_{m-2}(\calE/\calT_{m-2})\leq c_{m-2}-c_m\cdot \alpha_{m-2}\ .
	\end{align}
	If $\calT_{m-2}(\calE)$ is pure, i.e., if $\calT_{m-3}(\calE)=0$, by Lemma~\ref{lem:bound}, we get that $\calT_{m-2}(\calE)$ is supported on $Z^{(k_{m, m-2})}$, where $k_{m, m-2}\coloneqq \lfloor c_{m-2}-c_m\cdot \alpha_{m-2}  \rfloor$. Thus, the scheme-theoretic support of $\calE$ is contained in $\trunc{Z^{(k)}}$ as a closed subscheme, where $k \coloneqq k_m + k_{m, m-1}+k_{m, m-2}$. 	
	
	By iterating this argument, we get the claim.
\end{proof} 

Now, we are ready to prove the proposition.
\begin{proof}[Proof of Proposition~\ref{prop:HN-admissibility}]
	Let $k\geq 1$ be the integer, depending on $P(t)$ and $\alpha(t)$, whose existence is guaranteed by Proposition~\ref{prop:alpha-bound}. Consider $k\leq k_1\leq k_2$.
	
	Assume that $m\geq 1$. Let $i_{k_1, k_2}\colon \trunc{Z^{(k_1)}}\to \trunc{Z^{(k_2)}}$ be the canonical closed embedding. It is enough to prove that 
	\begin{align}
		\mathbf{i}_{k_1, k_2}\colon \frakU_\alpha(\trunc{Z^{(k_1)}};P(t))(S)\longrightarrow \frakU_\alpha(\trunc{Z^{(k_2)}};P(t))(S)
	\end{align}
	is an equivalence for any reduced scheme $S$. Let $\scrE$ be a $S$-flat family of coherent sheaves on $\trunc{Z^{(k_2)}}$ with fiber-wise Hilbert polynomial $P(t)$. Let
	\begin{align}
		\eta\colon \scrE\longrightarrow (\id_S\times i_{k_1, k_2})_\ast (\id_S\times i_{k_1, k_2})^\ast\scrE
	\end{align}
	be the canonical surjective morphism. It is enough to prove that $\eta$ is an isomorphism.
	
	Set $\scrF\coloneqq (\id_S\times i_{k_1, k_2})_\ast (\id_S\times i_{k_1, k_2})^\ast\scrE$. By Proposition~\ref{prop:alpha-bound}, the morphism $\eta_s$ is an isomorphism for any $s\in S$. Hence, the Hilbert polynomial of $\scrF_s$ coincides with the Hilbert polynomial of $\scrE_s$ for any $s\in S$. By \cite[Proposition~2.1.2]{HL_Moduli}, the Hilbert polynomial of $\scrF_s$ is locally constant as a function of $s\in S$. Since $S$ is reduced, \textit{loc.cit.} also implies that $\scrF$ is $S$-flat. Now, by \cite[Lemma~2.1.4]{HL_Moduli}, the kernel of $\eta_s$ if $\ker(\eta)_s$ for any $s\in S$. Thus, $\ker(\eta)$ has fiber-wise zero Hilbert polynomial, hence $\ker(\eta)=0$. 
	
	When $m=0$, the claim follows from a similar argument using Lemma~\ref{lem:bound}.
	
	Now, the first claim yields that $\frakU_\alpha(\widehat{X}_Z;P(t))$ is a qcqs indgeometric derived stack thanks to the description of $\widehat{X}_Z$ in Remark~\ref{rem:infinitesimal_neighborhood}.
\end{proof}

\begin{theorem}\label{thm:Cohnil_admissible}
	Let $X$ be a quasi-projective scheme over a field $k$ and let $Z$ be a closed subscheme of $X$ such that there exists a projective scheme $Y$, containing $X$ as an open subscheme, such that $Z$ is also closed in $Y$. Then, the indgeometric derived stack $\dstackCohpsnil(\widehat{X}_Z; P(t))$ is admissible for any polynomial $P(t)\in \Q[t]$.
\end{theorem}

\begin{proof}
	Since $\widehat{X}_Z\simeq \widehat{Y}_Z$, we can assume from the beginning that $X$ is projective. 
	
	Fix a sequence $\alpha_0(t)>\alpha_1(t)>\cdots$ of monic polynomials in $\Q[t]$. Then, the open substacks
	\begin{align}
		\frakU_k(\widehat{X}_Z; P(t))\coloneqq \frakU_{\alpha_k}(\widehat{X}_Z; P(t))
	\end{align}
	form an exhaustion of $\dstackCohpsnil(\widehat{X}_Z; P(t))$, which is admissible thanks to Proposition~\ref{prop:HN-admissibility} (cf.\ Remark~\ref{rem:cofinal_presentations}).
\end{proof}

We conclude this section with the following corollary.
\begin{corollary}\label{cor:Cohnil_as_formal_completion-2}
	Let $X$ be a quasi-projective scheme over a field $k$ and let $Z$ be a closed subscheme of $X$ such that there exists a projective scheme $Y$, containing $X$ as an open subscheme, such that $Z$ is also closed in $Y$. Then, the morphism \eqref{eq:nilpotent_into_ambient}
	\begin{align}\label{eq:closed-embedding}
		\Red{\dstackCohpsnil(\widehat{X}_Z)} \longrightarrow \dstackCohps(X) 
	\end{align}
	is a closed immersion.
\end{corollary}

\begin{proof}
	It is enough to prove that for every morphism $\calU \to \dstackCohps(X)$ representable by quasi-compact open immersions, setting
	\begin{align}
		\calV \coloneqq \dstackCohpsnil(\widehat{X}_Z) \times_{\dstackCohps(X)} \calU \ , 
	\end{align}
	the induced morphism $\Red{\calV} \to \calU$ is a closed immersion. Since $\dstackCohpsnil(\widehat{X}_Z)$ is admissible by Theorem~\ref{thm:Cohnil_admissible}. Thus, Corollary~\ref{cor:qc_open_in_admissible_is_qcqs} implies that $\calV$ is a qcqs indgeometric derived stack. In particular, Corollary~\ref{cor:qcqs_indgeometric_nilequivalences_in_every_presentation} guarantees that we can find a thickening $W \in \scrT_{Z \overunder X}$ such the morphism
	\begin{align}
		\calV \times_{\dstackCohpsnil(\widehat{X}_Z)} \dstackCohps(W) \longrightarrow \calV 
	\end{align}
	is a nil-equivalence. In particular, the map $\Red{\calV} \to \calU$ is canonically identified with the map
	\begin{align}
		\Red{\dstackCohps(W)} \times_{\dstackCohps(X)} \calU \longrightarrow \calU \ . 
	\end{align}
	It is therefore enough to argue that the map $\dstackCohps(W) \to \dstackCohps(X)$ is a closed immersion. Since $X$ satisfies the perfect resolution property (cf.\ Example~\ref{eg:quasi_projective_implies_perfect_resolution}), this follows from Proposition~\ref{prop:perfect_resolution_implies_closed_immersion_of_stacks}.
\end{proof}

\bigskip\section{Nilpotent cohomological Hall algebra}\label{sec:COHA}

In this section we introduce the cohomological Hall algebra associated to the category of coherent sheaves on a smooth quasi-projective surface over a field $k$ of characteristic zero, set-theoretically supported on a fixed closed subscheme. We shall call it the \textit{nilpotent} COHA. 

\subsection{The $\Lambda$-graded $2$-Segal structure on the Waldhausen construction}\label{subsec:2_Segal}

 Given an integer $n$ we let $[n]$ denote the linearly ordered poset $\{0 < 1 < \cdots < n\}$ and we set
\begin{align}
	\sfT_n \coloneqq \Fun([1],[n]) \ . 
\end{align}
The collection of the various $\sfT_n$ determines a functor
\begin{align}
	\sfT_\bullet \colon \bfDelta \longrightarrow \Cat_\infty \ . 
\end{align}
Given a stable $\infty$-category $\calC$, we set
\begin{align}
	\calS_n \subseteq \Fun(\sfT_n, \calC) 
\end{align}
be the full subcategory spanned by those functors $F \colon \sfT_n \to \calC$ satisfying the following two conditions:
\begin{enumerate}\itemsep=0.2cm
	\item $F(i,i) \simeq 0$;
	\item for every $0 \leqslant i < j \leqslant n-1$, the square
	\begin{align}
		\begin{tikzcd}[ampersand replacement=\&]
			F(i,j) \arrow{r} \arrow{d} \& F(i+1, j) \arrow{d} \\
			F(i,j+1) \arrow{r} \& F(i+1, j+1)
		\end{tikzcd}
	\end{align}
	is a pullback in $\calC$.
\end{enumerate}
The $\infty$-categories $\calS_n \calC$ depends simplicially on $n$, i.e., they assemble into a simplicial object
\begin{align}
	\calS_\bullet \calC \colon \bfDelta\op \longrightarrow \Cat_\infty \ , 
\end{align}
known as the \textit{Waldhausen construction} of $\calC$. One of the main observations of \cite{Dyckerhoff_Kapranov_Higher_Segal} is that this satisfies the $2$-Segal condition. The functoriality of this construction allow to replace $\calC$ by a functor
\begin{align}
	\scrC \colon \dAff_k\op \longrightarrow \Cat_\infty^{\mathsf{st}} \ , 
\end{align}
giving rise to a $2$-Segal object in derived prestacks
\begin{align}
	\calS_\bullet \scrC \colon \bfDelta\op \longrightarrow \PreSt_k \ , 
\end{align}
which we refer to as the Waldhausen construction of $\scrC$, Besides, this construction is obviously functorial in $\scrC$.

\begin{notation}
	Given a stable $\infty$-category $\calC$, we refer to $F \in \calS_n \calC$ as an \textit{$n$-flag of objects in $\calC$}.
	Given $0 \leqslant i \leqslant j \leqslant n$, we set
	\begin{align}\label{eq:ev}
		\ev_{i,j}(F) \coloneqq F(i,j) \ , 
	\end{align}
	We also denote by
	\begin{align}
		\partial_i(F) \in \calS_{n-1} \calC 
	\end{align}
	the $(n-1)$-flag obtained restricting $F$ along the morphism $[n-1] \to [n]$ in $\bfDelta$ that misses $i$.
	In particular, when $n = 2$ we have the overlap of notation
	\begin{align}\label{eq:partial}
		\partial_0 = \ev_{1,2} \ , \qquad \partial_1 = \ev_{0,2} \ , \qquad \partial_2 = \ev_{0,1} \ . 
	\end{align}
		We refer to the collection of objects $\{\ev_{i,i+1}(F)\}_{0 \leqslant i \leqslant n-1}$ as the \textit{diagonal of $F$}. Sending an $n$-flag to its diagonal gives rise to natural transformation
	\begin{align}
		\delta_{\calC, n} \colon \calS_n \calC \longrightarrow \prod_{i = 0}^{n-1} \calC \ . 
	\end{align}
	When $\calC$ is replaced by a functor $\scrC \colon \dAff_k\op \to \Cat_\infty$, this induces a natural transformation
	\begin{align}
		\delta_{\scrC, n} \colon \calS_n \scrC \longrightarrow \prod_{i = 0}^{n-1} \scrC 
	\end{align}
	of derived $\Cat_\infty$-valued prestacks.
\end{notation}

We apply this to the setting of \S\ref{sec:nilpotent_coherent_sheaves}. Specifically we fix a quasi-compact and quasi-separated derived formal scheme $\calX$ and we consider the functor
\begin{align}
	\dAff_k\op \longrightarrow \Cat_\infty 
\end{align}
sending $S$ to $\catCohnil_{S,\ps}(\calX \times S)$.
Applying the Waldhausen construction and passing to the maximal $\infty$-groupoid at the end, we obtain a $2$-Segal object
\begin{align}
	\calS_\bullet \widetilde{\dstackCoh}_{\ps}^{\mathsf{nil}}(\calX) \colon \bfDelta \longrightarrow \PreSt_\C \ .
\end{align}
We set
\begin{align}
	\calS_\bullet \dstackCohpsnil(\calX) \coloneqq \conv{\big( \calS_\bullet \widetilde{\dstackCoh}_{\ps}^{\mathsf{nil}}(\calX) \big)} \ . 
\end{align}
When $n = 1$, there is a canonical identification
\begin{align}
	\calS_1 \dstackCohpsnil(\calX) \simeq \dstackCohpsnil(\calX) \ . 
\end{align}
When $\calX \simeq \widehat{X}_Z$ is a formal completion, we saw that under mild assumptions on $X$ this derived stack is indgeometric or even admissible. We now generalize these results to arbitrary $n$. The key ingredient is the following:
\begin{theorem}\label{thm:embedded_2_Segal}
	Let $X$ be a quasi-compact and quasi-separated derived scheme over a field $k$ of characteristic zero and let $j \colon Z \hookrightarrow X$ be a closed immersion. Then for every $n \geqslant 1$ the square
	\begin{align}\label{eq:embedded_2_Segal}
		\begin{tikzcd}[ampersand replacement=\&]
			\calS_n \dstackCohpsnil(\widehat{X}_Z) \arrow{r} \arrow{d} \& \calS_n \dstackCoh_\ps(X) \arrow{d} \\
			\prod_{i = 1}^n \dstackCohpsnil(\widehat{X}_Z) \arrow{r} \& \prod_{i = 1}^n \dstackCoh_\ps(X)
		\end{tikzcd}
	\end{align}
	whose vertical morphisms send an $n$-flag to its diagonal, is a pullback.
\end{theorem}

\begin{proof}
	It follows from Theorem~\ref{thm:nilpotent_vs_set_theoretic} that both horizontal arrows are fully faithful on $S$ points for every $S \in \dAff_k$. Thus, the same goes for the induced morphism from $\calS_n \dstackCohpsnil(\widehat{X}_Z)$ to the fiber product. We are therefore left to argue for its essential surjectivity. Notice that since both source and target are convergent, we can restrict without loss of generality to the case where $S \in \evccndAff$. Unraveling the definition of the Waldhausen construction, we reduce ourselves to check that if
	\begin{align}
		\calF' \longrightarrow \calF \longrightarrow \calF'' 
	\end{align}
	is a fiber sequence in $\catCoh_{S,\ps}(X \times S)$ and both $\calF'$ and $\calF''$ belong to the essential image of $\catCohnil_{S,\ps}(\widehat{X}_Z)$, then the same goes for $\calF$, and this is immediate from Theorem~\ref{thm:nilpotent_vs_set_theoretic}, Corollary~\ref{cor:computing_Cohnil} and Lemma~\ref{lem:pushforward_coherent}.
\end{proof}

\begin{corollary}\label{cor:embedded_2_Segal}
	Let $X$ be a quasi-compact and quasi-separated derived scheme over a field $k$ of characteristic zero and let $j \colon Z \hookrightarrow X$ be a closed immersion. Then
	\begin{enumerate}\itemsep=0.2cm
		\item \label{item:embedded_2_Segal-1} if $X$ universally satisfies the perfect resolution property, then for every integer $n \geqslant 0$ the derived stack $\calS_n \dstackCohpsnil(\widehat{X}_Z)$ is indgeometric;
		
		\item\label{item:embedded_2_Segal-2} if $X$ is quasi-projective scheme admitting a projective compactification $Y$ that contains $Z$ as a closed subscheme, then for every integer $n \geqslant 0$ the indgeometric derived stack $\calS_n \dstackCohpsnil(\widehat{X}_Z)$ is admissible.
	\end{enumerate}
\end{corollary}

\begin{proof}	
	Assumption~\eqref{item:embedded_2_Segal-1} guarantees via Theorem~\ref{thm:indgeometricity} that the product $\prod_{i = 0}^{n-1} \dstackCohpsnil(\widehat{X}_Z)$ is indgeometric. Then Theorem~\ref{thm:embedded_2_Segal} guarantees that the pullback of a presentation for the above product via the diagonal $\delta_n$ provides a presentation for $\calS_n \dstackCohpsnil(\widehat{X}_Z)$.
	
	We now prove \eqref{item:embedded_2_Segal-2}. First, recall from \cite[Corollary~4.6]{Porta_Sala_Hall} that $\calS_n \dstackCohps(X)$ is geometric. Then \eqref{item:embedded_2_Segal-2} follows combining Theorems~\ref{thm:Cohnil_admissible} and \ref{thm:embedded_2_Segal} with Corollary~\ref{cor:admissible_stacks_finite_limits}.
\end{proof}

Fix now a free abelian group of finite rank $(\Lambda,+)$ and a map
\begin{align}
	v \colon \pi_0(\dstackCohnil_{\ps}(\calX)) \longrightarrow \Lambda \ , 
\end{align}
which we assume to be additive in the sense of \cite[\S7.3]{Alper-Halpern-Leistner-Heinloth}. We use the notion of $\Lambda$-graded $2$-Segal object of \cite[\S\ref*{torsion-pairs-sec:Lambda_graded}]{DPS_Torsion-pairs} to study the compatibility between the group structure of $\Lambda$ and the Hall multiplication. Given $\bfv \in \Lambda$, we let
\begin{align}
	\dstackCohnil_{\ps}(\calX;\bfv) 
\end{align}
be the corresponding union of connected components.
More generally, given an element $\underline{\bfv} = (\bfv_{0,1}, \bfv_{1,2}, \ldots, \bfv_{n-1,n}) \in \Lambda^n$, we let
\begin{align}
	\calS_n \dstackCohnil_{\ps}(\calX;\underline{\bfv}) 
\end{align}
be the open and closed substack of $\calS_n \dstackCohnil_{\ps}(\calX)$ parametrizing $n$-flags whose diagonal is of type $\underline{\bfv}$.
A routine verification similar to \cite[Lemma~4.1]{Porta_Sala_Hall} shows that this satisfies the $\Lambda$-graded $2$-Segal condition.
Since the convergent completion functor $\conv{(-)}$ commutes with limits (see Recollection~\ref{recollection:convergence}), it automatically follows that $\calS_\bullet \dstackCohpsnil(X)$ is again a $2$-Segal object in $\PreSt_\C$. In particular, \cite[Proposition~\ref*{torsion-pairs-prop:Lambda_graded_2_Segal}]{DPS_Torsion-pairs} shows that $\calS_\bullet \dstackCohpsnil(\calX)$ determines an $\E_1$-monoid in $\Corr^{\times_\Lambda}(\LambdadSt_\C)$.
In what follows, we tacitly consider $\calS_\bullet \dstackCohnil_{\ps}(\calX)$ as a $\Lambda$-graded $2$-Segal derived stack.

\subsection{Functoriality}\label{subsec:COHA-functoriality}

We now discuss the functorial dependency of $\calS_\bullet \dstackCohpsnil(\calX)$ on the derived formal scheme $\calX$. Ultimately, this relies on the functoriality of relative ind-coherent sheaves discussed in Appendix~\ref{appendix:relative_ind_coherent}.

\begin{construction}
	Given a morphism $f \colon \calX \to \calY$ of derived formal schemes over a derived affine $S$ we have by \ref{construction:relative_IndCoh_Kan_extension}  a well defined pushforward functor
	\begin{align}
		f_\ast \colon \catIndCoh(\calX/S) \longrightarrow \catIndCoh(\calY/S) \ . 
	\end{align}
	Assume that $f$ is ind-finite, that is that it admits a presentation as $\colim_\alpha f_\alpha$ where $f_\alpha \colon X_\alpha \to Y_\alpha$ is a finite morphism between derived quasi-compact and quasi-separated schemes (and $\{X_\alpha\}$ and $\{Y_\alpha\}$ are presentations for the formal schemes $\calX$ and $\calY$). In this case, Lemma~\ref{lem:pushforward_coherent} implies that $f_{S,\ast}$ descends to a morphism
	\begin{align}
		f_{S,\ast} \colon \catCohnil_{S,\ps}(\calX \times S) \longrightarrow \catCohnil_{S,\ps}(\calY \times S) 
	\end{align}
	for every $S \in \dAff$.
	In turn, this provides a well defined morphism
	\begin{align}
		\tensor*[^{(\bullet)}]{\bff}{_\ast} \colon \calS_\bullet \dstackCohpsnil(\calX) \longrightarrow \calS_\bullet \dstackCohpsnil(\calY) \ . 
	\end{align}
	When $\bullet = 1$, we simply write $\bff_\ast$ instead of $\tensor*[^{(1)}]{\bff}{_\ast}$.
\end{construction}

We will be particularly interested in the following example.
\begin{example}\label{eg:restriction_to_components}
	Let $X$ be a quasi-compact and quasi-separated derived scheme over a field $k$ of characteristic zero and let
	\begin{align}
		\begin{tikzcd}[column sep=small, ampersand replacement=\&]
			Z_1 \arrow{rr}{j_{12}} \arrow{dr}[swap]{j_1} \& \& Z_2 \arrow{dl}{j_2} \\
			\& X
		\end{tikzcd}
	\end{align}
	be closed immersions. This induces an ind-finite morphism
	\begin{align}
		\hat{\jmath}_{12} \colon \widehat{X}_{Z_1} \longrightarrow \widehat{X}_{Z_2} \ , 
	\end{align}
	and therefore a morphism
	\begin{align}
		\tensor*[^{(\bullet)}]{\bfjmathhat}{_{12,\ast}} \colon \calS_\bullet \dstackCohpsnil(\widehat{X}_{Z_1}) \longrightarrow \calS_\bullet \dstackCohpsnil(\widehat{X}_{Z_2}) \ . \tag*{\qedhere} 
	\end{align}
\end{example}

In the setting of the above example, the morphism of $2$-Segal derived stacks is strict in the following sense:
\begin{lemma}\label{lem:strict_2_Segal_morphism}
	In the setting of Example~\ref{eg:restriction_to_components} the square
	\begin{align}
		\begin{tikzcd}[column sep=large, ampersand replacement=\&]
			\calS_2 \dstackCohpsnil(\widehat{X}_{Z_1}) \arrow{r}{\tensor*[^{(2)}]{\bfjmathhat}{_{12,\ast}}} \arrow{d}{\delta} \& \calS_2 \dstackCohpsnil(\widehat{X}_{Z_2}) \arrow{d}{\delta} \\
			\dstackCohpsnil(\widehat{X}_{Z_1}) \times \dstackCohpsnil(\widehat{X}_{Z_1}) \arrow{r}{\bfjmathhat_{12,\ast} \times \bfjmathhat_{12,\ast}} \& \dstackCohpsnil(\widehat{X}_{Z_2}) \times \dstackCohpsnil(\widehat{X}_{Z_2})
		\end{tikzcd} 
	\end{align}
	is a pullback.
\end{lemma}

\begin{proof}
	This is an immediate consequence of Theorem~\ref{thm:embedded_2_Segal} and the transitivity property for pullback squares.
\end{proof}

\subsection{The main existence theorem}\label{subsec:main-existence-theorem}

We are now ready to establish the existence of the nilpotent cohomological Hall algebra. 

We start with some general considerations. Let $X$ be a $n$-dimensional smooth scheme over a field $k$ of characteristic zero and let $j \colon Z \hookrightarrow X$ be the inclusion of a closed subscheme. Assume that $X$ admits a projective compactification $\overline{X}$ that contains $Z$ as a closed subscheme. We let $N(X)$ be the numerical Grothendieck group of the category of properly supported coherent sheaves on $X$, and we set
\begin{align}
	N_{\leqslant 1}(X)\coloneqq \mathsf{Im}(K_0(\catCoh_{\leqslant 1}(X)) \longrightarrow N(X) ) \ .
\end{align}
We take $\Lambda \coloneqq N_{\leqslant 1}(X)$.
We have the following.
\begin{lemma}\label{lem:lci}
	\hfill
	\begin{enumerate}\itemsep=0.2cm
		\item \label{item:lci-1} the map
		\begin{align}
			\partial_0 \times \partial_2 \colon \calS_2 \dstackCohpsnil(\widehat{X}_Z) \longrightarrow \dstackCohpsnil(\widehat{X}_Z) \times \dstackCohpsnil(\widehat{X}_Z) 
		\end{align}
		is representable by linear stacks and its tor-amplitude is contained in $[1-n,1]$. In particular, if $X$ is a surface then $\partial_0 \times \partial_2$ is derived lci.
		
		\item \label{item:lci-2} Asume that $X$ is a surface. Then the map
		\begin{align}
			\partial_1 \colon \calS_2 \dstackCohpsnil(\widehat{X}_Z) \longrightarrow \dstackCohpsnil(\widehat{X}_Z) 
		\end{align}
		of $\Lambda$-graded derived stacks is locally rpas (in the sense of Definitions~\ref{def:rpas-connected} and \ref{def:admissible-rpas-connected}).
	\end{enumerate}
\end{lemma}

\begin{proof}
	Point \eqref{item:lci-1} is an immediate consequence of Theorem~\ref{thm:embedded_2_Segal} and \cite[Proposition~3.11]{Porta_Sala_Hall}.
	As for point \eqref{item:lci-2}, we already know that the analogous statement holds for $\dstackCoh_{\ps}(X)$, since in the latter case it is representable by Quot schemes, which are known to be proper as soon as the Hilbert polynomial is fixed. Therefore, it is enough to argue that the square
	\begin{align}
		\begin{tikzcd}[column sep=large, ampersand replacement=\&]
			\calS_2 \dstackCohpsnil(\widehat{X}_Z) \arrow{r}{\tensor*[^{(2)}]{\bfjmathhat}{_{12,\ast}}} \arrow{d}{\partial_1} \& \calS_2 \dstackCoh_{\ps}(X) \arrow{d}{\partial_1} \\
			\dstackCohpsnil(\widehat{X}_Z) \arrow{r}{\bfjmathhat_{12,\ast}} \& \dstackCoh_{\ps}(X)
		\end{tikzcd} 
	\end{align}
	is a pullback square. This can be done after evaluating on a test scheme $S \in \dAff$, and since all the stacks that appear are convergent, we can assume without loss of generality that $S \in \evccndAff$.
	It follows from Theorem~\ref{thm:nilpotent_vs_set_theoretic} that both the top and the bottom horizontal maps are fully faithful. So, unraveling the definitions we are reduced to show that if
	\begin{align}
		\calF' \longrightarrow \calF \longrightarrow \calF'' 
	\end{align}
	is a fiber sequence in $\catCoh_{S,\ps}(X \times S)$ and $\calF$ lies in the essential image of $\bfjmathhat_{12,\ast}$, the same goes for $\calF'$ and $\calF''$. Since $S$ is eventually coconnective, Theorem~\ref{thm:nilpotent_into_set_theoretic_moduli} reduces us to show that
	\begin{align}
		\calF' \vert_{X \times S \smallsetminus Z \times S} \simeq \calF'' \vert_{X \times S \smallsetminus Z \times S} \simeq 0 \ . 
	\end{align}
	Since both $\calF'$ and $\calF''$ are $S$-flat, this is automatic.
\end{proof}

Assume now that $X$ admits a projective compactification $\overline{X}$ that contains $Z$ as a closed subscheme (notice that this is automatic if $X$ is quasi-projective and $Z$ is itself proper). Then Corollary~\ref{thm:embedded_2_Segal} implies that each $\calS_n \dstackCohpsnil(\widehat{X}_Z)$ is admissible. Fixing a motivic formalism $\bfD^\ast$ as in \cite[\S\ref*{torsion-pairs-subsec:motivic_formalism}]{DPS_Torsion-pairs}, $\calA \in \CAlg(\bfD(\Spec(k)))$ and an abelian subgroup $\Gamma \subset \Pic(\bfD(S))$, we can consider
\begin{align}
	\HBMDGamma_0(\dstackCohpsnil(\widehat{X}_Z); \calA) \in \Pro(\mathsf{Ab}^{\gr}) \ . 
\end{align}
This is a pro-$\Gamma$-graded abelian group. In practice, we will consider it as a topological abelian group with the topology induced from the pro-structure. We refer to this topology as \textit{the quasi-compact topology} (see Remark~\ref{rem:quasi-compact_topology}).

\begin{theorem}\label{thm:nilpotent-COHA}
	Let $X$ be a smooth surface over a field $k$ of characteristic zero and let $j \colon Z \hookrightarrow X$ be the inclusion of a closed subscheme. Assume that $X$ admits a projective compactification $\overline{X}$ that contains $Z$ as a closed subscheme. Then for every choice of an $R$-linear motivic formalism $\bfD^\ast$ and every oriented $\calA \in \CAlg(\bfD(\Spec(k)))$ and every abelian subgroup $\Gamma \subset \Pic(\bfD(\Spec(k)))$ closed under Thom twists, the $\Lambda$-graded pro-$R$-module $\HBMDGamma_0(\dstackCohpsnil(\widehat{X}_Z); \calA)$ is an associative algebra, which we shall denote by $\DGammacoha_{\widehat{X}_Z,\calA}$, with multiplication
	\begin{align}
		(\partial_1)_\ast \circ (\partial_0\times \partial_2)^!\circ \boxtimes \ ,
	\end{align}
	where the maps $\partial_i$, for $i=0, 1, 2$, fit the correspondence
	\begin{align}
		\begin{tikzcd}[column sep=small, ampersand replacement=\&]
			\& \calS_2 \dstackCohpsnil(\widehat{X}_Z) \arrow{dr}{\partial_1} \arrow{dl}[swap]{\partial_0 \times \partial_2} \\
			\dstackCohpsnil(\widehat{X}_Z) \times \dstackCohpsnil(\widehat{X}_Z) \& \& \dstackCohpsnil(\widehat{X}_Z)
		\end{tikzcd}
	\end{align}
	and $\boxtimes$ is the exterior product. In particular, seeing $\DGammacoha_{\widehat{X}_Z,\calA}$ as a topological abelian group, the Hall multiplication defined above becomes continuous with respect to the quasi-compact topology. Moreover:
	\begin{enumerate}\itemsep0.2cm
		\item \label{item:nilpotent-COHA-1} let $j' \colon Z' \hookrightarrow X'$ be another closed embedding into a smooth surface over $k$ satisfying the above assumption.
		Then, an isomorphism of formal completions $\widehat{X'}_{Z'} \simeq \widehat{X}_Z$ induces a functorial and continuous isomorphism of algebras
		\begin{align}
			\DGammacoha_{\widehat{X'}_{Z'},\calA} \simeq \DGammacoha_{\widehat{X}_Z,\calA} \ . 
		\end{align}

		\item \label{item:nilpotent-COHA-2} Let $j' \colon Z' \hookrightarrow Z$ be a nested closed subscheme inside $X$. Then, the morphism
		\begin{align}
			\bfjmathhat'_\ast \colon \dstackCohpsnil(\widehat{X}_{Z'}) \longrightarrow \dstackCohpsnil(\widehat{X}_Z) 
		\end{align}
		induces a continuous morphism of algebras
		\begin{align}
			j'_\ast \colon \DGammacoha_{\widehat{X}_{Z'},\calA} \longrightarrow \DGammacoha_{\widehat{X}_Z,\calA} \ . 
		\end{align}
		
		\item \label{item:nilpotent-COHA-3} Let $(\bfD', \calA', \Gamma')$ be a second motivic formalism equipped with an oriented ring of coefficients $\calA'$ and an abelian subgroup $\Gamma'$ closed under Thom twists.
		A transformation
		\begin{align}
			(s,\phi) \colon (\bfD, \calA, \Gamma) \longrightarrow (\bfD', \calA', \Gamma') 
		\end{align}
		as in \cite[Theorem~\ref*{torsion-pairs-thm:functoriality_of_BM_homology}]{DPS_Torsion-pairs} induces a morphism of algebras
		\begin{align}
			\DGammacoha_{\widehat{X}_Z, \calA} \longrightarrow \coha^{\bfD',\Gamma'}_{\widehat{X}_Z, \calA'} \ .
		\end{align}
	\end{enumerate}
\end{theorem}

\begin{proof}
	We start by establishing the existence of the Hall multiplication. To begin with, using \cite[Proposition~\ref*{torsion-pairs-prop:Lambda_graded_2_Segal}]{DPS_Torsion-pairs} we can convert the $\Lambda$-graded $2$-Segal object $\calS_\bullet \dstackCohpsnil(\widehat{X}_Z)$ constructed in \S\ref{subsec:2_Segal} into an $\E_1$-monoid object in $\Corr^{\times_\Lambda}(\LambdaindGeomadm_k)$. Besides, Lemma~\ref{lem:lci} and the $2$-Segal property allow to promote this into an $\E_1$-monoid object in $\Corr^{\times_\Lambda}(\LambdaindGeomadm_k)_{\mathsf{qc.lci} \: \cap \: \mathsf{fconn},\mathsf{lrpas}}$. Since $\HBMDGamma_0(-;\calA)$ is a lax-monoidal functor in virtue of Theorem~\ref{thm:BM_Lambda_graded_functoriality}, we deduce that $\DGammacoha_{\widehat{X}_Z, \calA}$ acquires the structure of an associative algebra, whose product is given by the prescribed rule.
	
	\medskip
	
	Concerning the extra statements, \eqref{item:nilpotent-COHA-1} is simply a consequence of the fact that $\calS_\bullet \dstackCohpsnil(\widehat{X}_Z)$, as a derived stack, only depends on the formal completion $\widehat{X}_Z$. The ambient scheme $X$ was only used to guarantee, via Theorem~\ref{thm:embedded_2_Segal} and Lemma~\ref{lem:lci} that $\partial_0 \times \partial_2$ is representable by finitely connected (in the sense of Definitions~\ref{def:rpas-connected} and \ref{def:admissible-rpas-connected}), quasi-compact, lci derived geometric stacks, and that $\partial_1$ is locally rpas (in the sense of Definitions~\ref{def:rpas-connected} and \ref{def:admissible-rpas-connected}).
	
	\medskip
	
	We now prove \eqref{item:nilpotent-COHA-2}. To begin with, observe that Corollary~\ref{cor:Cohnil_as_formal_completion-2} implies that both
	\begin{align}
		\Red{\dstackCohpsnil(\widehat{X}_{Z'})} \qquad \text{and} \qquad \Red{\dstackCohpsnil(\widehat{X}_Z)} 
	\end{align}
	are closed inside $\dstackCoh_{\ps}(X)$. In particular, $\bfjmathhat_\ast$ is a nil-closed immersion, and therefore it induces a well defined continuous pushforward map
	\begin{align}
		j'_\ast \colon \DGammacoha_{\widehat{X}_{Z'},\calA} \longrightarrow \DGammacoha_{\widehat{X}_Z, \calA} \ . 
	\end{align}
	The compatibility with the Hall multiplication is a direct consequence of the functorial behavior of $\HBMDGamma_0(-;\calA)$ and of Lemma~\ref{lem:strict_2_Segal_morphism}. Finally, statement \eqref{item:nilpotent-COHA-3} follows directly from the last part of Theorem~\ref{thm:functoriality_of_BM_homology_admissible}.
\end{proof}

\begin{remark}\label{rem:ignoring_Lambda_grading}
	The Hall multiplication is compatible with the addition of $\Lambda$. In fact, the above proof exhibits $\DGammacoha_{\widehat{X}_Z,\calA}$ as an algebra object in $\Fun(\Lambda,\Modd_R)$, where the latter is equipped with Day's convolution. This is the ultimate reason we used the language of $\Lambda$-graded $2$-Segal objects.
\end{remark}

We are mostly interested in the following three examples.\begin{example}\label{ex:three-cohas}
	\hfill
	\begin{enumerate}\itemsep=0.2cm
	\item when $\bfD^\ast \coloneqq \mathbf{DM}_\Q^\ast$ is the rational Voevodsky's formalism of \cite[Example~\ref*{torsion-pairs-ex:genuine_motivic_formalism}]{DPS_Torsion-pairs} and $\calA \coloneqq \sfH \Q$ is the motivic Eilenberg-MacLane $\E_\infty$-ring spectrum, and $\Gamma \coloneqq \Z\langle 1 \rangle$, we write
	\begin{align}
		\motcoha_{X, Z} \coloneqq 	\DGammacoha_{\widehat{X}_Z,\sfH \Q} \ ,
	\end{align}
	whose underlying topological abelian group is
	\begin{align}
		\Hmotbullet(\dstackCohpsnil(\widehat{X}_Z); 0)\coloneqq \HBMDGamma_0(\dstackCohpsnil(\widehat{X}_Z); \sfH \Q) \ .
	\end{align}
	This recovers (motivically defined) Chow groups of $\dstackCohpsnil(\widehat{X}_Z)$. Thus, we refer to $\motcoha_{X, Z}$ as the \textit{motivic cohomological Hall algebra of $\dstackCohpsnil(\widehat{X}_Z)$}. 
	
	By replacing $\dstackCohpsnil(\widehat{X}_Z)$ with $\dstackCohps(X)$ and using the results of \cite[\S4]{Porta_Sala_Hall} with Theorem~\ref{thm:functoriality_of_BM_homology_admissible}, we also obtain the \textit{motivic cohomological Hall algebra of $\dstackCohps(X)$}.
	
	\item When $\bfD^\ast \coloneqq \mathbf{DM}_\Q^\ast$ is the rational Voevodsky's formalism of \cite[Example~\ref*{torsion-pairs-ex:genuine_motivic_formalism}]{DPS_Torsion-pairs} and $\calA \coloneqq \mathsf{KGL}^{\mathsf{et}}$ is the étale hypersheafification of the algebraic $K$-theory spectrum, and $\Gamma \coloneqq \Z\langle 1 \rangle$, we write
	\begin{align}
		\bfG_{X, Z} \coloneqq \DGammacoha_{\widehat{X}_Z,\mathsf{KGL}^{\mathsf{et}}}\ ,
	\end{align}
	whose underlying topological abelian group is
	\begin{align}
		\bfG(\dstackCohpsnil(\widehat{X}_Z))\coloneqq \HBMDGamma_0(\dstackCohpsnil(\widehat{X}_Z); \mathsf{KGL}^{\mathsf{et}}) \ .
	\end{align}
	This recovers the algebraic $G$-theory of $\dstackCohpsnil(\widehat{X}_Z)$. Taking its $\pi_0$, we obtain the \textit{K-theoretical Hall algebra of $\dstackCohpsnil(\widehat{X}_Z)$}. By replacing $\dstackCohpsnil(\widehat{X}_Z)$ with $\dstackCohps(X)$, we recover the (full) K-theoretical Hall algebra of $\dstackCohps(X)$, as defined in \cite[\S5]{Porta_Sala_Hall} (see also \cite{COHA_surface} for an equivalent realization of the K-theoretical Hall algebra of $\dstackCohps(X)$).
	
	\item When $\bfD^\ast \coloneqq \BettiD$ is the topological formalism of \cite[\S\ref*{torsion-pairs-ex:topological_formalism}]{DPS_Torsion-pairs}, $\calA \coloneqq \Q$, and $\Gamma \coloneqq \Z\langle 1/2 \rangle$ (see \cite[Remark~\ref*{torsion-pairs-rem:explicit_bigrading}]{DPS_Torsion-pairs}), we simply write
	\begin{align}
		\coha_{X, Z} \coloneqq  \DGammacoha_{\widehat{X}_Z, \Q}\ ,
	\end{align}
	whose a graded pro-$\Q$-vector space is 
	\begin{align}
		\HBMbullet(\dstackCohpsnil(\widehat{X}_Z))\coloneqq \HBMGamma_0(\dstackCohpsnil(\widehat{X}_Z); \Q) \ .
	\end{align}
	We shall call $\coha_{X, Z}$ the \textit{cohomological Hall algebra of $\dstackCohpsnil(\widehat{X}_Z)$}. 
	
	\medskip
	
	We will also be interested in taking $\Gamma \coloneqq \Z \langle 1 \rangle$, in which case we write
	\begin{align}
		\evencoha_{X, Z} \coloneqq  \DGammacoha_{\widehat{X}_Z, \Q}\quad\text{and}\quad \HBM_{\mathsf{even}}(\dstackCohpsnil(\widehat{X}_Z))\coloneqq \HBMDGamma_0(\dstackCohpsnil(\widehat{X}_Z); \Q)\ .
	\end{align}
	The natural inclusion $\Z\langle 1 \rangle \subset \Z\langle 1/2 \rangle$ induces a continuous morphism of algebras
	\begin{align}
		\evencoha_{X, Z} \longrightarrow \coha_{X, Z} \quad\text{and}\quad \HBM_{\mathsf{even}}(\dstackCohpsnil(\widehat{X}_Z))\longrightarrow \HBMbullet(\dstackCohpsnil(\widehat{X}_Z))\ .
	\end{align}
	By replacing $\dstackCohpsnil(\widehat{X}_Z)$ with $\dstackCohps(X)$, we recover the (even) cohomological Hall algebra of $\dstackCohps(X)$, which has been defined in \cite{COHA_surface}.
	
	Also, the natural transformation $\bfDM_\Q \to \BettiD_\Q$ of \cite[Remark~\ref*{torsion-pairs-rem:comparison_motivic_to_topological}]{DPS_Torsion-pairs} induces via Theorem~\ref{thm:nilpotent-COHA}--\eqref{item:nilpotent-COHA-3} a continuous morphism of algebras
	\begin{align}
		\motcoha_{X,Z} \longrightarrow \evencoha_{X,Z} \ , 
	\end{align}
	given by the cycle class map. \qedhere
\end{enumerate}
\end{example}

\begin{remark}\label{rem:equivariant-BM}
	Assume that there is a torus $T$ acting on $X$ such that $Z$ is $T$-invariant. Following the rough lines of \cite[\S4.3]{Porta_Sala_Hall} we can adapt the above results in the $T$-equivariant setting as follows. As observed in Variant~\ref{variant:nilpotent_qcoh_equivariant_case}, the $T$-action descends to a $T$-action on $\widehat{X}_Z$, and the structural morphism $[\widehat{X}_Z / T] \to \sfB T$ is representable by derived formal schemes. Thus, $\catQCohnil([\widehat{X}_Z/T])$ is well defined, which in turn allows to define a $2$-Segal object
	\begin{align}
		\calS_\bullet\catCohnil_{\ps,\sfB T}( [\widehat{X}_Z / T] ) \in \Fun_{2\textrm{-}\mathsf{Seg}}(\bfDelta\op, \dSt_{/\sfB T}) \ . 
	\end{align}
	By construction
	\begin{align}
		\Spec(k) \times_{\sfB T} \calS_\bullet\catCohnil_{\ps,\sfB T}( [\widehat{X}_Z / T] ) \simeq \calS_\bullet \catCohnil_{\ps}(\widehat{X}_Z) \ . 
	\end{align}
	Therefore, Corollary~\ref{cor:testing_admissibility_on_an_atlas} implies that each $\calS_\bullet \catCohnil_{\ps,\sfB T}( [\widehat{X}_Z / T] )$ is admissible. At this point, thus applying Theorem~\ref{thm:functoriality_of_BM_homology_admissible} we obtain a variant of Theorem~\ref{thm:nilpotent-COHA} in the equivariant setting. In particular, there exists an associative algebra structure on $\HBMDGamma_0(\dstackCohpsnil(\widehat{X}_Z)/T; \calA)$, which we shall denote by $\DGammaTcoha_{\widehat{X}_Z,\calA}$, and equivariant versions of statements~\eqref{item:nilpotent-COHA-1}, \eqref{item:nilpotent-COHA-2}, and \eqref{item:nilpotent-COHA-3} hold. In particular, we denote by $\motTcoha_{X, Z}$, $\bfG_{X, Z}^T$, and $\coha^T_{X, Z}$ the equivariant versions of the examples discussed in Example~\ref{ex:three-cohas}.
\end{remark}

\begin{remark}
	Let $C$ be a smooth projective curve over a field $k$ and let $X\coloneqq T^\ast C$ be its cotangent bundle. When $Z$ is the zero section of $X$, Theorem~\ref{thm:nilpotent-COHA} recovers the construction of the ($\G_m$-equivariant) (motivic) cohomological and K-theoretical Hall algebras of nilpotent Higgs sheaves on $C$ (see \cite{Sala_Schiffmann, Porta_Sala_Hall}, and \cite{Minets_Higgs} for the rank zero case).
	Furthermore, the generation theorem given in \cite[Theorem 5.1]{Sala_Schiffmann} admits a natural interpretation in this setting.
	Namely, with the notations introduced in \cite[\S5]{Sala_Schiffmann}, it says that the subalgebra generated by the subspaces $\sfH_\ast(\Lambda_{(\alpha)})$ is dense in $\coha_{T^\ast C, C}$.
\end{remark}

\bigskip\section{0-dimensional nilpotent cohomological Hall algebra}\label{sec:0-dimensional}

In this section, we provide an explicit characterization of the nilpotent cohomological Hall algebra of 0-dimensional coherent sheaves on a smooth surface in terms of W-algebras following the approach in \cite{MMSV} of Mellit, Minets, the fourth and fifth-named authors, under certain assumption on the surface $S$ and the closed subscheme $Z$.

\medspace

Let $X$ be a smooth projective surface over a field $k$ of characteristic zero. We denote by $\dstackCoh_0(X; n)$ be the derived stack of length $n$ zero-dimensional sheaves on $X$.

Now, let $i\colon Z\to X$ be a connected closed reduced subscheme, and $U\coloneqq X\smallsetminus Z$. Let $\dstackCohnil_0(\widehat{X}_Z; n)$ be the derived stack of length $n$ zero-dimensional sheaves on $X$ set-theoretically supported on $Z$. Let $\Red{\dstackCohnil_0(\widehat{X}_Z; n)}$ denote its reduced geometric classical stack, which is a closed substack of $\dstackCoh_0(X; n)$. 

We impose the following condition on $X$ and $Z$.
\begin{assumption}\label{assum:zero-dimensional}
	The pushforward map $\HBMbullet(Z)\to \HBMbullet(X)$ is injective. 
\end{assumption}

\begin{remark}
	Note that Assumption~\ref{assum:zero-dimensional} holds for example if $X$ is a minimal resolution of an ADE singularity and $Z$ is the exceptional curve of if $X\to \PP^1$ is a smooth projective elliptic surface and $Z$ is a singular fiber of affine type ADE.
\end{remark}

We now record two consequences of Assumption~\ref{assum:zero-dimensional}. First, it yields the following.
\begin{lemma}\label{lem:cohU}
	There is a short exact sequence 
	\begin{align}
		0\longrightarrow \sfH^k(X,U) \longrightarrow \sfH^k(X) \longrightarrow \sfH^k(U) \longrightarrow 0 \ ,
	\end{align}
	where $\sfH^k(X,U)\simeq \HBM_k(Z)$ by Poincaré duality. 
\end{lemma}

\begin{proposition}\label{prop:Sym}
	We have\footnote{We thank Ben Davison for suggesting the second formula.}
	\begin{align}\label{eq:0-dim-X}
	\HBMbullet( \dstackCoh_0(X); \Q )& \simeq \Sym \left( \bigoplus_{n\geq 1}\Hbullet(X; \Q^{\mathsf{vir}}) \otimes \Hbullet(B\C^\ast)\right) \ , \\ \label{eq:0-dim-Z}
	\HBMbullet( \dstackCohnil_0(\widehat{X}_Z); \Q )& \simeq \Sym \left( \bigoplus_{n\geq 1}\Hbullet_Z(X; \Q^{\mathsf{vir}}) \otimes \Hbullet(B\C^\ast)\right) \ .
	\end{align}
	as cohomologically graded mixed Hodge modules. Here, $\Hbullet(X; \Q^{\mathsf{vir}})$ (resp.\ $\Hbullet_Z(X; \Q^{\mathsf{vir}})$) is obtained by applying the derived global sections functor to $\Q_X\otimes \LL^{-1}$ (resp.\ $i^! \Q_X\otimes \LL^{-1}$), where $\LL$ is the Tate twist.
\end{proposition}

\begin{proof}
	Formula~\eqref{eq:0-dim-X} is proved in \cite[Proposition~9.1]{DHSM22}.
	
	We now prove Formula~\eqref{eq:0-dim-Z}. We denote by $\Sym^n(X)$ the $n$-th symmetric product of $X$ for a fixed integer $n\geq 1$, while by $\Sym(X)$ the disjoint union of all $\Sym^n(X)$'s by varying of $n$. We introduce similarly $\Sym^n(Z)$ and $\Sym(Z)$. Note that we have a pullback square in the category of reduced stacks
	\begin{align}
		\begin{tikzcd}[ampersand replacement=\&]
			\trunc{\dstackCohnil_0(\widehat{X}_Z; n)}_\red \ar{d}{\varpi_{Z, n}} \ar{r}{j_{Z, n}}\& \trunc{\dstackCoh_0(X; n)} \ar{d}{\varpi_n}\\
			\Sym^n(Z) \ar{r}{i_{Z, n}} \& \Sym^n(X)
		\end{tikzcd}\ .
	\end{align}
	We drop the subscripts $n$ when considering all possible lengths at once. Moreover,we set $i\coloneqq i_{Z, 1}$.
	
	Denote by $p$ and $p_Z$ the maps from $\dstackCoh_0(X; n)$ and $\dstackCohnil_0(\widehat{X}_Z)$ to a point, while by $q$ and $q_Z$ the maps from $\Sym(X)$ and $\Sym(Z)$, respectively, to a point. Furthermore, denote by $\Delta_n\colon X\to \Sym^n(X)$ and $\Delta_{Z, n}\colon Z\to \Sym^n(Z)$ the inclusions of the small diagonals. We have 
	\begin{align}
		\CBM ( \dstackCohnil_0(\widehat{X}_Z); \Q ) & \simeq (p_Z)_\ast p_Z^!\Q \simeq (q_Z\circ \varpi_Z)_\ast j_Z^! p^!\Q\simeq (q_Z)_\ast  i_Z^! \varpi_\ast p^!\Q\\
		&\simeq (q_Z)_\ast  i_Z^! \Sym_\boxdot \left( \bigoplus_{n\geq 1}(\Delta_n)_\ast p_X^!\Q \otimes \Hbullet(B\C^\ast)\right) \\
		&\simeq (q_Z)_\ast   \Sym_\boxdot \left( \bigoplus_{n\geq 1}(\Delta_{Z, n})_\ast i^! p_X^!\Q \otimes \Hbullet(B\C^\ast)\right) \ ,
	\end{align}
	where the isomorphism in the second row follows from \cite[Proposition~9.1]{DHSM22}. By taking derived global sections, the conclusion follows.
\end{proof}

This proposition together with Assumption~\ref{assum:zero-dimensional} implies the following.
\begin{proposition}\label{prop:cohaZA}
	Let $\zeta\colon \Red{\dstackCohnil_0(\widehat{X}_Z; n)}\to \dstackCoh_0(X; n)$ be the canonical closed embedding \eqref{eq:closed-embedding}. Then, under Assumption~\ref{assum:zero-dimensional}, the pushforward map 
	\begin{align}\label{eq:ZtoS} 
		\zeta_\ast\colon \HBM_\bullet(\dstackCohnil_0(\widehat{X}_Z; n))\longrightarrow \HBM_\bullet(\dstackCoh_0(X; n))
	\end{align}
	is injective for all $n\geq 0$. 
\end{proposition}

\begin{proof}
	Since $\HBMbullet(Z)\to \HBMbullet(X)$ is injective by Assumption~\ref{assum:zero-dimensional}, also
	\begin{align}
		\bigoplus_{n\geq 1}\HBMbullet(Z) \otimes \Hbullet(B\C^\ast)\longrightarrow \bigoplus_{n\geq 1}\HBMbullet(X) \otimes \Hbullet(B\C^\ast)
	\end{align}
	is injective. This lifts to an injective morphism 
	\begin{align}
		\bigoplus_{n\geq 1}\Hbullet_Z( X; \Q^{\mathsf{vir}})\otimes \Hbullet(B\C^\ast)\longrightarrow \bigoplus_{n\geq 1}\Hbullet( X; \Q^{\mathsf{vir}}) \otimes \Hbullet(B\C^\ast)
	\end{align}
	of cohomologically graded mixed Hodge structures. By applying $\Sym(-)$ and using Proposition~\ref{prop:Sym} the result follows.
\end{proof}

As in \cite[\S1.3]{MMSV}, set 
\begin{align}
	h_X(z,w) \coloneqq \sum_{n \in \N}\sum_{m\in \Z} \dim \HBM_m(\dstackCoh_0(X; n))\, (-z)^m w^n\ .
\end{align}
Define also 
\begin{align}
	h_Z(z,w) \coloneqq \sum_{n \in \N}\sum_{m\in \Z} \dim \HBM_m(\dstackCohnil_0(\widehat{X}_Z; n))\, (-z)^m w^n\ .
\end{align}
Finally, set
\begin{align}
	h_Y(z) \coloneqq \sum_{m\in \Z} \dim \HBM_m(Y) (-z)^m\ .
\end{align}
for $Y=X, Z$. Proposition~\ref{prop:Sym} yields that
\begin{align}\label{eq:hilbseriesA}
	h_Y(z,w) = \mathsf{Exp} \left( \frac{h_Y(z) z^{-2} w}{(1-z^{-2})(1-w)}\right) 
\end{align}
for $Y=X,Z$, where $\mathsf{Exp}(-)$ is the plethystic exponential. 

Now, let $\coha_{X;\, 0}$ be the cohomological Hall algebra of $0$-dimensional sheaves on $X$, while let $\coha_{X, Z;\, 0}$ denote the cohomological Hall algebra of $\dstackCohpsnil(\widehat{X}_Z; 0)$. Recall that these are $\N\times \Z$-graded associative algebras. We denote by
\begin{align}
	\coha_{X;\, 0}[n, \ell] \quad \text{and} \quad \coha_{X, Z;\, 0}[n, \ell]
\end{align}
the $(n, \ell)$-graded pieces, respectively.

Under Assumption~\ref{assum:zero-dimensional}, by Proposition~\ref{prop:cohaZA} we have an injective map of associative algebras 
\begin{align}
	\coha_{X, Z;\, 0} \longrightarrow \coha_{X;\, 0}\ .
\end{align}
Furthermore, as constructed in \cite[\S7.3]{MMSV}, we have a representation 
\begin{align}
	\Psi_X^{-}\colon \coha_{X; \,0}\op\longrightarrow \End( \bfV(X) ) 
\end{align}
where 
\begin{align}
	\bfV(X)\coloneqq \bigoplus_{n\geq 0} \HBM_\bullet( \Hilb_n(X) )\ .
\end{align}
By restriction, we obtain an action 
\begin{align}
	\Psi_Z^- \colon \coha_{X, Z;\, 0} \op \longrightarrow  \End( \bfV(X) ) \ .
\end{align}

Let $u\colon U \to X$ be the canonical open immersion and let $u^\ast\colon \Hbullet(X) \to \Hbullet(U)$ be the associated restriction map. By Lemma~\ref{lem:cohU}, there is natural  graded isomorphism $\ker(u^\ast) \simeq \HBM_\bullet(Z)$.  

As shown in \cite[Proposition~7.8]{MMSV}, for any $\lambda \in \Hbullet(X)$ and $n>0$, there exists an element $E_n(\lambda)\in \HBM_\bullet(\dstackCoh_0(X;n))$ so that $\Psi^-(E_n(\lambda))$ coincides with the Nakajima operator $\frakq_n(\lambda)$. We have the following.
\begin{lemma}\label{lem:ZNoperators}
	If $\lambda \in \ker(u^\ast)$, then $E_n(\lambda) \in \mathsf{Im}(\zeta_\ast)$, where $\zeta\colon \Red{\dstackCohnil_0(\widehat{X}_Z; n)} \to \dstackCoh_0(X; n)$ is the canonical closed immersion \eqref{eq:closed-embedding}. 
\end{lemma}  

\begin{proof} 
	We first recall the construction of $E_n(\lambda)$ in \cite[Proposition~7.8]{MMSV}. Let $t\colon \dstackCoh_0^{\mathsf{pt}}(X; n)\to \dstackCoh_0(X; n)$ be the closed substack defined by the pull-back square 
	\begin{align}
		\begin{tikzcd}[ampersand replacement=\&]
			\dstackCoh_0^{\mathsf{pt}}(X; n)\ar{r}{t} \ar{d}{s}\& \dstackCoh_0(X; n)\ar{d} \\
			X \ar{r}{\delta}\& \Sym^n(S)
		\end{tikzcd}\ .
	\end{align} 
	Note that the main results of \cite{Irred_quot_I, Irred_quot_II} imply that the reduced stack $\Red{\dstackCoh_0^{\mathsf{pt}}(X; n)}$ is irreducible, of expected dimension. Let $[\dstackCoh_0^{\mathsf{pt}}(X; n)]\in \HBM_\bullet(\dstackCoh_0^{\mathsf{pt}}(X; n))$ denote the fundamental cycle of the reduced stack $\Red{\dstackCoh_0^{\mathsf{pt}}(X; n)}$. Then one sets 
	\begin{align}
		E_n(\lambda) \coloneqq t_\ast(s^\ast(\lambda) \cap [\dstackCoh_0^{\mathsf{pt}}(X; n)])\ . 
	\end{align}
	Let $v\colon \dstackCoh_0^{\mathsf{pt}}(U; n) \hookrightarrow \dstackCoh_0^{\mathsf{pt}}(X; n)$ be the open substack defined by the pull-back square 
	\begin{align}
		\begin{tikzcd}[ampersand replacement=\&]
			\dstackCoh_0^{\mathsf{pt}}(U; n)\ar{r}{v} \ar{d}{}\&  \dstackCoh_0^{\mathsf{pt}}(X; n)\ar{d}{s} \\
			U \ar{r}{u}\& X
		\end{tikzcd}\ .
	\end{align} 
	Let also $\eta\colon \dstackCoh^{\mathsf{pt}}_0(\widehat{X}_Z; n) \to  \dstackCohnil_0(\widehat{X}_Z; n)$ be the substack defined by the pull-back square 
	\begin{align}\label{eq:Zpt}
		\begin{tikzcd}[ampersand replacement=\&]
			\dstackCoh^{\mathsf{pt}}_0(\widehat{X}_Z; n)\ar{d}{\eta}\ar{r}{t_Z}
			\&   \dstackCohnil_0(\widehat{X}_Z; n) \ar{d}{\zeta} \\
			\dstackCoh_{0}^{\mathsf{pt}}(X; n)\ar{r}{t}
			\& \dstackCoh_0(X; n)
		\end{tikzcd}\ ,
	\end{align} 
	in which all squares are pull-back. Then note that $(\dstackCoh_0^{\mathsf{pt}}(U; n), \Red{\dstackCoh^{\mathsf{pt}}_0(\widehat{X}_Z; n)})$ 
	form a complementary pair in $\dstackCoh^{\mathsf{pt}}_0(X; n)$. 
	
	Now assume that $u^\ast(\lambda)=0$. This implies that $v^\ast s^\ast(\lambda)=0$, hence the Borel-Moore homology class $s^\ast(\lambda)\cap [\dstackCoh_0^{\mathsf{pt}}(X;n)]$ belongs to the image of 
	\begin{align}
		\eta \colon \Red{\dstackCoh^{\mathsf{pt}}_0(\widehat{X}_Z; n)}\longrightarrow \Red{\dstackCoh^{\mathsf{pt}}_0(X; n)}\ . 
	\end{align}
	Given diagram~\eqref{eq:Zpt}, this further implies that $E_n(\lambda)$ belongs to $\mathsf{Im}(\zeta_\ast)$, as claimed. 
\end{proof} 

Recall that in \cite[Definition~3.1]{MMSV}, the \textit{positive W-algebra $W^{\geqslant}(X)$ of $X$} and the \textit{negative W-algebra $W^{\leqslant}(X)$ of $X$} are introduced as, respectively, the $\N\times \Z$-graded associative algebras generated by $T_n^\pm(\lambda)$, for $n\in \N$ and $\lambda\in \Hbullet(X)$, $\psi_n(\lambda)$, with $n\in\N$ and $\lambda\in \Hbullet(X)$, and the central element $C$ subject to certain relations. Let $W^0(X)$ be the subalgebra generated by $\psi_n(\lambda)$, with $n\in\N$ and $\lambda\in \Hbullet(X)$, and the central element $C$, while let $W^\pm(X)$ be the subalgebra generated by $T_n^\pm(\lambda)$, with $n\in\N$ and $\lambda\in \Hbullet(X)$.
In \cite[Proposition~3.8]{MMSV}, the authors also introduced elements $D_{m,n}(\lambda)\in W^\geqslant(X)$ and $D_{-m,n}(\lambda)\in W^\leqslant(X)$, for $m, n\in \N$. These elements are relevant to introduce the following definition.
\begin{definition}[{\cite[Definition~5.1]{MMSV}}]
	Let $I\subset \Hbullet(X)$ be a graded ideal and let $J\coloneqq \Hbullet(X)/I$. Let $W^\pm(I)\coloneqq W^\pm_\downarrow(I)$ be the smallest graded subalgebra of $W^\pm(X)$ containing the elements $D_{\pm n,\,0}(\lambda)$ for all $n\in \N$ and $\lambda\in I$ and stable under the operators $\mathsf{Ad}(\psi_\ell(\mu))$ for all $\ell>0$ and $\mu\in \Hbullet(X)$.
	
	Let $W^0(J)$ be the quotient of $W^0(X)$ by the ideal generated by $\psi_\ell(\lambda)$ with $\ell>0$ and $\lambda\in I$.
\end{definition}
We shall apply the above definition to $I\coloneqq \Hbullet(X, U)\simeq \HBMbullet(Z)$ and $J\coloneqq \Hbullet(X)/I^\perp$, where $I^\perp$ coincides with the compactly supported singular cohomology $\sfH^\bullet_{\mathsf{c}}(U)$ of $U$: set 
\begin{align}
	W^\pm(\widehat{X}_Z)\coloneqq W^\pm(I)\quad \text{and}\quad W^0(\widehat{X}_Z)\coloneqq W^0(J)\ .
\end{align}
By \cite[Proposition~5.3]{MMSV}, we have
\begin{align}\label{eq:hilbWZ} 
	P_{W^+(\widehat{X}_Z)}(z,w)=\Exp\left(\frac{h_Z(z^{-1}) z^2 w}{(1-z^2) (1-w)}\right)\ ,
\end{align}
where the left-hand-side denotes the Hilbert series of $W^+(\widehat{X}_Z)$. 

Thanks to \cite[Lemma~5.5]{MMSV}, we can also define the subalgebras
\begin{align}
	W^\leqslant(\widehat{X}_Z)\coloneqq W^0(\widehat{X}_Z) \ltimes W^-(\widehat{X}_Z) \subset W^\leqslant(X) \quad \text{and} \quad W^\geqslant(\widehat{X}_Z)\coloneqq W^0(\widehat{X}_Z) \ltimes W^+(\widehat{X}_Z) \subset W^\geqslant(X)\ .
\end{align}

Now, we relate the 0-dimensional nilpotent COHA $\coha_{X, Z; \, 0}$ with W-algebras readapting the arguments in \cite[\S7.6]{MMSV}. Let $\Phi^-\colon W^{-}(X)\longrightarrow \End(\bfV(X))$ be the representation defined in \cite[Proposition~6.8]{MMSV} and let $\Phi^-_Z$ be its restriction to $W^{-}(\widehat{X}_Z)$. Then, given Lemma~\ref{lem:ZNoperators}, by analogy with \cite[Theorem~7.6]{MMSV}, we obtain the main result of this section.
\begin{theorem}\label{thm:cohaZB} 
	Under Assumption~\ref{assum:zero-dimensional}, the diagram 
	\begin{align} 
		\begin{tikzcd}[ampersand replacement=\&]
			W^\leqslant(\widehat{X}_Z)\ar{r}{\Phi^-_Z} \& \End(\bfV(X))\& \ar{l}[swap]{\Psi^-_Z} (\coha_{X, Z; \, 0}^{\mathsf{ex}})\op
		\end{tikzcd}
	\end{align} 
	induces an isomorphism of algebras 
	\begin{align}
		\coha_{X, Z; \, 0} \simeq W^+(\widehat{X}_Z)\ .
	\end{align}
	Here, $\coha_{X, Z; \, 0}^{\mathsf{ex}}$ is the \textit{extended} COHA by the action of the tautological ring $\bS_{X, Z}$\footnote{This action is defined as in \cite[\S\ref*{COHA-Yangian-subsec:coha-rsv-tautological}]{DPSSV-3}.}.
\end{theorem} 

\begin{proof} 
	By \cite[\S~7.3]{MMSV}, $\Phi^-$ is injective, hence also $\Phi^-_Z$. Given Lemma~\ref{lem:ZNoperators}, \cite[Proposition~7.9 and Formula~(7.7)]{MMSV} imply that $\mathsf{Im}(W_Z^-(X))\subset \mathsf{Im}(\Psi_Z^-)$ coincides with $\Psi_Z^-\left((\coha_{X, Z; \, 0})_{\mathsf{sph}}\op\right)$, where $(\coha_{X, Z; \, 0})_{\mathsf{sph}}$ is the subalgebra of $\coha_{X, Z; \, 0}$ generated by $\HBMbullet(\dstackCohnil_0(\widehat{X}_Z; 1))$. This yields the sequence of inequalities 
	\begin{align} 
		\dim W^+(\widehat{X}_Z)[n,\ell] \leq \dim \Psi_Z^-\left((\coha_{X, Z; \, 0})_{\mathsf{sph}}\op\right)
		& \leq \dim (\coha_{X, Z; \, 0})_{\mathsf{sph}}\op[n,\ell]\\
		& \leq \dim \coha_{X, Z; \, 0}\op[n,\ell]\ .
	\end{align} 
	Then Proposition~\ref{prop:cohaZA} and Formula~\eqref{eq:hilbWZ} imply that all inequalities in the above Formula are equalities. Thus
	\begin{align}
		(\coha_{X, Z; \, 0})_{\mathsf{sph}} = \coha_{X, Z; \, 0}
	\end{align}
	and the assertion follows.
\end{proof} 

\begin{remark}\label{rem:equivariant}
	Let $T$ be an arbitrary algebraic torus. A $T$-equivariant version of Assumption~\ref{assum:zero-dimensional} together with the assumption that both $X$ and $Z$ are $T$-equivariantly formal ensure that Theorem~\ref{thm:cohaZB} holds also $T$-equivariantly, following the same reasoning as in the paragraph after \cite[Theorem~1.5]{MMSV}. If $X$ and $Z$ are cohomologically pure, then they are $T$-equivariantly formal. 
\end{remark}

\appendix

\bigskip\section{Ind-objects}\label{appendix:ind-objects}

In this section, we collect some standard material on ind-objects for which we do not know written references.

\subsection{Restricted presentations for ind-morphisms}

Fix an $\infty$-category $\scrC$ and let $P$ be a property of morphisms in $\scrC$. Recall from \cite[Proposition~5.3.5.15]{HTT} that for every $\infty$-category $\scrC$, there is a canonical equivalence
\begin{align}
	\Ind\big(\scrC^{\Delta^1}\big) \simeq \Ind(\scrC)^{\Delta^1} \ . 
\end{align}
Reviewing morphisms in $\scrC$ as objects in $\scrC^{\Delta^1}$, \cite[Definition~\ref*{torsion-pairs-def:restricted_presentations}]{DPS_Torsion-pairs} allows to talk about \textit{ind-$P$} morphisms in $\Ind(\scrC)$; concretely, a morphism $f \colon \scrX \to \scrY$ in $\Ind(\scrC)$ is ind-$P$ if there exists a presentation $\{f_\alpha \colon X_\alpha \to Y_\alpha\}_{\alpha \in I}$ for the morphism $f$ such that for every index $\alpha \in I$ the map $f_\alpha \colon X_\alpha \to Y_\alpha$ satisfies the property $P$.

\begin{example}\label{eg:representable_morphisms}
	\hfill
	\begin{enumerate}\itemsep=0.2cm
		\item \label{item:representable_morphisms-1} If the morphism $f \colon \scrX \to \scrY$ is representable by $P$-morphisms, in the sense that for every $Y \in \scrC$ and every morphism $Y \to \scrY$ in $\Ind(\scrC)$ the fiber product $Y \times_{\scrY} \scrX$ belongs to $\scrC$ and the canonical map $Y \times_{\scrY} \scrX \to Y$ satisfies the property $P$, then $f$ is ind-$P$.
		
		\item \label{item:representable_morphisms-2} Let $f \colon X \to Y$ be a morphism in $\scrC$ which, seen as a morphism in $\Ind(\scrC)$, is ind-$P$.
		Then $f$ is a retract of a $P$-morphism.
		In particular, if $P$-morphisms are stable under retracts, we see that a morphism in $\scrC$ satisfies $P$ if and only if it is ind-$P$ as a morphism in $\Ind(\scrC)$. \qedhere
	\end{enumerate}
\end{example}

\begin{definition}
	Let $\scrA$ be an $\infty$-category and let $F \colon \scrA \to \scrC$ be a diagram.
	We say that $F$ is a \textit{$P$-diagram} if it takes every arrow in $\scrA$ to $P$-morphisms in $\scrC$.
\end{definition}

The following result establishes the existence of \textit{simultaneous presentations} for specific diagrams.

\begin{proposition}\label{prop:simultaneous_presentations}
	Assume that $\scrC$ has finite limits and let $P$ be a property of morphisms in $\scrC$ which is stable under fiber products.
	For $k \in \{1,2\}$, let $\scrA \coloneqq \Lambda^2_k$ be the corresponding $2$-horn.
	For a diagram $F \colon \scrA \to \Ind(\scrC)$, the following statements are equivalent:
	\begin{enumerate}\itemsep=0.2cm
		\item \label{item:simultaneous_presentations-1} $F$ is an ind-$P$-diagram;
		
		\item \label{item:simultaneous_presentations-2} there exists a filtered diagram $I \to \Fun(\scrA, \scrC)_{/F}$, sending $\alpha \in \scrA$ to a $P$-diagram $F_\alpha$ such that the canonical map
		\begin{align}
			\colim_\alpha F_\alpha \longrightarrow F 
		\end{align}
		is an equivalence in $\Ind(\scrC)$.
	\end{enumerate}
\end{proposition}

\begin{proof}
	We deal with the case $k = 1$.
	The case $k = 2$ is similarly dealt with.
	First of all, the implication \eqref{item:simultaneous_presentations-2} $\Rightarrow$ \eqref{item:simultaneous_presentations-1} is trivial.
	Let us prove the converse.
	Let therefore $F$ be an ind-$P$-diagram, which we represent as
	\begin{align}
		\begin{tikzcd}[column sep = 15pt, ampersand replacement=\&]
			\scrX \arrow{r}{f} \& \scrY \arrow{r}{g} \& \scrZ
		\end{tikzcd} \ ,
	\end{align}
	where both $f$ and $g$ are ind-$P$ morphisms.
	Using \cite[Proposition 5.3.5.15]{HTT}, we find a canonical equivalence $\Ind(\Fun(\Lambda^2_1, \scrC)) \simeq \Fun(\Lambda^2_1, \Ind(\scrC))$.
		Applying \cite[Lemma~\ref*{torsion-pairs-lem:Q_presentations}]{DPS_Torsion-pairs}, we are reduced to check that the following factorization problem has always a solution: for every solid commutative diagram
	\begin{align}\label{eq:factorization_problem}
		\begin{tikzcd}[ampersand replacement = \&]
			X_0 \arrow{rr}{f_0} \arrow{dd} \arrow[dashed]{dr} \& \& Y_0 \arrow{rr}{g_0} \arrow[dashed]{dr} \& \& Z_0 \arrow[dashed]{dr} \\
			\& \overline{X} \arrow[dashed]{dl} \arrow{rr}[near start]{\overline{f}} \& \& \overline{Y} \arrow[dashed]{dl} \arrow{rr}[near start]{\overline{g}} \& \& \overline{Z} \arrow[dashed]{dl} \\
			\scrX \arrow{rr}{f} \& \& \scrY \arrow{rr}{g} \arrow[leftarrow, crossing over]{uu} \& \& \scrZ \arrow[leftarrow, crossing over]{uu}
		\end{tikzcd}\ ,
	\end{align}
	where the top line is a compact object in $\Fun(\Lambda^2_1, \Ind(\scrC))$ (that is, $X_0$, $Y_0$ and $Z_0$ are compact in $\Ind(\scrC)$), there exists a $P$-diagram $\overline{F} \colon \Lambda^2_1 \to \scrC$ (represented as the middle line in the above diagram) and a factorization as indicated.
	Since $g$ is ind-$P$, we can fix a $P$-presentation as filtered colimit of a diagram $\{g_\alpha \colon Y_\alpha \to Z_\alpha\}_{\alpha \in J}$ of $P$-morphisms.
	Since $g_0 \colon Y_0 \to Z_0$ is compact in $\Ind(\scrC)^{\Delta^1}$, we deduce the existence of an index $\alpha$ and a factorization of the front right square in \eqref{eq:factorization_problem} as
	\begin{align}
		\begin{tikzcd}[ampersand replacement=\&]
			Y_0 \arrow{d}{g_0} \arrow{r} \& Y_\alpha \arrow{d}{g_\alpha} \arrow{r} \& \scrY \arrow{d}{g} \\
			Z_0 \arrow{r} \& Z_\alpha \arrow{r}{f_\alpha} \& \scrZ
		\end{tikzcd} \ .
	\end{align}
	On the other hand, $f \colon \scrX \to \scrY$ is also an ind-$P$-morphism.
	Thus, we can fix a $P$-presentation $\{f_\alpha \colon X_\beta \to \widetilde{Y}_\beta\}_{\beta \in J'}$.
	Since source and target of the composite morphism $f_0' \colon X_0 \to Y_\alpha$ are compact in $\Ind(\scrC)$, we deduce the existence of an index $\beta \in J'$ and of a factorization
	\begin{align}
		\begin{tikzcd}[ampersand replacement=\&]
			X_0 \arrow{r} \arrow{d}{f_0'} \& X_\beta \arrow{d}{f_\beta} \arrow{r} \& \scrX \arrow{d}{f} \\
			Y_\alpha \arrow{r} \& \widetilde{Y}_\beta \arrow{r} \& \scrY
		\end{tikzcd}  \ .
	\end{align}
	Set $X_{\alpha,\beta} \coloneqq Y_\alpha \times_{\widetilde{Y}_\beta} X_\beta$.
	Since $\scrC$ has finite limits, this is an object in $\scrC$; moreover since $P$-morphisms are stable under pullbacks, the morphism $f_{\alpha,\beta} \colon X_{\alpha,\beta} \to Y_\alpha$ satisfies $P$.
	Combining what we have found so far, we see that the diagram
	\begin{align}
		\begin{tikzcd}[column sep = 20pt, ampersand replacement=\&]
			X_{\alpha,\beta} \arrow{r}{f_{\alpha,\beta}} \& Y_\alpha \arrow{r}{g_\alpha} \& Z_\alpha
		\end{tikzcd}
	\end{align}	
	solves the factorization problem \eqref{eq:factorization_problem}.
\end{proof}

\begin{variant}\label{variant:simultaneous_presentations}
	In the setting of the previous proposition, let $P$ and $Q$ be two different properties of morphisms of $\scrC$, and assume that at least one of the two is stable under fiber products.
	Say that a diagram $F \colon \Lambda^2_1 \to \scrC$ depicted as
	\begin{align}
		\begin{tikzcd}[column sep = 20pt, ampersand replacement=\&]
			X \arrow{r}{f} \& Y \arrow{r}{g} \& Z
		\end{tikzcd} 
	\end{align}
	is a $(P,Q)$-diagram if $f$ satisfies $P$ and $g$ satisfies $Q$.
	Then for a diagram
	\begin{align}
		\begin{tikzcd}[column sep = 20pt, ampersand replacement=\&]
			\scrX \arrow{r}{f} \& \scrY \arrow{r}{g} \& \scrZ
		\end{tikzcd}
	\end{align} 
	in $\Ind(\scrC)$ the following statements are equivalent:
	\begin{enumerate}\itemsep=0.2cm
		\item $f$ is ind-$P$ and $g$ is ind-$Q$;
		
		\item the above diagram admits a presentation as filtered colimit of $\scrC$-valued $(P,Q)$-diagrams.
	\end{enumerate}
	A similar statement holds replacing $\Lambda^2_1$ by $\Lambda^2_2$.
\end{variant}

\begin{remark}\label{rem:pushout}
	If $\scrC$ admits finite colimits and $P$-morphisms are closed under pushout, then a dual argument shows that in the statements of Proposition~\ref{prop:simultaneous_presentations} and Variant~\ref{variant:simultaneous_presentations} one can replace $\Lambda^2_2$ by $\Lambda^2_0$.
	However, these assumptions will not be satisfied in our framework.
\end{remark}

\begin{corollary}\label{cor:ind_morphisms_properties}
	Let $\scrC$ be an $\infty$-category and let $P$ be a property of morphisms in $\scrC$.
	Then:
	\begin{enumerate}\itemsep=0.2cm
		\item \label{item:ind_morphisms_properties-1} If for every $X \in \scrC$ the identity $\id_X$ satisfies $P$, then for every $\scrX \in \Ind(\scrC)$ the identity $\id_{\scrX}$ is ind-$P$.
	\end{enumerate}
	Assume furthermore that $\scrC$ has finite limits and that $P$-morphisms are stable under pullback.
	Then:
	\begin{enumerate}\setcounter{enumi}{1}\itemsep=0.2cm
		\item  \label{item:ind_morphisms_properties-2} Ind-$P$-morphisms are stable under pullbacks.
		
		\item  \label{item:ind_morphisms_properties-3} If $P$-morphisms are closed under composition, then the same goes for ind-$P$-morphisms.
		
		\item  \label{item:ind_morphisms_properties-4} Assume that $P$-morphism have the following closure property: given a commutative triangle
		\begin{align}
			\begin{tikzcd}[column sep = small, ampersand replacement=\&]
				X \arrow{rr}{f} \arrow{dr}[swap]{h} \& \& Y \arrow{dl}{g} \\
				{} \& Z 
			\end{tikzcd} \ ,
		\end{align}
		then if $g$ satisfies $P$, then $f$ satisfies $P$ if and only if $h$ satisfies $P$.
		Then ind-$P$ morphism have the same closure property.
	\end{enumerate}
\end{corollary}

\begin{proof}
	If $\id_X$ satisfies $P$ for every $X \in \scrC$, then it is obvious that $\id_{\scrX}$ is ind-$P$ for every $\scrX \in \Ind(\scrC)$.
	This proves  \eqref{item:ind_morphisms_properties-1}.
	
	For  \eqref{item:ind_morphisms_properties-2}, let
	\begin{align}
		\begin{tikzcd}[ampersand replacement=\&]
			\scrW \arrow{r} \arrow{d} \& \scrX \arrow{d}{f} \\
			\scrZ \arrow{r}{g} \& \scrY
		\end{tikzcd} 
	\end{align}
	be a pullback square in $\Ind(\scrC)$ where $f$ is in ind-$P$.
	Applying Variant~\ref{variant:simultaneous_presentations} (with $Q = \mathsf{all}$ the class of all morphisms of $\scrC$), we can find a filtered diagram 
	\begin{align}
		\begin{tikzcd}[ampersand replacement=\&]
			\Big\{Z_\alpha \ar{r}{g_\alpha} \& Y_\alpha \ar[leftarrow]{r}{f_\alpha} \& X_\alpha\Big\}_{\alpha \in I}
		\end{tikzcd}
	\end{align}
	whose colimit coincides with 
	\begin{align}
		\begin{tikzcd}[ampersand replacement=\&]
			\scrZ \ar{r}{g}\& \scrY \ar[left]{r}{f} \& \scrX
		\end{tikzcd}
	\end{align}
	and for which $f_\alpha$ satisfies $P$ for every index $\alpha \in I$.
	Thus, the map $\scrW \to \scrZ$ can be presented as the colimit of the $P$-morphisms $Z_\alpha \times_{Y_\alpha} X_\alpha \to Z_\alpha$, which satisfy $P$ by assumptions.
	
	We now prove \eqref{item:ind_morphisms_properties-3}.
	Let $f \colon \scrX \to \scrY$ and $g \colon \scrY \to \scrZ$ be two ind-$P$ morphisms.
	Applying Proposition~\ref{prop:simultaneous_presentations} with $k = 1$, we can find a simultaneous presentation 
	\begin{align}
		\begin{tikzcd}[ampersand replacement=\&]
			\Big\{ X_\alpha \ar{r}{f_\alpha}  \&Y_\alpha \ar{r}{g_\alpha} \& Z_\alpha\Big\}_{\alpha \in I}
		\end{tikzcd}
	\end{align}
	for the two morphisms $f$ and $g$ such that $f_\alpha$ and $g_\alpha$ both satisfy $P$ for every index $\alpha \in I$.
	The conclusion is therefore obvious.
	
	For  \eqref{item:ind_morphisms_properties-4}, we proceeds similarly, applying Proposition~\ref{prop:simultaneous_presentations} with $k = 2$.
\end{proof}

\subsection{Restricted presentations for ind-objects}

In this section, we discuss presentations of ind-objects, following \cite[\S\ref*{torsion-pairs-appendix:ind_objects}]{DPS_Torsion-pairs}.

We fix as usual an $\infty$-category $\scrC$ and a property $P$ of morphisms in $\scrC$. We make the following global assumption:

\begin{assumption}\label{assumption:P_morphisms}
	Let $\scrC$ be an $\infty$-category with finite limits and let $P$ be a property of morphisms in $\scrC$ such that:
	\begin{enumerate}\itemsep=0.2cm
		\item \label{item:P_morphisms-1} every identity of $\scrC$ satisfies $P$;
		
		\item \label{item:P_morphisms-2} $P$-morphisms are stable under composition and retracts;
		
		\item \label{item:P_morphisms-3} given morphisms $f \colon X \to Y$ and $g \colon Y \to Z$ such that $g$ satisfies $P$, then $f$ satisfies $P$ if and only if $g \circ f$ satisfies $P$.
	\end{enumerate}
\end{assumption}

Notice that \eqref{item:P_morphisms-1} guarantees that every object $X \in \scrC$ is ind-$P$.

\begin{remark}\label{rem:P_presentations_objects_and_morphisms}
	Under the above assumption, let $\scrC_P$ be the associated (non full) subcategory of $\scrC$ and consider the induced functor $\Ind(\scrC_P) \to \Ind(\scrC)$.
	Unraveling the definition, we see that an object belongs to the essential image of this functor if and only if it is ind-$P$.
	Moreover, this functor take every morphism to an ind-$P$ morphism.
	In other words, it factors through the non-full subcategory $\Ind_P(\scrC)_{\ind\textrm{-}P}$.
\end{remark}

\begin{lemma}\label{lem:morphism_restricted_presentations}
	Let $\scrX \in \Ind(\scrC)$ be an ind-object and let $\{X_\alpha\}_{\alpha \in I}$ be a presentation for $\scrX$.
	The following statements are equivalent:
	\begin{enumerate}\itemsep=0.2cm
		\item $\{X_\alpha\}_{\alpha \in I}$ is a $P$-presentation;
		
		\item for every $\alpha \in I$, the morphism $X_\alpha \to \scrX$ is an ind-$P$ morphism;
	\end{enumerate}
\end{lemma}

\begin{proof}
	First assume that $\{X_\alpha\}_{\alpha \in I}$ is a $P$-presentation.
	Then for every $\alpha \in I$, the map $X_\alpha \to \scrX$ can be written as colimit of the maps $\{X_\alpha \to X_\beta\}_{\beta \in I_{\alpha/}}$.
	Since every transition map satisfies $P$, we conclude that $X_\alpha \to \scrX$ is an ind-$P$ morphism.
	
	Conversely, suppose that for every $\alpha \in I$, the map $X_\alpha \to \scrX$ is ind-$P$.
	Let $\alpha \to \beta$ be a morphism in $\scrC$, and consider the commutative triangle
	\begin{align}
		\begin{tikzcd}[column sep = small, ampersand replacement=\&]
			X_\alpha \arrow{rr} \arrow{dr} \& \& X_\beta \arrow{dl} \\
			{}\& \scrX 
		\end{tikzcd}\ .
	\end{align}
	Since $P$ is stable under pullbacks and satisfies the $2$-out-of-$3$ property, Corollary~\ref{cor:ind_morphisms_properties}--\eqref{item:ind_morphisms_properties-4} guarantees that the map $X_\alpha \to X_\beta$ is ind-$P$.
	Since $P$-morphisms are closed under retracts, Example~\ref{eg:representable_morphisms}--\eqref{item:representable_morphisms-2} guarantees that $X_\alpha \to X_\beta$ satisfies $P$ as well.
	Thus, $\{X_\alpha\}_{\alpha \in P}$ is a $P$-presentation.
\end{proof}

\begin{corollary}\label{cor:morphism_restricted_presentations}
	For an ind-object $\scrX \in \Ind(\scrC)$, the following statements are equivalent:
	\begin{enumerate}\itemsep=0.2cm
		\item \label{item:morphism_restricted_presentations-1} $\scrC$ is an ind-$P$ object;
		\item \label{item:morphism_restricted_presentations-2} the full subcategory $\scrC_{/\scrX}^{\ind\textrm{-}P}$ of $\scrC \times_{\Ind(\scrC)} \Ind(\scrC)_{/\scrX}$ spanned by ind-$P$ morphisms $X \to \scrX$ is filtered and the inclusion $\scrC_{/\scrX}^{\ind\textrm{-}P} \hookrightarrow \scrC_{/\scrX}$ is colimit-cofinal.
	\end{enumerate}
\end{corollary}

\begin{proof}
	The implication \eqref{item:morphism_restricted_presentations-2} $\Rightarrow$ \eqref{item:morphism_restricted_presentations-1} follows directly from Lemma~\ref{lem:morphism_restricted_presentations}.
	
	Conversely, fix a $P$-presentation $\{X_\alpha\}_{\alpha \in I}$.
	Reasoning as in \cite[Lemma~\ref*{torsion-pairs-lem:Q_presentations}]{DPS_Torsion-pairs}, we reduce to the following factorization statement: for every compact object $\overline{X} \in \Ind(\scrC)$ and every morphism $f \colon \overline{X} \to \scrX$, there exists a factorization of $f$ as
	\begin{align}
		\begin{tikzcd}[column sep = small, ampersand replacement=\&]
			\overline{X} \arrow{r}{f'} \& X' \arrow{r}{f''} \& \scrX
		\end{tikzcd} \ ,
	\end{align}
	where $X' \in \scrC$ and $f''$ is ind-$P$.
	Since $\overline{X}$ is compact, there exists $\alpha \in I$ such that $f$ factors through the structural map $X_\alpha \to \scrX$.
	On the other hand, Lemma~\ref{lem:morphism_restricted_presentations} guarantees that this map is ind-$P$.
	Thus, the conclusion follows.
\end{proof}

\begin{proposition}\label{prop:simultaneous_presentations_II}
	Under Assumption~\ref{assumption:P_morphisms}, the induced functor
	\begin{align}\label{eq:P_presentations_objects_and_morphisms}
		\Ind(\scrC_P) \to \Ind_P(\scrC)_{\ind\textrm{-}P}
	\end{align}
	is an equivalence.
	In particular, if $f \colon \scrX \to \scrY$ is an ind-$P$ morphism between ind-$P$ objects, then there exists a presentation $\{f_\alpha \colon X_\alpha \to Y_\alpha\}_{\alpha \in I}$ such that for every $\alpha \in I$ the morphism $f_\alpha$ satisfies $P$ and both $\{X_\alpha\}_{\alpha \in I}$ and $\{Y_\alpha\}_{\alpha \in I}$ are $P$-presentations for $\scrX$ and $\scrY$, respectively.
\end{proposition}

\begin{proof}
	In virtue of Remark~\ref{rem:P_presentations_objects_and_morphisms}, we see that the functor under consideration is essentially surjective.
	It is therefore enough to check that it is fully faithful.
	Observe that the inclusion $\scrC_P \hookrightarrow \scrC$ is faithful, that is for every $X, Y \in \scrC_P$, the morphism
	\begin{align}
		\Map_{\scrC_P}(X,Y) \longrightarrow \Map_{\scrC}(X,Y) 
	\end{align}
	is $(-1)$-truncated. 
	Let now $\scrX$ and $\scrY$ be two objects in $\Ind(\scrC_P)$ and fix two $P$-presentations $\{X_\alpha\}_{\alpha \in I}$ and $\{Y_\beta\}_{\beta \in J}$.
	Then we have
	\begin{align}
		\Map_{\Ind(\scrC_P)}(\scrX, \scrY) \simeq \lim_\alpha \colim_\beta \Map_{\scrC_P}(X_\alpha, Y_\beta) 
	\end{align}
	and
	\begin{align}
		\Map_{\Ind(\scrC)}(\scrX, \scrY) \simeq \lim_\alpha \colim_\beta \Map_{\scrC}(X_\alpha, Y_\beta) \ . 
	\end{align}
	Since $(-1)$-truncated morphisms are stable under arbitrary limits and under filtered colimits, it follows that the map
	\begin{align}
		\Map_{\Ind(\scrC_P)}(\scrX, \scrY) \longrightarrow \Map_{\Ind(\scrC)}(\scrX, \scrY) 
	\end{align}
	is $(-1)$-truncated as well, i.e.\ $\Ind(\scrC_P) \to \Ind(\scrC)$ is a faithful functor.
	Since the same is true for $\Ind(\scrC)_{\ind\textrm{-}P} \hookrightarrow \Ind(\scrC)$, we deduce that the functor \eqref{eq:P_presentations_objects_and_morphisms} is faithful as well.
	To complete the proof, it is therefore enough to check that the functor in consideration is full.
	
	\medskip
	
	Let therefore $f \colon \scrX \to \scrY$ be an ind-$P$ morphism between ind-$P$ objects.
	Let $(\scrC^{\Delta^1})_{/f}^{\ind\textrm{-}P}$ be the full subcategory of $(\scrC^{\Delta^1})_{/f} \coloneqq \scrC^{\Delta^1} \times_{\Ind(\scrC^{\Delta^1})} \Ind(\scrC^{\Delta^1})_{/f}$ spanned by those squares
	\begin{align}
		\begin{tikzcd}[ampersand replacement=\&]
			X \arrow{d}{g} \arrow{r} \& \scrX \arrow{d}{f} \\
			Y \arrow{r} \& \scrY
		\end{tikzcd} 
	\end{align}
	whose horizontal morphisms are ind-$P$.
	Since ind-$P$ morphism satisfy the $2$-out-of-$3$ property and since $P$-morphisms are closed under retracts, it automatically follows in this situation that $g \colon X \to Y$ satisfies $P$.
	It follows that the forgetful functor
	\begin{align}
		(\scrC^{\Delta^1})_{/f}^{\ind\textrm{-}P} \longrightarrow \scrC^{\Delta^1}
	\end{align}
	can be canonically factored through $(\scrC_P)^{\Delta^1}$.
	Thus, it is enough to check that $(\scrC^{\Delta^1})_{/f}^{\ind\textrm{-}P}$ is filtered and that the inclusion
	\begin{align}
		(\scrC^{\Delta^1})_{/f}^{\ind\textrm{-}P} \hookrightarrow (\scrC^{\Delta^1})_{/f} 
	\end{align}
	is colimit-cofinal.
	As usual, it is enough to show that the following factorization problem can always be solved: for every solid diagram
	\begin{align}
		\begin{tikzcd}[ampersand replacement=\&]
			\overline{X} \arrow{dd}{g} \arrow{rr} \arrow[dashed]{dr} \& \& \scrX \arrow{dd}{f} \\
			\& X' \arrow[dashed]{ur}{u} \\
			\overline{Y} \arrow{rr} \arrow[dashed]{dr} \& \& \scrY  \\
			\& Y' \arrow[dashed]{ur}{v} \arrow[leftarrow, crossing over]{uu}[fill = white,swap]{g'}
		\end{tikzcd}\ ,
	\end{align}
	where $\overline{X}$ and $\overline{Y}$ are compact objects in $\Ind(\scrC)$, there exists a factorization as indicated in the diagram, where furthermore $u$ and $v$ are ind-$P$.
	To see this, start by choosing a $P$-presentation $\{X_\alpha\}_{\alpha \in I}$ for $\scrX$.
	Since $\overline{X}$ is compact, we can find $\alpha \in I$ such that the map $\overline{X} \to \scrX$ factors through $X_\alpha \to \scrX$.
	Fix at the same time a $P$-presentation $\{Y_\beta\}_{\beta \in J}$ for $\scrY$.
	Since both $\overline{Y}$ and $X_\alpha$ are compact in $\Ind(\scrC)$, we can find a factorization of the original square as
	\begin{align}
		\begin{tikzcd}[ampersand replacement=\&]
			\overline{X} \arrow{r} \arrow{d} \& X_\alpha \arrow{r}{u} \arrow{d} \& \scrX \arrow{d} \\
			\overline{Y} \arrow{r} \& Y_\beta \arrow{r}{v} \& \scrY
		\end{tikzcd} \ .
	\end{align}
	Since we started with $P$-presentations for $\scrX$ and $\scrY$, Lemma~\ref{lem:morphism_restricted_presentations} guarantees that $u$ and $v$ are ind-$P$ morphisms, whence the conclusion.
\end{proof}

Remarkably, this result allows to improve Corollary~\ref{cor:ind_morphisms_properties}--\eqref{item:ind_morphisms_properties-4} without assuming the existence of pushouts in $\scrC$ (cf.\ Remark~\ref{rem:pushout}).

\begin{corollary}\label{cor:two_out_of_three}
	Assume in addition to Assumption~\ref{assumption:P_morphisms} that $P$-morphisms satisfy the $2$-out-of-$3$ property.
	Then the same goes for ind-$P$ morphisms.
\end{corollary}

\begin{proof}
	Let
	\begin{align}\label{eq:two_out_of_three}
		\begin{tikzcd}[column sep = small,ampersand replacement = \&]
			\scrX \arrow{rr}{f} \arrow{dr}[swap]{h} \& \& \scrY \arrow{dl}{g} \\
			{} \& \scrZ
		\end{tikzcd}
	\end{align}
	be a commutative triangle in $\Ind(\scrC)$.
	Assume that $\scrX$, $\scrY$ and $\scrZ$ are ind-$P$.
	If $g$ is ind-$P$, then Corollary~\ref{cor:ind_morphisms_properties}--\eqref{item:ind_morphisms_properties-4} shows that $f$ is ind-$P$ if and only if $h$ is ind-$P$.
	We are therefore left to check that if both $h$ and $f$ are ind-$P$, then the same goes for $g$.
	To do this, we apply \cite[Lemma~\ref*{torsion-pairs-lem:Q_presentations}]{DPS_Torsion-pairs} with $\scrE = \scrC^{\Delta^2}$, saying a $2$-simplex of the form \eqref{eq:two_out_of_three} in $\scrC$ satisfies $Q$ if $f$, $g$ and $h$ satisfy the property $P$.
	We are therefore led to consider the following factorization problem: given a morphism $\phi$ in $\Ind(\scrC^{\Delta^2})$
	\begin{align}\label{eq:factorization_problem_II}
		\begin{tikzcd}[column sep = small,ampersand replacement =\&]
			\& \& \& \overline{X} \arrow{dlll}[swap,xshift = -1.5ex, yshift = -0.75ex]{\phi_0} \arrow{rr} \arrow{dr} \& \& \overline{Y} \arrow{dl} \arrow[crossing over]{dlll}[swap,xshift = -1.5ex, yshift = -0.75ex]{\phi_1} \\
			\scrX \arrow{rr}[xshift = 0.5ex]{f} \arrow{dr}[swap,near start]{h} \& \& \scrY \arrow{dl}[xshift = -1.5ex, yshift = -0.75ex]{g} \& \& \overline{Z} \arrow{dlll}{\phi_2} \\
			\& \scrZ 
		\end{tikzcd}
	\end{align}
	where $\overline{X}$, $\overline{Y}$ and $\overline{Z}$ are compact in $\Ind(\scrC)$, we can factor $\alpha$ through a $2$-simplex in $\scrC$ satisfying $Q$.
	Since $f$ is ind-$P$, Proposition~\ref{prop:simultaneous_presentations_II} provides a presentation $\{f_\alpha \colon X_\alpha \to Y_\alpha\}_{\alpha \in I}$ such that every $f_\alpha$ satisfies $P$ and for which the structural morphisms $X_\alpha \to \scrX$ and $Y_\alpha \to \scrY$ are ind-$P$.
	Compactness of $\overline{X}$ and $\overline{Y}$ implies that there exists an index $\alpha \in I$ such that the morphism $\phi$ factors through the $2$-simplex
	\begin{align}
		\begin{tikzcd}[column sep = small, ampersand replacement=\&]
			X_\alpha \arrow{rr}{f_\alpha} \arrow{dr}[swap, near start]{h_\alpha} \& \& Y_\alpha \arrow{dl}[near start]{g_\alpha} \\
			{} \& \scrZ 
		\end{tikzcd} \ .
	\end{align}
	Fix now a $P$-presentation $\{Z_\beta\}_{\beta \in J}$ for $\scrZ$.
	Compactness of $X_\alpha$, $Y_\alpha$, $\overline{X}$, $\overline{Y}$ and $\overline{Z}$ implies therefore that there exists an index $\beta \in J$ and a commutative diagram
	\begin{align}
		\begin{tikzcd}[column sep = small,ampersand replacement =\&]
			\& \& \& \overline{X} \arrow{dlll} \arrow{rr} \arrow{dr} \& \& \overline{Y} \arrow{dl} \arrow[crossing over]{dlll} \\
			X_\alpha \arrow{rr}{f_\alpha} \arrow{dr}[swap,near start]{h_{\alpha,\beta}} \& \& Y_\alpha \arrow{dl}[near start,fill=white]{g_{\alpha,\beta}} \& \& \overline{Z} \arrow{dlll} \\
			\& Z_\beta
		\end{tikzcd}
	\end{align}
	factorizing the original diagram \eqref{eq:factorization_problem_II}.
	By construction, $f_\alpha$ satisfies $P$.
	To conclude, it is enough to argue that both $h_{\alpha,\beta}$ and $g_{\alpha,\beta}$ satisfy $P$.
	Since $P$-morphism satisfy the $2$-out-of-$3$ property, it is enough to argue that $h_{\alpha,\beta}$ satisfies $P$.
	However, the diagram
	\begin{align}
		\begin{tikzcd}[column sep = small, ampersand replacement=\&]
			X_\alpha \arrow{rr}{h_{\alpha,\beta}} \arrow{dr}[swap,near start]{h_{\alpha}} \& \& Z_\beta \arrow{dl}[near start]{j_\beta} \\
			{} \& \scrZ
		\end{tikzcd}
	\end{align}
	is commutative, and both $h_\alpha$ and $j_\beta$ are ind-$P$.
	Thus, Corollary~\ref{cor:ind_morphisms_properties}--\eqref{item:ind_morphisms_properties-4} guarantees that $h_{\alpha,\beta}$ satisfies $P$.
	The conclusion follows.
\end{proof}

\begin{corollary}\label{cor:indP_objects_fiber_products}
	Assume in addition to Assumption~\ref{assumption:P_morphisms} that $\scrC$ has finite limits and that $P$ is stable under fiber products.
	Then $\Ind_P(\scrC)$ is closed under fiber products in $\scrC$.
\end{corollary}

\begin{proof}
	Consider a span
	\begin{align}
		\scrX \stackrel{f}{\longrightarrow} \scrY \stackrel{g}{\longleftarrow} \scrZ 
	\end{align}
	in $\Ind_P(\scrC)$.
	We claim that whenever given a solid commutative diagram
	\begin{align}\label{eq:factorization_problem_III}
		\begin{tikzcd}[ampersand replacement = \&]
			X_0 \arrow{rr}{f_0} \arrow{dd} \arrow[dashed]{dr} \& \& Y_0 \arrow[dashed]{dr} \& \& Z_0 \arrow[dashed]{dr} \arrow{ll}[swap]{g_0} \\
			\& \overline{X} \arrow[dashed]{dl} \arrow{rr}[near start]{\overline{f}} \& \& \overline{Y} \arrow[dashed]{dl} \& \& \overline{Z} \arrow[dashed]{dl} \arrow{ll}[near end,swap]{\overline{g}} \\
			\scrX \arrow{rr}{f} \& \& \scrY \arrow[leftarrow, crossing over]{uu} \& \& \scrZ \arrow[leftarrow, crossing over]{uu} \arrow{ll}[swap]{g}
		\end{tikzcd}\ ,
	\end{align}
	where $X_0$, $Y_0$ and $Z_0$ are compact in $\Ind(\scrC)$, there exists the dashed factorization through a span in $\scrC$ with the property that the morphisms $\overline{X} \to \scrX$, $\overline{Y} \to \scrY$ and $\overline{Z} \to \scrZ$ are ind-$P$.
	Combining \cite[Lemma~\ref*{torsion-pairs-lem:Q_presentations}]{DPS_Torsion-pairs} and Lemma~\ref{lem:morphism_restricted_presentations}, this implies that the original span admits a simultaneous presentation
	\begin{align}
		\{ X_\alpha \stackrel{f_\alpha}{\longrightarrow} Y_\alpha \stackrel{g_\alpha}{\longleftarrow} Z_\alpha \}_{\alpha \in I} 
	\end{align}
	which restricts to $P$-presentations $\{X_\alpha\}_{\alpha \in I}$, $\{Y_\alpha\}_{\alpha \in I}$ and $\{Z_\alpha\}_{\alpha \in I}$ for $\scrX$, $\scrY$ and $\scrZ$ respectively.
	At this point, since $\scrC$ has finite limits and $P$ is closed under fiber products, we deduce that
	\begin{align}
		\{X_\alpha \times_{Y_\alpha} Z_\alpha\}_{\alpha \in I} 
	\end{align}
	is a $P$-presentation for $\scrX \times_{\scrY} \scrZ$.
	
	We are therefore left to prove the above claim. Choose $P$-presentations $\{\overline{X}_\alpha\}_{\alpha \in I_1}$, $\{\overline{Y}_\beta\}_{\beta \in I_2}$ and $\{\overline{Z}_\gamma\}_{\gamma \in I_3}$ for $\scrX$, $\scrY$ and $\scrZ$, respectively. Since $X_0$ and $Z_0$ are compact, we can find indexes $\alpha \in I_1$ and $\gamma \in I_3$ together with factorizations of $X_0 \to \scrX$ and of $Z_0 \to \scrZ$ through $\overline{X}_\alpha$ and $\overline{Z}_\gamma$, respectively. Furthermore, since $Y_0$, $\overline{X}_\alpha$ and $\overline{Z}_\gamma$ are compact, we can find $\beta \in I_2$ such that the natural morphisms $\overline{X}_\alpha \longrightarrow \scrY$, $Y_0 \longrightarrow \scrY$ and $\overline{Z}_\gamma \longrightarrow \scrY$ factor through $\overline{Y}_\beta$. Thus, the induced span $\overline{X}_\alpha \to \overline{Y}_\beta \leftarrow \overline{Z}_\gamma$ solves the factorization problem \eqref{eq:factorization_problem_III}.
\end{proof}

\bigskip\section{Relative ind-coherent sheaves}\label{appendix:relative_ind_coherent}

In this section we briefly explain how to adapt Gaitsgory's machinery of indcoherent sheaves \cite{Gaitsgory_IndCoh} to the relative setting. This is important for us to study \textit{families} of nilpotent sheaves over formal schemes. We note that this theory has also been addressed recently, in the setting of solid quasi-coherent sheaves in \cite[\S2.8]{Hansen_Mann}; our treatment is more down-to-earth, and since it plays an important role in the main body of the paper, we decided to keep this appendix for the sake of completeness and readability.

We fix a field of characteristic zero $k$ and $S \in \dAff_k$. Accordingly to the convention of the main body, this implicitly means that $S$ is laft (in the sense of Definition~\ref{def:laft-convergent}).
\begin{definition}
	Let $X \to S$ be a morphism of derived schemes.
	A quasi-coherent sheaf $\calF \in \catQCoh(S)$ is said to be \textit{relatively bounded coherent with respect to $S$} if it is almost perfect on $X$ and has finite tor-amplitude with respect to $S$.
	
	We denote the resulting $\infty$-category by $\catCohb(X/S)$.
\end{definition}

\begin{definition}
	Let $X \to S$ be a laft morphism of qcqs derived schemes, with $S$ affine. We define the $\infty$-category of \textit{indcoherent sheaves on $X$ relative to $S$} to be
	\begin{align}
		\catIndCoh(X/S) \coloneqq \Ind( \catCohb(X/S) ) \ . \tag*{\qedhere}  
	\end{align}
\end{definition}

This notion interpolates between indcoherent and quasicoherent sheaves, as the following example illustrates:
\begin{example}\label{eg:relative_IndCoh_extreme_cases}
	\hfill
	\begin{enumerate}\itemsep=0.2cm
		\item Assume $S = \Spec(K)$, with $K$ field. In this case the relative finite tor-amplitude condition is equivalent to being bounded with respect to the standard $t$-structure on $\catQCoh(X)$. So, in this case $\catCohb(X/S) = \catCohb(X)$ and, accordingly, $\catIndCoh(X/S) = \catIndCoh(X)$.
		
		\item Assume at the other extreme that $X = S$. In this case the finite tor-amplitude condition is equivalent to say that $\catCohb(X) = \catPerf(X)$. Accordingly, $\catIndCoh(X/X) = \catQCoh(X)$. \qedhere
	\end{enumerate}
\end{example}

Most of the results of \cite{Gaitsgory_IndCoh} carry out, without much change, to this setting.
\begin{lemma}\label{lem:relative_IndCoh_base_change}
	Let
	\begin{align}
		\begin{tikzcd}[ampersand replacement=\&]
			Y \arrow{d}{q} \arrow{r}{g} \& X \arrow{d}{p} \\
			T \arrow{r}{f} \& S
		\end{tikzcd}
	\end{align}
	be a pullback square in $\dAff_k$. Then, the functor $g^\ast \colon \catQCoh(X) \longrightarrow \catQCoh(Y)$ restricts to a functor
	\begin{align}
		g^\ast \colon \catCohb(X/S) \longrightarrow \catCohb(Y/T) 
	\end{align}
	and therefore induces a well defined functor in $\PrLomega$
	\begin{align}
		g^\ast \colon \catIndCoh(X/S) \longrightarrow \catIndCoh(Y/T) \ . \tag*{\qedhere} 
	\end{align}
\end{lemma}

\begin{proof}
	Let $\calF \in \catCohb(X/S)$ and $M \in \catQCoh^\heartsuit(T)$. Since $g$ is affine, $g_\ast \colon \catQCoh(Y) \to \catQCoh(X)$ is conservative and $t$-exact. Therefore, $q^\ast(M) \otimes g^\ast(F)$ is in amplitude $[a,b]$ if and only if the same goes for
	\begin{align}
		g_\ast( q^\ast(M) \otimes g^\ast(F) ) \simeq g_\ast(q^\ast(M)) \otimes F \simeq p^\ast f_\ast(M) \otimes F \ . 
	\end{align}
	Since $f_\ast$ is $t$-exact, we have $f_\ast(M) \in \catQCoh^\heartsuit(S)$ and therefore the conclusion follows.
\end{proof}

\begin{notation}
	In the setting of Lemma~\ref{lem:relative_IndCoh_base_change}, we denote by
	\begin{align}
		g_\ast \colon \catIndCoh(Y/T) \longrightarrow \catIndCoh(X/S) 
	\end{align}
	the right adjoint to $g^\ast$. Notice that since $g^\ast$ preserves compact objects, $g_\ast$ is again in $\PrL$.
\end{notation}

\begin{construction}\label{construction:relative_IndCoh_t_structure}
	Let $p \colon X \to S$ be a morphism of derived affine schemes.
	We let
	\begin{align}
		\catCohb(X/S)_{\geqslant 0} \coloneqq \catQCoh(X)_{\geqslant 0} \cap \catCohb(X/S) 
	\end{align}
	and we set
	\begin{align}
		\catIndCoh(X/S)_{\geqslant 0} \coloneqq \Ind(\catCohb(X/S)_{\geqslant 0}) \ . \tag*{\qedhere}  
	\end{align}
\end{construction}

\begin{warning}
	Notice that $\catCohb(X/S)_{\geqslant 0}$ is not the connective part of a $t$-structure on $\catCohb(X/S)$. Indeed, when $X = S$, this recovers $\catPerf(X)_{\geqslant 0}$, which is not the connective part of a $t$-structure unless $X$ is smooth.
\end{warning}

\begin{lemma}\label{lem:relative_IndCoh_t_structure}
	In the setting of Construction~\ref{construction:relative_IndCoh_t_structure}, the category $\catIndCoh(X/S)_{\geqslant 0}$ is a full subcategory of $\catIndCoh(X/S)$ is closed under colimits and extensions. In particular, it defines the connective part of a $t$-structure.
\end{lemma}

\begin{proof}
	The fact that it is a full subcategory follows from the definition and \cite[Proposition~5.3.5.11]{HTT}. Moreover, the definition shows as well that the induced inclusion functor commutes with colimits.
	Let now
	\begin{align}
		\calF' \longrightarrow \calF \longrightarrow \calF''  
	\end{align}
	be a fiber sequence in $\catIndCoh(X/S)$, with $\calF', \calF'' \in \catIndCoh(X/S)_{\geqslant 0}$.
	Rotating the triangle, we have
	\begin{align}
		\calF \simeq \fib( \calF'' \to \calF'[1] ) \ . 
	\end{align}
	Using \cite[Proposition 5.3.5.15]{HTT}, we can choose a simultaneous presentation for this morphism of the form
	\begin{align}
		\colim_{\alpha \in I} \Big( \calF''_\alpha \to \calF'_\alpha[1] \Big) \ , 
	\end{align}
	with $\calF'_\alpha, \calF''_\alpha \in \catCohb(X/S)_{\geqslant 0}$.
	Set
	\begin{align}
		\calF_\alpha \coloneqq \fib( \calF''_\alpha \to \calF'_\alpha[1] ) \ . 
	\end{align}
	Then $\calF_\alpha \in \catCohb(X/S)_{\geqslant 0}$ and
	\begin{align}
		\calF \simeq \colim_{\alpha \in I} \calF_\alpha \ , 
	\end{align}
	which shows that $\calF \in \catIndCoh(X/S)_{\geqslant 0}$.
	
	The last part follows directly from \cite[Proposition~1.4.4.11]{Lurie_HA}.
\end{proof}

\begin{example}
	Continuing Example~\ref{eg:relative_IndCoh_extreme_cases}, we see that the $t$-structure of Lemma~\ref{lem:relative_IndCoh_t_structure} recovers the usual $t$-structure on $\catIndCoh(X)$ when $S = \Spec(K)$, and the one on $\catQCoh(X)$ when $X = S$.
\end{example}

\begin{construction}
	The natural embedding $\catCohb(X/S) \hookrightarrow \catQCoh(X)$ provides a functor
	\begin{align}
		\Psi_{X/S} \colon \catIndCoh(X/S) \longrightarrow \catQCoh(X) = \catIndCoh(X/X) \tag*{\qedhere}  
	\end{align}
\end{construction}

\begin{lemma}\label{lem:relative_IndCoh_eventually_coconnective}
	The functor $\Psi_{X/S}$ restricts for every $n \in \Z$ to an equivalence
	\begin{align}
		\Psi_{X/S} \colon \catIndCoh(X/S)_{\leqslant n} \stackrel{\sim}{\longrightarrow} \catQCoh(X)_{\leqslant n} \ . 
	\end{align}
\end{lemma}

\begin{proof}
	The proof of \cite[Proposition~1.2.4]{Gaitsgory_IndCoh} applies \textit{verbatim}.
\end{proof}

\begin{lemma}\label{lem:relative_IndCoh_relative_Serre_regularity}
	If $p \colon X \to S$ is smooth and $S \in \evccndAff_k$, then the functor $\Psi_{X/S}$ is an equivalence.
\end{lemma}

\begin{proof}
	By definition of $\Psi_{X/S}$ it suffices to show that the full subcategories $\catCohb(X/S)$ and $\catPerf(X)$ of $\catQCoh(X)$ coincide. Since $p$ is smooth, the inclusion $\catPerf(X) \subset \catCohb(X/S)$ is tautological. For the converse implication, let $\calF \in \catCohb(X/S)$. Since by definition $\calF$ is almost perfect, it suffices to show that $\calF$ has finite tor-amplitude. Since $S$ is eventually coconnective, $\calF$ is cohomologically bounded, and it is therefore finitely $m$-presented for some $m \geqslant 0$ by \cite[Remark~2.7.1.4]{Lurie_SAG}. Therefore, the condition of having finite tor-amplitude is open on $X$ by \cite[Corollary 6.1.4.6]{Lurie_SAG}. Thus, it suffices to argue that for every $s \colon T \coloneqq \Spec(K) \to S$ with $K$ a field, the pullback of $\calF$ to the (derived) fiber $X_s$ has finite tor-amplitude. By Lemma~\ref{lem:relative_IndCoh_base_change}, this pullback belongs to $\catCohb(X_s/T)$. Since $X_s$ is smooth, this follows from Serre's regularity criterion.
\end{proof}

\begin{construction}\label{construction:relative_IndCoh_operations}
	Let
	\begin{align}
		\begin{tikzcd}[column sep=small, ampersand replacement=\&]
			X \arrow{rr}{f} \arrow{dr}[swap]{p} \& \& X' \arrow{dl}{q} \\
			\& S
		\end{tikzcd}
	\end{align}
	be a diagram in $\dAff_k$.
	\begin{enumerate}\itemsep=0.2cm
		\item \label{item:construction-relative_IndCoh_operations-1} The functor
		\begin{align}
			f_\ast \colon \catQCoh(X) \longrightarrow \catQCoh(X') 
		\end{align}
		restricts to a functor
		\begin{align}
			f_\ast \colon \catCohb(X/S) \longrightarrow \catQCoh(X'/S)^+ \simeq \catIndCoh(X'/S)^+ \hookrightarrow \catIndCoh(X'/S) \ , 
		\end{align}
		where the middle equivalence is induced by the functor $\Psi_{X'/S}$ via Lemma~\ref{lem:relative_IndCoh_eventually_coconnective}. By left Kan extension, this defines a functor
		\begin{align}
			f_\ast \colon \catIndCoh(X/S) \longrightarrow \catIndCoh(X'/S) \ . 
		\end{align}
		By construction, this functor is in $\PrL$.
		
		\item \label{item:construction-relative_IndCoh_operations-2} If $f$ has finite tor-amplitude, the functor $f^\ast$ restricts to a functor
		\begin{align}
			f^\ast \colon \catCohb(X'/S) \longrightarrow \catCohb(X/S) 
		\end{align}
		and therefore it induces a well defined functor
		\begin{align}
			f^\ast \colon \catIndCoh(X'/S) \longrightarrow \catIndCoh(X/S) \ . 
		\end{align}
		When $f$ is an open immersion, we simply denote this operation by $f^!$. \qedhere
	\end{enumerate}
\end{construction}

\begin{proposition}[Zariski descent]\label{prop:relative_IndCoh_Zariski_descent}
	The $\infty$-category $\catIndCoh(X/S)$ satisfies both Zariski descent in $X$ (with respect to the functoriality described in Construction~\ref{construction:relative_IndCoh_operations}--\eqref{item:construction-relative_IndCoh_operations-2}) and in $S$ (with respect to the functoriality described in Lemma~\ref{lem:relative_IndCoh_base_change}).
\end{proposition}

\begin{proof}
	Notice that the first statement implies the second: indeed, if
	\begin{align}
		\begin{tikzcd}[ampersand replacement=\&]
			V \arrow{r}{v} \arrow{d} \& X \arrow{d} \\
			U \arrow{r}{u} \& S
		\end{tikzcd} 
	\end{align}
	is a pullback square with $u$ a Zariski open immersion, then the functor $v^\ast$ (notation as in Lemma~\ref{lem:relative_IndCoh_base_change}) coincide with the functor $v^!$ (notation as in Construction~\ref{construction:relative_IndCoh_operations}--\eqref{item:construction-relative_IndCoh_operations-2}). As for the first statement, the proof of \cite[Proposition~4.2.1]{Gaitsgory_IndCoh} applies \textit{verbatim}.
\end{proof}

\begin{construction}\label{construction:relative_IndCoh_proper_uppershriek}
	Let $S \in \dAff_k$ and let
	\begin{align}
		\begin{tikzcd}[column sep=small, ampersand replacement=\&]
			X \arrow{rr}{f} \arrow{dr}[swap]{p} \& \& X' \arrow{dl}{q} \\
			\& S
		\end{tikzcd} 
	\end{align}
	be a diagram in $\dSch_{/S}$. If $f$ is proper, the functor $f_\ast \colon \catQCoh(X) \to \catQCoh(X')$ has finite cohomological dimension and preserves almost perfect complexes. Together with the projection formula, this implies that $f_\ast$ restricts to a well defined functor
	\begin{align}
		f_\ast \colon \catCohb(X/S) \longrightarrow \catCohb(X'/S) \ . 
	\end{align}
	In particular,
	\begin{align}
		f_\ast \colon \catIndCoh(X'/S) \longrightarrow \catIndCoh(X/S) 
	\end{align}
	preserves compact objects and therefore its right adjoint
	\begin{align}
		f^! \colon \catIndCoh(X'/S) \longrightarrow \catIndCoh(X/S) 
	\end{align}
	is again in $\PrL$.
\end{construction}

\begin{proposition}[Proper base-change]\label{prop:relative_IndCoh_proper_basechange}
	Fix $S \in \dAff_k$ and let
	\begin{align}
		\begin{tikzcd}[ampersand replacement=\&]
			X' \arrow{r}{p} \arrow{d}{q} \& X \arrow{d}{f} \\
			Y' \arrow{r}{g} \& Y
		\end{tikzcd} 
	\end{align}
	be a pullback square in $\dSch_{/S}$. Assume that $f$ is proper. Then the canonical Beck-Chevalley morphism
	\begin{align}
		p_\ast q^! \longrightarrow f^! g_\ast 
	\end{align}
	is an equivalence.
\end{proposition}

\begin{proof}
	First, one deals with the case where $Y$ is a derived scheme. In this case, leveraging on Lemma~\ref{lem:relative_IndCoh_eventually_coconnective}, the proof of \cite[Proposition~3.4.2]{Gaitsgory_IndCoh} applies \textit{verbatim}. Next, one bootstraps on this case, writing $Y$ as colimit of derived affines over $S$ and using \cite[Corollary~4.7.4.18-(2)]{Lurie_HA}.
\end{proof}

\begin{construction}\label{construction:relative_IndCoh_Kan_extension}
	Fix $S \in \dAff_k$. Using the open-proper Nagata factorization (which can be applied to schemes over $S$ by first forgetting the structural map to the base), we extend the $(-)^!$-functoriality of $\catIndCoh(-/S)$ just as in \cite{Gaitsgory_IndCoh}. The $\infty$-categorical tools needed to perform the extension were developed in \cite{Gaitsgory_Rozenblyum_Study_I}, and their method applies verbatim. We can therefore extend the construction $\catIndCoh(-/S)$ to a functor
	\begin{align}
		\PreSt_{/S}\op \longrightarrow \PrL 
	\end{align}
	via the $(-)^!$ functoriality of Construction~\ref{construction:relative_IndCoh_operations}--\eqref{item:construction-relative_IndCoh_operations-1}. It follows from Proposition~\ref{prop:relative_IndCoh_Zariski_descent} that this construction satisfies Zariski descent and from Proposition~\ref{prop:relative_IndCoh_proper_basechange} that the base-change for the $(-)_\ast \dashv (-)^!$ adjunction holds for proper schematic morphisms.
\end{construction}

\begin{remark}\label{rem:relative_IndCoh_on_formal_schemes}
	In this paper we do not need the full $(-)^!$-functoriality of the previous construction. Indeed, we only need to evaluate $\catIndCoh(-/S)$ on qcqs derived schemes over $S$ and on relative formal schemes. For the latter, we can always find a presentation
	\begin{align}
		\calX \simeq \colim_{W \in \scrT_{\Red{\calX} \overunder \calX}} W \ , 
	\end{align}
	in terms of laft nilpotent thickenings of $\Red{\calX}$, see Recollection~\ref{recollection:presentation_formal_completion}. All the morphisms in $\scrT_{\Red{\calX} \overunder \calX}$ are closed immersions, so the limit
	\begin{align}
		\lim_{W \in \scrT_{\Red{\calX} \overunder \calX}\op} \catIndCoh(W/S) 
	\end{align}
	can equally be computed as
	\begin{align}
		 \colim_{W \in \scrT_{\Red{\calX} \overunder \calX}} \catIndCoh(W/S) 
	\end{align}
	in $\PrL$, with respect to the $(-)_\ast$-functoriality of Construction~\ref{construction:relative_IndCoh_operations}--\eqref{item:construction-relative_IndCoh_operations-1}. In this situation, the proper base-change of Proposition~\ref{prop:relative_IndCoh_proper_basechange} can be proven directly provided that both $f$ and $g$ are schematic and $f$ is proper.
\end{remark}

\begin{proposition}[Localization]\label{prop:relative_IndCoh_localization}
	Let $S \in \dAff_k$ and let $i \colon U \hookrightarrow X$ be an open Zariski embedding in $\dSch_{/S}$. Then the presentable $\infty$-category
	\begin{align}
		\catIndCoh(X/S)_U \coloneqq \ker( i^! \colon \catIndCoh(X/S) \to \catIndCoh(U/S) ) 
	\end{align}
	is compactly generated and
	\begin{align}
		\catIndCoh(X/S)_U^\omega \coloneqq \ker( i^! \colon \catCohb(X/S) \to \catCohb(U/S) ) \ .
	\end{align}
\end{proposition}

\begin{proof}
	When $S = \Spec(K)$ is a field, this is proven in \cite[Proposition~4.1.7]{Gaitsgory_IndCoh}. His argument bootstraps on the analogous property for $\catQCoh$, and using Lemma~\ref{lem:relative_IndCoh_eventually_coconnective}, the same argument applies in the relative setting as well.
\end{proof}

\bigskip\section{Set-theoreticity and pure sheaves}\label{appendix:set-theoreticity-pure}

Let $X$ be an projective scheme over a field $k$ and let $Z$ be a nonempty closed subscheme of $X$. 

The goal of this appendix is the proof of the following result.
\begin{proposition}\label{prop:filtration-scheme-theoretically-appendix}
	Let $\calE$ be a pure coherent sheaf on $X$ of dimension $m\geq 1$, set-theoretically supported on $Z$. Then, there exists a filtration 
	\begin{align}
		0\eqqcolon\calE_{\ell+1}\subset \calE_\ell\subset \calE_{\ell-1} \subset \cdots \subset \calE_0 \coloneqq \calE \ ,
	\end{align}
	for $\ell\geq 1$, so that each subquotient is a pure $m$-dimensional sheaf with scheme-theoretic support contained in $Z$.
\end{proposition}
Before proving the proposition, we need some preliminary results. 

Let $\calI_Z\subset \scrO_X$ be the defining ideal sheaf of $Z$. For any $k\geq 1$, let $Z_{\mathsf{cl}}^{(k)}\subset X$ be the closed subscheme defined by the ideal sheaf $\calI_Z^k \subset \scrO_X$. By convention, we set $\calI_Z^0\coloneqq\scrO_X$. Then, $Z_{\mathsf{cl}}^{(k)}=\trunc{Z^{(k)}}$ for any $k\in \N$, with $k\geq 1$.

Let $\calE$ be a nontrivial coherent sheaf on $X$ of pure dimension $m\geq 1$, with set-theoretic support contained in $Z$. Let $W\subset X$ be the scheme-theoretic support of $\calE$, i.e., the closed subscheme determined by the annihilator ideal $\mathsf{Ann}(\calE)\subset \scrO_X$. Since $\calE$ is set-theoretically supported on $Z$, it follows that $\Red{W}$ is contained in $Z$ as a closed subscheme. 

Note that there exists a unique positive integer $k_\calE\geq 1$ so that $\calI_Z^k\subset \mathsf{Ann}(\calE)$ for all $k\geq k_\calE$, while $\calI_Z^k\nsubseteq \mathsf{Ann}(\calE)$ for all $k < k_\calE$. 

For any $k \geq 1$, let $\calI_Z^k \calE$ be the image of the canonical multiplication map $\calI_Z^k \otimes \calE \to \calE$. Then one has a canonical exact sequence 
\begin{align}\label{eq:Ekseq} 
	0\longrightarrow \calI_Z^k \calE\longrightarrow \calE \longrightarrow \calE\otimes \scrO_{Z_{\mathsf{cl}}^{(k)}}\longrightarrow 0\ .
\end{align}
Furthermore, the subsheaf $\calI_Z^k \calE\subset \calE$ is a nonzero purely $m$-dimensional subsheaf of $\calE$ for any $1\leq k < k_\calE$.
For any $k \geq 1$, let $\calQ_k \subset \calE\otimes \scrO_{Z^{(k)}}$ be the maximal subsheaf so that $\dim \mathsf{Supp}(\calQ_k)<m$. Set $\calF_k \coloneqq \calE\otimes \scrO_{Z_{\mathsf{cl}}^{(k)}}/\calQ_k$. 
\begin{lemma}\label{lem:tensorsupport}
	The following hold:
	\begin{enumerate}\itemsep0.2cm
		\item \label{item:tensorsupport-1} $\calF_k$ is a nonzero pure coherent sheaf of dimension $m$ for all $k\geq 1$. 
		
		\item \label{item:tensorsupport-2} For any $1\leq k < k_\calE-1$ there is an epimorphism $f_{k+1}\colon \calF_{k+1}\to \calF_k$, which fits into a commutative diagram 
		\begin{align}\label{eq:tensordiagA}
			\begin{tikzcd}[ampersand replacement=\&]
				\calE\otimes \scrO_{Z_{\mathsf{cl}}^{(k+1)}} \ar[r] \ar[d] \& \calF_{k+1}  \ar{d}{f_{k}}  \\
				\calE\otimes \scrO_{Z_{\mathsf{cl}}^{(k)}} \ar[r]  \& \calF_{k} 
			\end{tikzcd}\ .
		\end{align} 
		Moreover, $\ker(f_{k})$ is a nonzero subsheaf of $\calF_{k+1}$ annihilated by $\calI_Z$. 
	\end{enumerate} 
\end{lemma} 

\begin{proof} 
	We start by proving \eqref{item:tensorsupport-1}. Clearly, it suffices to prove the claim for $1\leq k \leq k_\calE$. Moreover, given the construction of $\calF_k$, it suffices to prove that the set-theoretic support of $\calE\otimes \scrO_{Z_{\mathsf{cl}}^{(k)}}$ coincides with the set-theoretic support of $\calE$ for any $1\leq k \leq k_\calE$. Assume that this is not the case and let $p$ be a point in $\mathsf{Supp}(\calE)\smallsetminus \mathsf{Supp}(\calE\otimes \scrO_{Z_{\mathsf{cl}}^{(k)}})$ for some $1\leq k \leq k_\calE$. Let $\calJ \subset \scrO_{Z_{\mathsf{cl}}^{(k_{\calE})}}$ be the defining ideal of the closed subscheme $Z_{\mathsf{cl}}^{(k)}\subset Z_{\mathsf{cl}}^{(k_\calE)}$. Note that $\calJ$ is nilpotent. By assumption, the stalk of the tensor product $\calE\otimes  \scrO_{Z^{(k)}}$ at $p$ vanishes. This implies that $\calJ\calE_p = \calE_p$. Since $\calJ$ is nilpotent, Nakayama's lemma  \cite[Tag~07RC, Lemma~10.20.1]{stacks-project} implies that stalk $\calE_p$ vanishes, leading to a contradiction. In conclusion, each quotient $\calF_k$ is nonzero, pure, and of dimension $m$.
	
	Now, we prove \eqref{item:tensorsupport-2}. First, note that claim \eqref{item:tensorsupport-1} implies that the composition 
	\begin{align}
		\calE \otimes \scrO_{Z_{\mathsf{cl}}^{(k+1)}}\longrightarrow \calE\otimes \scrO_{Z^{(k)}}\longrightarrow \calF_k
	\end{align}
	factors through the epimorphism $\calE \otimes \scrO_{Z_{\mathsf{cl}}^{(k+1)}}\to \calF_{k+1}$. Therefore, one obtains indeed a commutative diagram as in \eqref{eq:tensordiagA} for any $k \geq 1$, where the induced morphism $f_k \colon\calF_{k+1}\to \calF_k$ is surjective. 
	
	We next show that $\ker(f_k)\neq 0$ for all $1\leq k < k_\calE$. Assume that this is not the case for some $1\leq k < k_\calE$. Then, $f_k$ is an isomorphism. Moreover, recall that $\calI^k\calE$ is a nonzero pure $m$-dimensional subsheaf of $\calE$. Hence, its set-theoretic support is a non-empty $m$-dimensional closed subspace of $\vert \Red{W}\vert$. Let 
	\begin{align}
		U\coloneqq \mathsf{Supp}(\calI_Z^k\calE) \smallsetminus \mathsf{Supp}(\calQ_k)\ ,
	\end{align}
	and note that $U$ is non-empty since $\dim\mathsf{Supp}(\calQ_i) < m$ for $i=k, k+1$.	Let $p \in U$ be an arbitrary point. By construction, the stalk of $\calF_i$ at $p$ coincides with the stalk of $\calE\otimes \scrO_{Z_{\mathsf{cl}}^{(i)}}$ for $i=k, k+1$. Therefore, the canonical epimorphism 
	\begin{align}
		\calE_p \otimes \scrO_{Z_{\mathsf{cl}}^{(i+1)},p}\longrightarrow \calE_p \otimes \scrO_{Z_{\mathsf{cl}}^{(i)},p}
	\end{align}
	is an isomorphism. Given the canonical commutative diagram
	\begin{align}\label{eq:Ekdiag}
		\begin{tikzcd}[ampersand replacement=\&]
			0\ar[r] \&  \calI_Z^{k+1} \calE\ar[r]\ar[d] \&  \calE \ar[r] \ar{d}{\id}\&  \calE\otimes \scrO_{Z_{\mathsf{cl}}^{(k+1)}}\ar[r]\ar[d] \&  0\\
			0\ar[r] \&  \calI_Z^k \calE\ar[r] \&  \calE \ar[r] \&  \calE\otimes \scrO_{Z_{\mathsf{cl}}^{(k)}}\ar[r] \&  0
		\end{tikzcd}\ ,
	\end{align}
	this implies that the natural injection 
	\begin{align}
		\calI_p^{k+1}\calE_p \longrightarrow \calI_p^k\calE_p
	\end{align}
	is an isomorphism. Then, by Nakayama's lemma \cite[Tag~07RC, Lemma~10.20.1]{stacks-project}, there exists $f \in \calI_p$ so that $(1-f) \calI_p^k \calE_p=0$. Since $p\in \mathsf{Supp}(\calI_Z^k\calE)$, one has $\calI_p^k\calE_p\neq 0$. Hence $(1-f)$ is not a unit in  $\scrO_{X,p}$. This implies that it must belong to the maximal ideal $\frakm_p \subset \scrO_{X,p}$ by \cite[Tag~00E9,Lemma~10.18.2]{stacks-project}. By \textit{loc.cit.}, one also has $\calI_p \subset \frakm_p$, which leads to $1 \in \frakm_p$; clearly a contradiction. In conclusion, $\ker(f_k) \neq 0$ for all $1\leq k < k_\calE$. 
	
	Furthermore, applying the snake lemma to diagram~\eqref{eq:Ekdiag}, one obtains an isomorphism 
	\begin{align}
		\ker(f_k) \longrightarrow \calI_Z^k\calE/\calI_Z^{k+1}\calE \ .
	\end{align}
	This implies that $\ker(f_k)$ is annihilated by $\calI_Z$. 
\end{proof} 

\begin{proof}[Proof of Proposition~\ref{prop:filtration-scheme-theoretically-appendix}] 
	For any $k\geq 1$, let $\calE_k \subset \calE$ be the kernel of the canonical epimorphism $\calE\to \calF_k$. Then, the filtration
	\begin{align}
		0\eqqcolon \calE_{k_\calE} \subset \cdots \subset \calE_0 \coloneqq \calE 
	\end{align}
	is so that each subquotient $\calE_k/\calE_{k-1}$ is a nonzero pure sheaf of dimension $m$, annihilated by $\calI_Z$. Indeed, for any $1\leq k < k_\calE$, diagram~\eqref{eq:tensordiagA} in Lemma~\ref{lem:tensorsupport} yields the commutative diagram 
	\begin{align}
		\begin{tikzcd}[ampersand replacement=\&]
			0\ar[r] \& \calE_{k+1} \ar[r] \ar[d] \& \calE \ar[r] \ar{d}{\id} \& \calF_{k+1} \ar[r] \ar{d}{f_k} \& 0 \\
			0\ar[r] \& \calE_{k}  \ar[r]  \& \calE \ar[r]  \& \calF_k \ar[r]  \& 0
		\end{tikzcd}\ ,
	\end{align}
	where the rows are exact and $\calE_{k+1} \subset \calE_k$ as subsheaves of $\calE$. Moreover, applying the snake lemma, one obtains an isomorphism $\ker(f_k)\simeq \calE_{k}/\calE_{k+1}$. Then, the claim follows from Lemma~\ref{lem:tensorsupport}--\eqref{item:tensorsupport-2}.
\end{proof}



\begin{thebibliography}{DPS{\etalchar{+}}25b}
	
	\bibitem[AHLH23]{Alper-Halpern-Leistner-Heinloth}
	J.~Alper, D.~Halpern-Leistner, and J.~Heinloth, \textit{Existence of moduli
		spaces for algebraic stacks}, Invent. Math. \textbf{234} (2023), no.~3,
	949--1038.
	
	\bibitem[Bar97]{Irred_quot_II}
	V.~Baranovsky, \textit{On punctual {Q}uot schemes for algebraic surfaces},
	\href{https://arxiv.org/pdf/alg-geom/9703038.pdf}{\sf
		arXiv:alg-geom/9703038}, 1997.
				
	\bibitem[BP21]{Binda_Porta}
	F.~Binda and M.~Porta, \textit{Descent problems for derived Azumaya algebras}, \href{https://arxiv.org/pdf/2107.03914.pdf}{\sf arXiv:2107.03914},
	2021.	
	
	\bibitem[BD23]{BD_Okounkov}
	T.~M. Botta and B.~Davison, \textit{Okounkov's conjecture via {BPS} {L}ie
		algebras}, \href{https://arxiv.org/pdf/2312.14008.pdf}{\sf arXiv:2312.14008},
	2023.
	
	\bibitem[BMN25]{Brantner_Magidson_Nuiten}
	L.~Brantner, K.~Magidson, and J.~Nuiten, \textit{Formal integration of derived
		foliations}, \href{https://arxiv.org/pdf/2502.05257.pdf}{\sf
		arXiv:2502.05257}, 2025.
	
	\bibitem[CPT{\etalchar{+}}17]{CPTVV}
	D.~Calaque, T.~Pantev, B.~To\"{e}n, M.~Vaqui\'{e}, and G.~Vezzosi,
	\textit{Shifted {P}oisson structures and deformation quantization}, J. Topol.
	\textbf{10} (2017), no.~2, 483--584.
	
	\bibitem[CW23]{CW-Ind-Geometric}
	S.~Cautis and H.~Williams, \textit{Ind-geometric stacks},
	\href{https://arxiv.org/pdf/2306.03043.pdf}{\sf arXiv:2306.03043}, 2023.
	
	\bibitem[DHSM22]{DHSM22}
	B.~Davison, L.~Hennecart, and S.~Schlegel Mejia, \textit{BPS Lie algebras for totally negative 2-Calabi-Yau categories and nonabelian Hodge theory for stacks}, 
	\href{https://arxiv.org/pdf/2212.07668.pdf}{\sf arXiv:2212.07668}, 2022.
	
	\bibitem[DPS22]{DPS_Torsion-pairs}
	D.-E. Diaconescu, M.~Porta, and F.~Sala, \textit{Cohomological Hall algebras,
		their categorification, and their representations via torsion pairs},
	\href{https://arxiv.org/pdf/2207.08926.pdf}{\sf arXiv:2207.08926}, 2022.
	
	\bibitem[DPS23]{DPS_McKay}
	\bysame, \textit{Mc{K}ay correspondence, cohomological {H}all algebras and
		categorification}, Represent. Theory \textbf{27} (2023), 933--972.
	
	\bibitem[DPS{\etalchar{+}}25]{DPSSV-2}
	D.-E. Diaconescu, M.~Porta, F.~Sala, O.~Schiffmann, and E.~Vasserot,
	\textit{{A}malgamation and a {PBW}-type theorem for nilpotent cohomological
		{H}all algebras}, available at:
	\href{https://web.dm.unipi.it/sala/assets/pdf/COHA_Amalgamation.pdf}{\texttt{link}},
	2025.
	
	\bibitem[DPS{\etalchar{+}}26]{DPSSV-3}
	\bysame, \textit{Cohomological {H}all algebras of one-dimensional sheaves on
		surfaces and {Y}angians}, \href{https://arxiv.org/pdf/2603.03386.pdf}{\sf
		arXiv:2603.03386}, 2026.
	
	\bibitem[DPSZ]{DPSZ}
	D.-E. Diaconescu, M.~Porta, F.~Sala, and Y.~Zhao, \textit{Hecke operators of
		one-dimensional sheaves on a smooth surface}, in preparation.
	
	\bibitem[DEL95]{Donagi_Ein_Lazarsfeld}
	R.~Donagi, L.~Ein, and R.~Lazarsfeld, \textit{A non-linear deformation of the
		{H}itchin dynamical system}, \href{https://arxiv.org/alg-geom/9504017.pdf}{\sf
		arXiv:alg-geom/9504017}, 1995.
	
	\bibitem[DK19]{Dyckerhoff_Kapranov_Higher_Segal}
	T.~Dyckerhoff and M.~Kapranov, \textit{Higher {S}egal spaces}, Lecture Notes in
	Mathematics, vol. 2244, Springer, Cham, 2019.
	
	\bibitem[EL99]{Irred_quot_I}
	G.~Ellingsrud and M.~Lehn, \textit{Irreducibility of the punctual quotient scheme
		of a surface}, Ark. Mat. \textbf{37} (1999), no.~2, 245--254.
	
	\bibitem[FGI{\etalchar{+}}05]{FGA_Explained}
	B.~Fantechi, L.~G\"{o}ttsche, L.~Illusie, S.~L. Kleiman, N.~Nitsure, and
	A.~Vistoli, \textit{Fundamental algebraic geometry}, Mathematical Surveys and
	Monographs, vol. 123, American Mathematical Society, Providence, RI, 2005,
	Grothendieck's FGA explained.
	
	\bibitem[Gai13]{Gaitsgory_IndCoh}
	D.~Gaitsgory, \textit{Ind-coherent sheaves}, Mosc. Math. J. \textbf{13}
	(2013), no.~3, 399--528, 553. 
	
	\bibitem[GR14]{Gaitsgory_Rozenblyum_DG_indschemes}
	D.~Gaitsgory and N.~Rozenblyum, \textit{D{G} indschemes}, Perspectives in
	representation theory, Contemp. Math., vol. 610, Amer. Math. Soc.,
	Providence, RI, 2014, pp.~139--251.
	
	\bibitem[GR17a]{Gaitsgory_Rozenblyum_Study_I}
	\bysame, \textit{A study in derived algebraic geometry. {V}ol. {I}.
		{C}orrespondences and duality}, Mathematical Surveys and Monographs, vol.
	221, American Mathematical Society, Providence, RI, 2017.
	
	\bibitem[GR17]{Gaitsgory_Rozenblyum_Study_II}
	\bysame, \textit{A study in derived algebraic geometry. {V}ol. {II}.
		{D}eformations, {L}ie theory and formal geometry}, Mathematical Surveys and
	Monographs, vol. 221, American Mathematical Society, Providence, RI, 2017.
	
	\bibitem[GKV95]{GKV_Langlands}
	V.~Ginzburg, M.~Kapranov, and E.~Vasserot, \textit{Langlands reciprocity for
		algebraic surfaces}, Math. Res. Lett. \textbf{2} (1995), no.~2, 147--160.
	
	\bibitem[HLP23]{HL_Categorical_properness}
	D.~Halpern-Leistner and A.~Preygel, \textit{Mapping stacks and categorical
		notions of properness}, Compos. Math. \textbf{159} (2023), no.~3, 530--589.
	
	\bibitem[HM26]{Hansen_Mann}
	D.~Hansen and L.~Mann, \textit{The categorical local Langlands conjecture}, 
	\texttt{arXiv:2606.00983}, 2026.
	
	\bibitem[HL10]{HL_Moduli}
	D.~Huybrechts and M.~Lehn, \textit{The geometry of moduli spaces of sheaves},
	second ed., Cambridge Mathematical Library, Cambridge University Press,
	Cambridge, 2010.
	
	\bibitem[KV23]{COHA_surface}
	M.~Kapranov and E.~Vasserot, \textit{The cohomological {H}all algebra of a
		surface and factorization cohomology}, J. Eur. Math. Soc. (JEMS) \textbf{25}
	(2023), no.~11, 4221--4289.
	
	\bibitem[Kha19]{Khan_VFC}
	A.~A. Khan, \textit{Virtual fundamental classes of derived stacks {I}},
	\href{https://arxiv.org/pdf/1909.01332.pdf}{\sf arXiv:1909.01332}, 2019.
	
	\bibitem[Kha21]{Khan_Voevodsky_criterion}
	\bysame, \textit{Voevodsky’s criterion for constructible categories of
		coefficients}, available on the author's webpage:
	\href{https://www.preschema.com/papers/six.pdf}{\sf link}, 2021.
	
	\bibitem[LJZ25]{LJZ_Blown-up}
	W.-P. Li, Q.~Jiang, and Y.~Zhao, \textit{Instantons on the blown-up surface and
		the affine vertex algebra}, \href{https://arxiv.org/abs/2511.18959}{\sf
		arXiv:2511.18959}, 2025.
	
	\bibitem[Lur09]{HTT}
	J.~Lurie, \textit{Higher topos theory}, Annals of Mathematics Studies, vol. 170,
	Princeton University Press, Princeton, NJ, 2009.
	
	\bibitem[Lur17]{Lurie_HA}
	\bysame, \textit{Higher algebra}, available at
	\href{https://www.math.ias.edu/~lurie/papers/HA.pdf}{his webpage}, 2017.
	
	\bibitem[Lur18]{Lurie_SAG}
	\bysame, \textit{Spectral algebraic geometry}, available at
	\href{http://www.math.harvard.edu/~lurie/papers/SAG-rootfile.pdf}{his
		webpage}, 2018.
	
	\bibitem[MO19]{MO_Yangian}
	D.~Maulik and A.~Okounkov, \textit{Quantum groups and quantum cohomology},
	Ast\'{e}risque (2019), no.~408, ix+209.

	\bibitem[MMSV23]{MMSV}
	A.~Mellit, A.~Minets, O.~Schiffmann, and E.~Vasserot, \textit{Coherent sheaves on
		surfaces, {COHA}s and deformed $\mathsf{W}_{1+\infty}$-algebras},
	\href{http://arxiv.org/abs/2311.13415}{\sf arXiv:2311.13415}, 2023.
	
	\bibitem[Min20]{Minets_Higgs}
	A.~Minets, \textit{Cohomological {H}all algebras for {H}iggs torsion sheaves,
		moduli of triples and sheaves on surfaces}, Selecta Math. (N.S.) \textbf{26}
	(2020), no.~2, Paper No. 30, 67.
	
	\bibitem[Nak96]{Nakajima1996}
	H.~Nakajima, \textit{Instantons and affine lie algebras}, Nuclear Phys. B Proc.
	Suppl. \textbf{46} (1996), 154--161.
	
	\bibitem[PS23]{Porta_Sala_Hall}
	M.~Porta and F.~Sala, \textit{Two-dimensional categorified {H}all algebras}, J.
	Eur. Math. Soc. (JEMS) \textbf{25} (2023), no.~3, 1113--1205.
	
	\bibitem[PY16]{Porta_Yu_Higher_Analytic_Stacks}
	M.~Porta and T.~Y. Yu, \textit{Higher analytic stacks and {GAGA} theorems}, Adv.
	Math. \textbf{302} (2016), 351--409.
		
	\bibitem[Rou08]{Rouquier_Dimensions_of_triangulated}
	R.~Rouquier, \textit{Dimensions of triangulated categories}, J. K-Theory
	\textbf{1} (2008), no.~2, 193--256.
		
	\bibitem[Ryd08]{Rydh}	
	D.~Rydh, \textit{Families of cycles and the Chow scheme}, 2008, PhD Thesis, KTH.	
		
	\bibitem[SS20]{Sala_Schiffmann}
	F.~Sala and O.~Schiffmann, \textit{Cohomological {H}all algebra of {H}iggs
		sheaves on a curve}, Algebr. Geom. \textbf{7} (2020), no.~3, 346--376.
		
	\bibitem[SSS25]{SSS}
	F.~Sala, O.~Schiffmann, and P.~Shimpi, \textit{Kleinian orbifolds, Cohomological Hall Algebras, and Yangians}, \href{https://arxiv.org/pdf/2511.08576}{\sf arXiv:2511.08576}, 2025.
	
	\bibitem[SV13]{SV_Cherednik}
	O.~Schiffmann and E.~Vasserot, \textit{Cherednik algebras, {W}-algebras and the
		equivariant cohomology of the moduli space of instantons on
		{$\mathbb{A}^2$}}, Publ. Math. Inst. Hautes \'{E}tudes Sci. \textbf{118}
	(2013), 213--342.
	
	\bibitem[SV20]{SV_generators}
	\bysame, \textit{On cohomological {H}all algebras of quivers: {G}enerators}, J.
	Reine Angew. Math. \textbf{760} (2020), 59--132.
	
	\bibitem[SV23]{SV_YangiansCOHA}
	\bysame, \textit{Cohomological {H}all algebras of quivers and {Y}angians},
	\href{https://arxiv.org/abs/2312.15803}{\sf arXiv:2312.15803}, 2023.
	
	\bibitem[Sim94]{Moduli_pi_one}
	C.~T. Simpson, \textit{Moduli of representations of the fundamental group of a
		smooth projective variety. {I}}, Inst. Hautes \'{E}tudes Sci. Publ. Math.
	(1994), no.~79, 47--129.
	
	\bibitem[{Sta}25]{stacks-project}
	The {Stacks project authors}, \textit{The stacks project}, 2025.
	
	\bibitem[TV08]{TV_HAG-II}
	B.~To\"{e}n and G.~Vezzosi, \textit{Homotopical algebraic geometry. {II}.
		{G}eometric stacks and applications}, Mem. Amer. Math. Soc. \textbf{193}
	(2008), no.~902, x+224.
	
	\bibitem[Zha25]{Zhao_Deformation_normal_cone}
	Y.~Zhao, \textit{Moduli of sheaves and deformation to the normal cone}, \href{https://arxiv.org/pdf/2511.17700}{\sf arXiv:2511.17700} 2025.
	
\end{thebibliography}

\newcommand{\etalchar}[1]{$^{#1}$}
\providecommand{\bysame}{\leavevmode\hbox to3em{\hrulefill}\thinspace}
\providecommand{\MR}{\relax\ifhmode\unskip\space\fi MR }
\providecommand{\MRhref}[2]{%
	\href{http://www.ams.org/mathscinet-getitem?mr=#1}{#2}
}
\providecommand{\href}[2]{#2}

\end{document}